\documentclass[xypic,amssymb,rawfonts,amsthm,amsmath,amsfonts,syntonly,11pt]{amsart}

\usepackage{amssymb}
\usepackage{hyperref}
\usepackage{mathrsfs}
\usepackage[notref,notcite,final]{showkeys}

\newcounter{thecounter}
\numberwithin{thecounter}{section}
\newtheorem{lemma}[thecounter]{Lemma}
\newtheorem{prop}[thecounter]{Proposition}
\newtheorem{thm}[thecounter]{Theorem}
\newtheorem{cor}[thecounter]{Corollary}

\theoremstyle{definition}

\newtheorem{rem}[thecounter]{Remark}
\newtheorem{recol}[thecounter]{Recollection}

\newtheorem{construction}[thecounter]{Construction}

\numberwithin{equation}{section}

\newcommand*\InsertTheoremBreak{
  \begingroup 
    \setlength\itemsep{0pt}
    \setlength\parsep{0pt}
    \item[\vbox{\null}]
  \endgroup
}

\newcommand{\map}{\operatorname{map}}

\newcommand{\colim}{\operatorname{colim}}
\newcommand{\aut}{\operatorname{aut}}

\newcommand{\Aut}{\operatorname{Aut}}
\newcommand{\Rep}{\operatorname{Rep}}

\newcommand{\Out}{\operatorname{Out}}
\newcommand{\End}{\operatorname{End}}

\newcommand{\Hom}{\operatorname{Hom}}
\newcommand{\im}{\operatorname{im}}

\newcommand{\res}{\operatorname{res}}
\newcommand{\Der}{\operatorname{Der}}

\newcommand{\diag}{\operatorname{diag}}

\newcommand{\rank}{\operatorname{rank}}

\newcommand{\hocolim}{\operatorname{hocolim}}

\newcommand{\cC}{{\mathcal{C}}}

\newcommand{\bO}{{\mathbf O}}

\newcommand{\Spaces}{\operatorname{{\bf Spaces}}}

\newcommand{\pcom}{\hat{{}_p}}

\newcommand{\twocom}{\hat{{}_2}}

\newcommand{\GL}{\operatorname{GL}}
\newcommand{\SO}{\operatorname{SO}}
\newcommand{\PU}{\operatorname{PU}}
\newcommand{\DI}{\operatorname{DI}}

\newcommand{\SU}{\operatorname{SU}}
\newcommand{\Sp}{\operatorname{Sp}}
\newcommand{\PSO}{\operatorname{PSO}}
\newcommand{\PSp}{\operatorname{PSp}}

\newcommand{\SL}{\operatorname{SL}}

\newcommand{\Spin}{\operatorname{Spin}}

\newcommand{\Q}{{\mathbb {Q}}}

\newcommand{\F}{{\mathbb {F}}}
\newcommand{\Z}{{\mathbb {Z}}}

\newcommand{\D}{{\mathbf D}}
\newcommand{\CP}{{\mathbb{C}\operatorname{P}}}
\newcommand{\RP}{{\mathbb{R}\operatorname{P}}}

\newcommand{\A}{{\mathbf A}}

\newcommand{\W}{{\mathcal{W}}}
\newcommand{\N}{{\mathcal{N}}}

\newcommand{\C}{{\mathbf{C}}}
\newcommand{\cZ}{{\mathcal{Z}}}

\newcommand{\onto}{\twoheadrightarrow}

\newcommand{\beq}{\begin{eqnarray*}}
\newcommand{\eeq}{\end{eqnarray*}}

\newcommand{\tuborg}{\left\{\begin{array}{ll}}
\newcommand{\sluttuborg}{\end{array}\right.}

\newcommand{\cd}{\operatorname{cd}}

\newcommand{\U}{\mathrm{U}}

\input xypic
\xyoption{all}
\CompileMatrices       

\newcommand{\dN}{\breve \N}
\newcommand{\dZ}{\breve \cZ}
\newcommand{\dT}{\breve T}

\newcommand{\AM}{\Phi}

\renewcommand{\phi}{\varphi}
\renewcommand{\tilde}{\widetilde}

\begin{document}
\title{The classification of $2$-compact groups} 


\author[K. Andersen]{Kasper K. S. Andersen}
\author[J. Grodal]{Jesper Grodal}

\thanks{The second-named author was partially supported by NSF grant DMS-0354633 and an Alfred P. Sloan Research Fellowship}


\address{Department of Mathematical Sciences, University of Aarhus}
\email{kksa@imf.au.dk}

\address{Dept.~of~Mathematics, University of Chicago, and Dept.~of~Math.~Sci., Univ.\ of Copenhagen}
\email{jg@math.uchicago.edu}

\begin{abstract}
We prove that any connected $2$-compact group is classified by its
$2$-adic root datum, and in particular the exotic
$2$-compact group $\DI(4)$, constructed by Dwyer-Wilkerson, is the
only simple $2$-compact group not arising as the $2$-completion of a
compact connected Lie group. Combined with our earlier work with
M{\o}ller and Viruel for $p$ odd, this establishes the full classification
of $p$-compact groups, stating that, up to isomorphism, there is a
one-to-one correspondence between connected $p$-compact groups and root
data over the $p$-adic integers. As a consequence we prove the maximal torus
conjecture, giving a one-to-one correspondence between compact Lie groups
and finite loop spaces admitting a maximal torus. Our proof is a general
induction on the dimension of the group, which works for all
primes. It refines the Andersen-Grodal-M{\o}ller-Viruel methods to
incorporate the theory of root data over the $p$-adic integers, as
developed by Dwyer-Wilkerson and the authors, and we show that certain
occurring obstructions vanish, by relating them to obstruction
groups calculated by Jackowski-McClure-Oliver in the early 1990s.
\end{abstract}

\maketitle



\section{Introduction}

In this paper we prove that any connected $2$-compact group $X$ is
classified, up to isomorphism, by its root datum $\D_X$ over the
$2$-adic integers $\Z_2$. This, combined with our previous work with
M{\o}ller and Viruel \cite{AGMV03} for odd primes, finishes the proof
of the classification of $p$-compact groups for all primes $p$. 
The classification states that, up to isomorphism, there is a
one-to-one correspondence between connected $p$-compact groups and
root data over $\Z_p$. Hence the classification of $p$-compact groups
is completely parallel to the classification of compact connected Lie
groups \cite[\S 4, no.\ 9]{bo9}, just with $\Z$ replaced by $\Z_p$. The
classification of $2$-compact groups has the following consequence,
which captures most of the classification statement. 

\begin{thm}[Classification of $2$-compact groups; splitting version] \label{splitobj}
Let $X$ be a connected $2$-compact group. Then $BX \simeq BG\twocom
\times B\DI(4)^s$, $s\geq 0$, where $G$ is a compact connected Lie
group, and $B\DI(4)$ is the classifying space of the exotic
$2$-compact group $\DI(4)$ constructed by Dwyer-Wilkerson
\cite{DW93}.
\end{thm}

This is the traditional form of the classification conjecture e.g.,
stated by Dwyer in his 1998 ICM address as \cite[Conjs.~5.1 and
5.2]{dwyer98}. A $p$-compact group, introduced by Dwyer-Wilkerson
\cite{DW94}, can be defined as a pointed, connected, $p$-complete
space $BX$ with $H^*(\Omega BX;\F_p)$ finite over $\F_p$, and $X$ is
then the pointed loop space $\Omega BX$. The space $BX$ is hence the
{\em classifying space} of the {\em loop space} $X$, justifying the
convention of referring to the $p$-compact group simply by $X$.
A $p$-compact group is
called {\em connected} if the space $X$ is connected, and two
$p$-compact groups are said to be {\em isomorphic} if their
classifying spaces are homotopy equivalent. For more background on
$p$-compact groups, including details on the history of the
classification conjecture, we refer to \cite{AGMV03} and the
references therein---we also return to it later in this introduction. 

To make the classification more precise,
we now recall the notion of a root datum over $\Z_p$. For $p=2$
this theory provides a key new input to our proofs, and was developed in
the paper \cite{DW05} of Dwyer-Wilkerson and in our paper
\cite{AG05automorphisms}  (see also
Section~\ref{rootdatumappendix}); we say more about this later in the
introduction where we give an outline of the proof of
Theorem~\ref{conn-classification}.

For a principal ideal domain $R$, an {\em $R$-root datum} $\D$ is
a triple $(W,L,\{Rb_\sigma\})$, where $L$ is a finitely generated free
$R$-module, $W \subseteq \Aut_R(L)$ is a finite
subgroup generated by reflections (i.e., elements $\sigma$ such that
$1-\sigma \in \End_{R}(L)$ has rank one), and $\{Rb_\sigma\}$
is a collection of rank one submodules of $L$, indexed by the
reflections $\sigma$ in $W$, satisfying the two conditions
$$
\im(1-\sigma) \subseteq Rb_\sigma \text{ and } w (Rb_\sigma) =
Rb_{w\sigma w^{-1}} \text{ for all } w \in W.
$$
The element $b_\sigma \in L$, determined up to a unit in $R$, is
called the {\em coroot} corresponding to $\sigma$,  and together with
$\sigma$, it determines a {\em root} $\beta_\sigma: L \to R$ via the
formula $\sigma(x) = x + \beta_\sigma(x)b_\sigma$. There is a
one-to-one correspondence between $\Z$-root data and classically
defined root data, by associating $(L,L^*,\{\pm b_\sigma \},\{\pm
\beta_\sigma \})$ to $(W,L,\{\Z b_\sigma\})$; see \cite[Prop.~2.16]{DW05}.
For both $R = \Z$ or $\Z_p$, one can,
instead of $\{Rb_\sigma\}$, equivalently
consider their span, the {\em coroot lattice}, $L_0 =
+_\sigma Rb_\sigma \subseteq L$, the definition given in
\cite[\S 1]{AGMV03} (under the name ``$R$-reflection datum'').
For $R = \Z_p$, $p$ odd, the notion of an
$R$-root datum agrees with that of an $R$-reflection group $(W,L)$; see
Section~\ref{rootdatumappendix}.
Given two $R$-root data $\D = (W,L,\{Rb_\sigma\})$ and $\D' =
(W',L',\{Rb'_{\sigma'}\})$, an {\em isomorphism} between
$\D$ and $\D'$ is an isomorphism $\varphi : L \to L'$ such that $\varphi
W \varphi^{-1} = W'$ as subgroups of $\Aut(L')$ and 
$\varphi(Rb_\sigma) = R b'_{\varphi\sigma\varphi^{-1}}$ for every reflection
$\sigma\in W$. We let $\Aut(\D)$ be the automorphism group of $\D$, and
we define the {\em outer automorphism group} as $\Out(\D) =
\Aut(\D)/W$. A classification of $\Z_p$-root data is given as
Theorems~\ref{rootdata-classification} and
\ref{data-quotient-structure}.

We now explain how to associate a $\Z_p$-root datum to a connected $p$-compact
group. By a theorem of Dwyer-Wilkerson \cite[Thm.~8.13]{DW94} any
$p$-compact group $X$ has a
maximal torus, which is a map $i: BT = (BS^1\pcom)^r \to BX$
satisfying that the homotopy fiber has finite $\F_p$-cohomology and
non-trivial Euler characteristic. Replacing $i$ by an equivalent
fibration, we define the {\em Weyl space} ${\W}_X(T)$ as the
topological monoid of self-maps $BT \to BT$ over $i$. The {\em Weyl
  group} is defined as $W_X(T) = \pi_0({\W}_X(T))$ and the classifying
space of the {\em maximal torus normalizer} is defined as the Borel
construction $B\N_X(T) = BT_{h\W_X(T)}$. By definition, $W_X$
acts on $L_X = \pi_2(BT)$ and if $X$ is connected, this gives a faithful
representation of $W_X$ on $L_X$ as a finite $\Z_p$-reflection group
\cite[Thm.~9.7(ii)]{DW94}.
There is also an easy formula for the coroots $b_\sigma$ in terms of
the maximal torus normalizer $\N_X$, for which we refer to
Section~\ref{rootdatumsub} (or \cite{DW05},
\cite{AG05automorphisms}). We define $\D_X = (W_X, L_X, \{\Z_p
b_\sigma\})$.

We are now ready to state the precise version of our main theorem.

\begin{thm}[Classification of $p$-compact groups] \label{conn-classification}
The assignment which to a connected $p$-compact group $X$ associates
its $\Z_p$-root datum $\D_X$ gives a one-to-one correspondence
between connected $p$-compact groups, up to isomorphism, and
$\Z_p$-root data, up to isomorphism.
Furthermore the map $\AM: \Out(BX)  \to \Out(\D_X)$ is an
isomorphism and more generally
$$
B\Aut(BX) \xrightarrow{\simeq} (B^2\cZ(\D_X))_{h\Out(\D_X)}
$$
where the action of $\Out(\D_X)$ on $B^2\cZ(\D_X)$ is the canonical one.
\end{thm}

Here $\Aut(BX)$ is the space of self-homotopy equivalences of $BX$,
$\Out(BX)$ is its component group, $B^2\cZ(\D_X)$ is the double
classifying space of the center of the $\Z_p$-root datum $\D_X$ (see
Proposition~\ref{centercentralizer}\eqref{centerpart1}) and $\AM$ is
the standard
Adams-Mahmud map \cite[Lem.~4.1]{AGMV03} given by lifting a
self-equivalence $BX$ to $BT$, see Recollection~\ref{adamsmahmud}. We
remark that the existence of the map $B\Aut(BX) \to
(B^2\cZ(\D_X))_{h\Out(\D_X)}$ in the last part of the theorem, requires
the knowledge that the fibration $B\Aut(BX) \to B\Out(BX)$ splits,
which was established in \cite[Thm.~A]{AG05automorphisms}.
Theorem~\ref{conn-classification} implies that connected $p$-compact
groups are classified by their maximal torus normalizer, the
classification conjecture in
\cite[Conj.~5.3]{dwyer98}. For $p$ odd, Theorem~\ref{conn-classification}
is \cite[Thm.~1.1]{AGMV03} (with an improved description of
$B\Aut(BX)$ by \cite[Thm.~A]{AG05automorphisms}). Our proof here is
written so that it is independent of the prime $p$---see the outline
of proof later in this introduction for a further discussion.

\smallskip

The main theorem has a number of important corollaries. The ``maximal
torus conjecture'', gives a purely homotopy theoretic
characterization of compact Lie groups amongst finite loop spaces:

\begin{thm}[Maximal torus conjecture] \label{maxtorus}
The classifying space functor, which to a compact Lie group
$G$ associates the finite loop space $(G,BG,e: G \xrightarrow{\simeq}
\Omega BG)$ gives a one-to-one correspondence between isomorphism classes of
compact Lie groups and finite loop spaces with a maximal torus.
Furthermore, for $G$ connected,  $B\Aut(BG) \simeq (B^2
Z(G))_{h\Out(G)}$.
\end{thm}

The automorphism statement above is included for completeness, but follows 
easily by combining previous work of Jackowski-McClure-Oliver,
Dwyer-Wilkerson, and de Siebenthal (cf., \cite[Cor.~3.7]{JMO95},
\cite[Thm.~1.4]{dw:center}, and \cite[Ch.~I, \S 2, no.\ 2]{desiebenthal56}).
The maximal torus conjecture seems to first have made it into print in
1974, where Wilkerson \cite{wilkerson74} described it as a ``popular
conjecture towards which the author is biased''.

The ``Steenrod problem'' from around 1960 (see Steenrod's papers
\cite{steenrod61,steenrod71}), asks which graded polynomial
algebras are realized as the polynomial ring of some space $X$? The
problem was solved with $\F_p$-coefficients, for $p$ ``large
enough'', by Adams-Wilkerson \cite{AW80} in 1980, extending work of
Clark-Ewing \cite{CE74}, and for all odd $p$, by
Notbohm \cite{notbohm99} in 1999. The case $p=2$ is different from odd primes, 
for instance since generators can appear in odd degrees.

\begin{thm}[Steenrod's problem for $\F_2$] \label{Steenrod}
Suppose that $P^*$ is a graded polynomial algebra over $\F_2$ in
finitely many variables. If $H^*(Y;\F_2)\cong P^*$ for some
space $Y$, then $P^*$ is isomorphic, as a graded
algebra, to
$$
H^*(BG;\F_2) \otimes H^*(B\DI(4);\F_2)^{\otimes r}\otimes
H^*(\RP^\infty;\F_2)^{\otimes s} \otimes H^*(\CP^\infty;\F_2)^{\otimes t}
$$
for some $r, s, t \geq 0$, where $G$ is a compact connected Lie group
with finite center and $\RP^\infty$ and $\CP^\infty$ denotes infinite
dimensional real and complex projective space, respectively.
In particular if $P^*$ has all generators in degree $\geq 3$ then
$P^*$ is a tensor product of the following graded algebras:
\begin{align*}
H^*(B\SU(n);\F_2) & \cong \F_2[x_4, x_6, \ldots, x_{2n}],\\
H^*(B\Sp(n);\F_2) & \cong \F_2[x_4, x_8, \ldots, x_{4n}],\\
H^*(B\Spin(7);\F_2) & \cong \F_2[x_4, x_6, x_7, x_8],\\
H^*(B\Spin(8);\F_2) & \cong \F_2[x_4, x_6, x_7, x_8, y_8],\\
H^*(B\Spin(9);\F_2) & \cong \F_2[x_4, x_6, x_7, x_8, x_{16}],\\
H^*(BG_2;\F_2) & \cong \F_2[x_4, x_6, x_7],\\
H^*(BF_4;\F_2) & \cong \F_2[x_4, x_6, x_7, x_{16}, x_{24}],\\
H^*(B\DI(4);\F_2) & \cong \F_2[x_8, x_{12}, x_{14}, x_{15}].
\end{align*}
\end{thm} 

Since the classification of $p$-compact groups is a space-level
statement, it also gives which graded polynomial algebras over the Steenrod
algebra can occur as the cohomology rings of a space; e.g., the
decomposition in Theorem~\ref{Steenrod} where $P^*$ is assumed to have
generators in degrees $\geq 3$ also holds over the Steenrod algebra.
It should also be possible to give a more concrete list even
without the degree $\geq 3$ assumption, by finding all polynomial
rings which occur as $H^*(BG;\F_2)$ for $G$ a compact connected Lie
group with finite center;  for $G$ simple a list can be found in
\cite[Thm.~5.2]{kono75}, cf.\ Remark~\ref{konoremark}.
In a short companion paper \cite{AG06steenrod} we show
how the theory of $p$-compact groups in fact allows for a solution to
the Steenrod problem with coefficients in an arbitrary
commutative ring $R$.

We can also determine to which extent the realizing space is unique:
Recall that two spaces $Y$
and $Y'$ are said to be $\F_p$-equivalent
if there exists a space $Y''$ and a zig-zag $Y \to Y'' \leftarrow Y'$
inducing isomorphisms on $\F_p$-homology.  The statement below 
also holds verbatim when $p$ is odd, 
where the result is due to Notbohm \cite[Cors.~1.7 and
1.8]{notbohm99}---again complications arise for $p=2$, e.g., due to
the possibility of generators in odd degrees.

\begin{thm}[Uniqueness of spaces with polynomial
  $\F_2$-cohomology] \label{cohomologicaluniqueness} 
If $A^*$ is a graded polynomial
$\F_2$-algebra over the Steenrod algebra ${\mathcal A}_2$, in finitely
many variables, all in degrees $\geq 3$, then there exists at most one
space $Y$, up to $\F_2$-equivalence, with $H^*(Y;\F_2) \cong A^*$, as
graded $\F_2$-algebras over the Steenrod algebra.

If $P^*$ is a finitely generated graded polynomial
$\F_2$-algebra, then there exists at most finitely many spaces $Y$ up
to $\F_2$-equivalence such that $H^*(Y;\F_2) \cong P^*$ as graded
$\F_2$-algebras.
\end{thm}

The early uniqueness results on $p$-compact groups starting with
\cite{DMW86}, which predate root data, or even the formal definition
of a $p$-compact group, were formulated in this language---we give a
list of earlier classification results later in the introduction.
The assumption that all generators are in degrees $\geq 3$ for the
first statement cannot be dropped since for instance $B(S^1 \times
\SU(p^2))$ and $B((S^1 \times \SU(p^2))/C_p)$ have isomorphic
$\F_p$-cohomology algebras over the Steenrod algebra, but are not
$\F_p$-equivalent. Also, the same graded polynomial $\F_p$-algebra can
of course often have multiple Steenrod algebra structures, the option
left open in the second statement: $B\SU(2) \times B\SU(4)$ and
$B\Sp(2) \times B\SU(3)$ have isomorphic $\F_p$-cohomology algebras,
but with different Steenrod algebra structures at all primes.

\smallskip

Bott's theorem on the cohomology of $X/T$, the Peter-Weyl Theorem, and
Borel's characterization of when centralizers of elements of order $p$
are connected, given for $p$ odd as Theorems 1.5, 1.6, and 1.9 of the
paper \cite{AGMV03}, also hold verbatim for $p=2$ as a direct
consequence of the classification (see Remark~\ref{applrem}).
Likewise \cite[Thm.~1.8]{AGMV03}, giving different formulations of being
$p$-torsion free, holds verbatim except that condition (3) should be
removed, cf.\ \cite[Rem.~10.10]{AGMV03}.
We also remark that the classification together with results of
Bott for compact Lie groups, gives that $H^*(\Omega X;\Z_p)$ is
$p$-torsion free and concentrated in even degrees for all $p$-compact
groups. This result was first proved by Lin and Kane, in fact in the
more general setting of finite mod $p$ $H$-spaces, in a series of
celebrated, but highly technical, papers
\cite{lin76,lin78,lin82,kane86}, using completely different
arguments.

\smallskip

Theorem~\ref{conn-classification} also implies a classification for
non-connected $p$-compact groups, though, just as for compact Lie
groups, the classification is less calculationally explicit: Any
disconnected $p$-compact group $X$ fits into a fibration sequence
$$
BX_1 \to BX \to B\pi
$$
with $X_1$ connected, and since our main theorem also includes an
identification of the classifying space of such a fibration
$B\Aut(BX_1)$ with the algebraically defined space
$(B^2\cZ(\D_{X_1}))_{h\Out(\D_{X_1})}$,
this allows for a description of the moduli space of $p$-compact
groups with component group $\pi$ and whose identity component has
$\Z_p$-root datum $\D$. More precisely we have the following
theorem, which in the case where $\pi$ is the trivial group recovers
our classification theorem in the connected case.

\begin{thm}[Classification of non-connected $p$-compact groups]
\label{nonconnclassification}
Let $\D$ be a $\Z_p$-root datum, $\pi$ a finite $p$-group and set
$B\aut(\D) = (B^2\cZ(\D))_{h\Out(\D)}$. The space
$$
M = (\map(B\pi,B\aut(\D)))_{h\Aut(B\pi)}
$$
classifies $p$-compact groups whose identity component has $\Z_p$-root
datum isomorphic to $\D$ and component group isomorphic to $\pi$,
in the following sense:
\begin{enumerate}
\item
There is a one-to-one correspondence between isomorphism classes of
$p$-compact groups $X$ with $\pi_0(X) \cong \pi$ and $\D_{X_1} \cong
\D$, and  components of $M$, given by associating to $X$ the component
of $M$ given by the classifying map $B\pi \to B\Aut(BX_1)
\xrightarrow{\simeq} B\aut(\D_{X_1})$. In particular the set of
isomorphism classes of such $p$-compact groups identifies with
the set of $\Out(\pi)$-orbits on $[B\pi,B\aut(\D)]$, which is finite.
\item
For each $p$-compact group $X$ the corresponding component of $M$ has
the homotopy type of $B\Aut(BX)$ via the zig-zag
$$
B\Aut(BX) \xleftarrow{\simeq}
(\map(B\pi,B\Aut(BX_1))_{C(f)})_{h\Aut(B\pi)} \xrightarrow{\simeq}
(\map(B\pi,B\aut(\D_{X_1}))_{C(f)})_{h\Aut(B\pi)}
$$
where $C(f)$ denote the $\Out(\pi)$-orbit on $[B\pi,B\Aut(BX_1)]$ of
the element classifying $f: BX \to B\pi$.
\end{enumerate}
\end{thm}

\smallskip

Finally, remark that the uniqueness part of the classification
Theorem~\ref{conn-classification} can be reformulated as an {\em
  isomorphism theorem} stating that the isomorphisms, up to
conjugation, between two arbitrary connected $p$-compact groups are
exactly the isomorphisms, up to conjugation, between their root data.
For algebraic groups the {\em isomorphism theorem} can be strengthened
to an {\em isogeny theorem} stating that isogenies of algebraic groups
correspond to isogenies of root data; see e.g., \cite{steinberg99}. In
another companion paper \cite{AG06isogeny} we deduce from our
classification theorem that the same is true for $p$-compact groups:
Homotopy classes of maps $BX \to BX'$ which induce isomorphism in
rational cohomology (the notion of an isogeny for
$p$-compact groups) are in one-to-one
correspondence with the conjugacy classes of isogenies between the
associated $\Z_p$-root data, sending isomorphisms to
isomorphisms. Here an isogeny of $\Z_p$-root data $(W, L, \{\Z_p
b_\sigma\}) \to (W', L', \{\Z_p b'_\sigma\})$ is a $\Z_p$-linear
monomorphism $\phi: L \to L'$ with finite cokernel, such that the
induced isomorphism $\phi: \Aut(L \otimes_{\Z_p} \Q_p) \to
\Aut(L'\otimes_{\Z_p} \Q_p)$ sends $W$ isomorphically to $W'$
and such that $\phi(\Z_p b_\sigma) = \Z_p b'_{\phi(\sigma)}$, for
every reflection $\sigma \in W$ where the corresponding
factor of $W$ has
order divisible by $p$.
As a special case this theorem also  contains the description of
rational self-equivalences of $p$-completed classifying spaces of
compact connected Lie groups,  the most general of the theorems obtained by
Jackowski-McClure-Oliver in \cite{JMO95}, illuminating their result.

\subsubsection*{Structure of the paper and outline of the proof of the classification.}

Our proof of the classification of $2$-compact groups, written to work
for any prime, follows the same overall structure as our proof for $p$
odd with M{\o}ller and Viruel in \cite{AGMV03}, but with significant
additions and modifications. Most importantly we draw on the theory of
root data \cite{DW05} \cite{AG05automorphisms} and have a different way of
dealing with the obstruction group problem. We outline our strategy
below, and also refer the reader to \cite[Sec.~1]{AGMV03} where we
discuss the proof for $p$ odd.

An inspection of the classification of $\Z_p$-root data,
Theorem~\ref{rootdata-classification}, shows that all $\Z_p$-root data
have already been realized as root data of $p$-compact groups by
previous work, so only uniqueness is an issue. (For $p=2$ the root
datum of $\DI(4)$ \cite{DW93} is the only irreducible $\Z_2$-root
datum not coming from a $\Z$-root datum; for $p$ odd see \cite{AGMV03}.)

The proof that two connected $p$-compact groups with isomorphic
$\Z_p$-root data are isomorphic, is divided into a prestep and three
steps, spanning Sections 2--6. Before describing these steps, we have to recall
some necessary results about root data and maximal torus normalizers
from \cite{DW05} and \cite{AG05automorphisms}: The first thing to show is that
the maximal torus normalizer $\N_X$ can be explicitly {\em
  constructed} from $\D_X$. For $p$ odd this follows by a theorem of
the first-named author \cite{kksa:thesis} stating that the maximal
torus normalizer $\N_X$ is always split (i.e., the fibration $BT \to
B\N_X \to B\W_X$ splits). For $p=2$ this is not
necessarily the case, and the situation is more subtle: Recently
Dwyer-Wilkerson \cite{DW05} showed how to extend part of the classical
paper of Tits \cite{tits66} to $p$-compact groups, in particular
reconstructing $\N_X$ from $\D_X$. Since the automorphisms of $\N_X$
differ from those of $X$, one however for classification purposes has
to consider an additional piece of data, namely certain ``root
subgroups'' $\{\N_\sigma\}$, which one can define algebraically for
each reflection $\sigma$.
We constructed these for $p$-compact groups in \cite{AG05automorphisms},
see also Section~\ref{rootdatumappendix}, and described a candidate
algebraic model for the space $B\Aut(BX)$.  (In the
setting of algebraic groups, this ``root subgroup'' $\N_\sigma$ will
be the maximal torus normalizer of $\langle U_\alpha,
U_{-\alpha}\rangle$, where $U_\alpha$ is the root subgroup in the
sense of algebraic groups corresponding to the root $\alpha$ dual to
the coroot $b_\sigma$; see \cite[Rem.~3.1]{AG05automorphisms}.)
For the connoisseur we note that the reliance in \cite{DW05} on a
classification of connected $2$-compact groups of rank $2$ was
eliminated in \cite{AG05automorphisms}.

In \cite{AG05automorphisms}, recalled in
Recollection~\ref{adamsmahmud}, we furthermore showed that the ``Adams-Mahmud'' map,
which to a homotopy equivalence of $BX$ associates a homotopy
equivalence of $B\N$, factors
$$
\AM: \Out(BX) \to \Out(B\N,\{B\N_\sigma\}) \xrightarrow{\cong}
\Out(\D_X)
$$
where $\Out(B\N,\{B\N_\sigma\})$ is the set of homotopy classes of
self-homotopy equivalences of $B\N$ permuting the root subgroups
$B\N_\sigma$ (see Recollection~\ref{adamsmahmud} for the precise
definition).
As a part of our proof we will show that $\AM$ is a isomorphism by induction.

The main argument proceeds by induction on the cohomological dimension
of $X$, which can be determined from $(W_X,L_X)$ alone. (Almost
equivalently one could do induction on the order of $W_X$.) It is
divided into a prestep (Section~\ref{reductionsection}) and three
steps (Sections~4--6).

\smallskip

\noindent{\em Prestep (Section~\ref{reductionsection}):}
The first step is to reduce to the case of simple, center-free
groups. For this we use rather general arguments with fibrations and
their automorphisms, in spirit similar to the arguments in
\cite{AGMV03}. However the theory of root data and root subgroups is
needed both for the statements and results for $p=2$, and these are
incorporated throughout.

With this in hand, we can assume that we have two connected simple,
center-free $p$-compact groups $X$ and $X'$ with isomorphic $\Z_p$-root data
$\D_X$ and $\D_{X'}$. As explained in the discussion above,
\cite{DW05} \cite{AG05automorphisms} implies  that the corresponding maximal
torus normalizers and root subgroups are isomorphic, and we can hence
assume that they both equal $(B\N,\{B\N_\sigma\})$ embedded via maps $j$
and $j'$ in $X$ and $X'$,
$$
\xymatrix{
 & (B\N,\{B\N_\sigma\}) \ar[dr]^-{j'} \ar[dl]_-{j} & \\
                BX \ar@{..>}[rr] & & BX'
}
$$
where the dotted arrow is the one that we want to construct.

\smallskip

\noindent{\em Step 1 (Section~\ref{proofI}):}
Using that in a connected $p$-compact group every element of order $p$ can be
conjugated into the maximal torus, uniquely up to conjugation in $\N$,
we have, for every element $\nu: B\Z/p \to BX$, of order $p$ in $X$, a
diagram of the form
$$
\xymatrix{
 & (B\cC_\N(\nu),\{B{\cC_\N(\nu)}_\sigma\}) \ar[dr] \ar[dl] & \\
                B\cC_X(\nu) & & B\cC_{X'}(\nu).
}
$$
We can furthermore take covers of this diagram with respect to the
fundamental group $\pi_1(\D)$ of the root datum, which we indicate by
adding a tilde $\widetilde{(\cdot)}$. (This uses the formula for the
fundamental group of a $p$-compact group \cite{DW06}, but see also
Theorem~\ref{fundgrp}.)
In Section~\ref{proofI} we prove that one can use the induction hypothesis to
construct a homotopy equivalence between $B\widetilde{\cC_X(\nu)}$ and
$B\widetilde{\cC_{X'}(\nu)}$ under $B\widetilde{\cC_\N(\nu)}$. The
tricky point here is that these centralizers need not themselves be
connected, so one first has to construct the map on the identity
component $B\cC_X(\nu)_1$ and then show that it extends, and this in
turn requires
that one has control of the space of self-equivalences of
$B\cC_X(\nu)_1$.

Now for a general elementary abelian $p$-subgroup $\nu: BE \to BX$ of $X$ we
can pick an element of order $p$ in $E$, and restriction provides a
map
$$
B\widetilde{\cC_X(\nu)} \to B\widetilde{\cC_{X}(\Z/p)} \to
B\widetilde{\cC_{X'}(\Z/p)} \to B\widetilde{X'}.
$$

\smallskip

\noindent{\em Step 2 (Section~\ref{independence}):}
To make sure that these maps are chosen in a compatible way, one has
to show that this map does not depend on the choice of rank one subgroup of $E$. In
Section~\ref{independence} we show that this lift does not depend on
the choices, relying on techniques developed in  \cite{AGMV03}. We
furthermore show that they combine to form an element
$$
[\vartheta] \in \lim{}_{\nu \in
  \A(X)}^0 [B\widetilde{\cC_X(\nu)},B\widetilde{X'}]
$$
where $\A(X)$ is the Quillen category of $X$, with objects
the elementary
abelian $p$-subgroups of $X$ and morphisms induced by conjugation in $X$.

\smallskip 

\noindent{\em Step 3 (Section~\ref{rigidification-section}):}
The construction of the element $[\vartheta]$ basically guarantees that $X$
and $X'$ have the same $p$-fusion, and the last step, which we carry
out in Section~\ref{rigidification-section}, deals with the
rigidification question, where our approach differs significantly
from \cite{AGMV03}. In particular, since we work on universal covers
throughout, we are able to relate our obstruction groups to groups
already calculated by Jackowski-McClure-Oliver in \cite{JMO92}.
Since the exotic $p$-compact groups (only $\DI(4)$ for $p=2$) are
easily dealt with, we can assume that $BX = BG\pcom$ for some simple,
center-free Lie group. It was shown in \cite{JMO92} that
$BG\pcom$ can be expressed as a homotopy colimit of certain
subgroups $P$ of $G$, the so-called $p$-radical (also known as $p$-stubborn)
subgroups. For a $p$-radical subgroup $P$ of $G$, our element
$[\vartheta]$ above gives maps
$$
B\widetilde{P}\pcom \to B\widetilde{C_G({}_pZ(P))}\pcom \to B\widetilde{X'},
$$
where ${}_pZ(P)$ denotes the subgroup of elements of order at most $p$
in $Z(P)$. These maps combine to form an element in
$$
\lim{}_{\widetilde G/\widetilde P \in \bO_p^r(\widetilde
  G)}^0 [B\widetilde P\pcom,B\widetilde X']
$$
where $\bO_p^r(\widetilde G)$ is the full subcategory of the $p$-orbit
category of $\widetilde G$ with objects $\widetilde G/\widetilde P$
for $\widetilde P$ a $p$-radical subgroup of $\widetilde G$.

The obstructions to rigidifying this to get a map on the homotopy colimit
$$
\hocolim_{\widetilde G/\widetilde P \in \bO_p^r(\widetilde G)}
(E\tilde G \times_{\tilde G} \tilde G/\widetilde P)\pcom \to B\widetilde X'
$$
lies in obstruction groups which identify with
$$
\lim{}_{\widetilde G/\widetilde P \in \bO_p^r(\widetilde
  G)}^* \pi_*(Z(\widetilde P)\pcom).
$$
Using extensive case-by-case calculations, Jackowski-McClure-Oliver
showed in \cite{JMO92} that these obstructions in fact vanish.
Hence we have constructed a map 
$$
B\widetilde G\pcom \xleftarrow{\simeq} (\hocolim_{\widetilde
  G/\widetilde P \in \bO_p^r(\widetilde G)} (E\tilde G \times_{\tilde
  G} \tilde G / \widetilde P)\pcom)\pcom \to B\widetilde{X'}
$$
which is easily seen to be an equivalence. Passing to a quotient we get the
wanted equivalence $BG\pcom \to BX'$, finishing the proof of
uniqueness. The remaining statements of
Theorem~\ref{conn-classification} also fall out of this approach.

\smallskip

\noindent Section~\ref{consequences-section} proves the stated
consequences of the classification and the appendix
Section~\ref{rootdatumappendix} is used to establish a number of
general properties of root data over $\Z_p$ used throughout the paper.

Here is the contents in table form:
\tableofcontents

\subsubsection*{Related work and acknowledgments.}

We refer to the introduction of our paper \cite{AGMV03} with M{\o}ller
and Viruel for a detailed discussion of the history of the
classification for odd primes.
The first classification results for $2$-compact groups were obtained
by Dwyer-Miller-Wilkerson \cite{DMW86} twenty years ago, in the
fundamental cases $\SU(2)$ and $\SO(3)$. Notbohm \cite{notbohm94} and
M{\o}ller-Notbohm \cite{MN98} covered $\SU(n)$. Viruel covered $G_2$
\cite{viruel98}, Viruel-Vavpeti{\v{c}} covered $F_4$ and $\Sp(n)$
\cite{VV02} \cite{VVsymplectic}, Morgenroth and Notbohm handled
$\SO(2n+1)$ and $\Spin(2n+1)$ \cite{morgenroththesis}
\cite{notbohm02orthogonal}, and Notbohm proved the result for $\DI(4)$
\cite{notbohm00di4preprint}, all using arguments specific to the case
in question.

Obviously this paper owes a great debt to our earlier work
with M{\o}ller and Viruel \cite{AGMV03} for odd primes. Jesper
M{\o}ller introduced us to the induction-on-centralizers approach to
the classification, and Antonio Viruel gave us the idea of trying to
compare the centralizer and $p$-stubborn decomposition, a method he
had used in the paper \cite{VV02} in a special
case. We are very grateful to them for sharing their insights. We
would furthermore like to thank Bill Dwyer and Clarence Wilkerson for
helpful correspondence, and for sharing an early version of their
manuscript \cite{DW06} on fundamental groups of $p$-compact groups
with us, and Haynes Miller for useful questions. The results of
this paper were announced in Spring 2005 \cite{AGskyetalk}.
Independently of our results, Jesper M{\o}ller has announced a proof
of the classification of connected $2$-compact groups
(Theorem~\ref{splitobj}) using computer algebra \cite{ndet2compactgroups.v2}.
We benefited from the hospitality of Aarhus University, University of
Copenhagen, and University of Chicago while writing this paper.


\section{Reduction to center-free simple $p$-compact groups}
\label{reductionsection}

In this section we reduce the classification of connected $p$-compact
groups to the case of simple center-free groups, in the sense that the
classification statement, Theorem~\ref{conn-classification}, holds for
a connected $p$-compact group if it holds for the simple factors
occurring in the corresponding adjoint (center-free) $p$-compact group
(see Propositions~\ref{productreduction} and
\ref{centerred}).
We do this by extending the proofs given in
\cite[Sec.~6]{AGMV03} for $p$ odd, to all primes by incorporating
$\Z_p$-root data, and root subgroups. Since this additional data
requires restructuring of most of the proofs, we present this
reduction in some detail.

As in \cite{AGMV03} we make the following working definition: A connected
$p$-compact group $X$ is said to be {\em determined by} its $\Z_p$-root datum
$\D_X$ if any connected $p$-compact group $X'$ with $\D_{X'} \cong
\D_X$ is isomorphic to $X$. (Theorem~\ref{conn-classification} will
eventually show that this always holds.)

\begin{prop}[Product Reduction]\label{productreduction}
Suppose $X = X_1 \times \cdots \times X_k$ is a product of simple
$p$-compact groups.
\begin{enumerate}
\item \label{prod1}
If $\AM: \Out(BX_i) \to \Out(\D_{X_i})$ is injective for each $i$,
then so is $\AM:  \Out(BX) \to \Out(\D_{X})$.
\item \label{prod2}
If $\AM: \Out(BX_i) \to \Out(\D_{X_i})$ is surjective and $X_i$ is
determined by $\D_{X_i}$ for each $i$, then $\AM:  \Out(BX) \to
\Out(\D_{X})$ is surjective.
\end{enumerate}
\end{prop}

\begin{proof}
The proof of \eqref{prod1} is identical to the proof of
\cite[Lem.~6.1(2)]{AGMV03} (the key fact is knowing that a map out of
a connected $p$-compact which is trivial when restricted to the
maximal torus is in fact trivial, which e.g., follows from
\cite[Thm.~6.1]{moller96}). The statement in \eqref{prod2} is a direct
consequence of the description of $\Out(\D)$ in
Proposition~\ref{rootdatumauto} together with the assumption that if
$\D_{X_i}$ is isomorphic to $\D_{X_j}$ then $X_i$ is isomorphic to
$X_j$.
\end{proof}

The reader might want to note that conversely to
Proposition~\ref{productreduction}\eqref{prod2}, if $\AM: \Out(BX) \to
\Out(\D_X)$ is surjective for all connected $p$-compact groups $X$,
then all connected $p$-compact groups are determined by their root
datum, as is seen by considering products.

\begin{construction}[Quotients of $p$-compact groups] \label{adjointconstruction}
For explicitness we recall the quotient construction for $p$-compact
groups, and describe when a self-homotopy equivalence induces a
homotopy equivalence on quotients, since this will be used in what
follows:

Let $X$ be a $p$-compact group, $A$ an abelian $p$-compact group and
$i: BA \to BX$ a central homomorphism. By assumption $BA$ is homotopy
equivalent to $\map(BA,BA)_1$ and $\map(BA,BX)_{i}
\xrightarrow{\text{ev}} BX$ is an equivalence. Recall that the
quotient $BX/A$ is defined as the Borel construction of the
composition action of $\map(BA,BA)_1$ on $\map(BA,BX)_{i}$, cf.\
\cite[Pf.~of~Prop.~8.3]{DW94}. This action and the
resulting quotient space $BX/A$ only depends on the (free) homotopy
class of $i$, even on the point-set level, and we have a canonical
quotient map $q: BX \xleftarrow{\simeq} \map(BA,BX)_i \to BX/A$, well
defined up to homotopy.

Now suppose we have a self-homotopy equivalence $f: BX \to BX$ such
that there exists a homotopy equivalence $g: BA \to BA$ making the diagram
$$
\xymatrix{
BA \ar[r]^-g \ar[d]_-i & BA \ar[d]^-i \\
BX \ar[r]^-f & BX
}
$$
commute up to homotopy. We claim that $f$ naturally induces a map on quotients:

First, by using the bar construction model for $BA$, we can without
restriction assume that $g$ is induced by a group homomorphism and has
a strict inverse $g^{-1}$.
Next note that in general, if $\phi: G \to G'$ is a map of monoids,
$h: Y \to Y'$ is a map from a $G$-space $Y$ to a $G'$-space $Y'$, which is
$\phi$-equivariant in the sense that $h(g\cdot y) = \phi(g) \cdot
h(y)$, then there is a canonical induced map on Borel constructions
$Y_{hG} \to Y'_{hG'}$ under $h: Y \to Y'$ and over $B\phi: BG \to
BG'$, e.g., by viewing the Borel construction as a homotopy colimit
via the one-sided bar construction.

In the above setup take  $\phi = c_g$, the monoid automorphism of
$\map(BA,BA)_1$ given by $c_g(\alpha) = g\circ \alpha\circ
g^{-1}$, and $h$ the self-map of $\map(BA,BX)_{i}$ given by
$\beta \mapsto f \circ \beta\circ g^{-1}$.
Then $f$ induces a map $\bar f: BX/A \to BX/A$, which fits into the
homotopy commutative diagram
$$
\xymatrix{
BX \ar[d]_-q \ar[r]^-f & BX \ar[d]^-q \\
BX/A \ar[r]^-{\overline{f}} & BX/A.
}
$$
The quotient construction furthermore behaves naturally with respect
to the maximal torus: If $j:BT \to BX$ is a maximal torus of $X$, and
$h: BT \to BT$ is a lifting of $f:BX \to BX$, then $i:BA \to BX$
factors through $j$ and $g: BA \to BA$ lifts $h$
\cite[Lem.~6.5]{dw:center}, and the above construction produces a
homotopy commutative diagram
$$
\xymatrix{
BT/A \ar[r]^-{\overline{h}} \ar[d]_-{j/A} & BT/A \ar[d]^-{j/A} \\
BX/A \ar[r]^-{\overline{f}} & BX/A
}
$$
up to homotopy under the diagram one has before taking quotients.
\end{construction}

\begin{lemma} \label{autoinjprop}
Let $X$ be a connected $p$-compact group with center $i: B\cZ \to BX$
and let $q: BX \to BX/\cZ$ denote the quotient map. Then any
self-homotopy equivalence $\varphi:
BX\to BX$ fitting into a homotopy commutative diagram
$$
\xymatrix{
 & B\cZ \ar[dr]^-i \ar[dl]_-i & \\
BX \ar[dr]_-q \ar[rr]^-\varphi & & BX \ar[dl]^-q \\
& BX/\cZ
}
$$
is homotopic to the identity.
\end{lemma}

\begin{proof}
Let $q: BX \to BX/\cZ$ denote the quotient map, turned into a
fibration. By changing $\phi$ up to homotopy, we can assume that
$\phi$ is a map strictly over $q$. By \cite{DKS89} (see also
\cite[Prop.~11.9]{dw:center}) the homotopy class of $\phi$ as a map
over $q$ corresponds to an element $[\phi] \in
\pi_1(\map(BX/\cZ,B\Aut(B\cZ))_{f})$, where $f: BX/\cZ \to
B\Aut(B\cZ)$ is the map classifying the fibration $q$. Note that the
class $[\phi]$ could a priori depend on how we choose $\phi$, although this
turns out not to be the case.

The composite $\bar f: BX/\cZ \to B\Aut(B\cZ) \to B\Out(B\cZ)$ is
null-homotopic since $\pi_1(BX/\cZ) = 0$, and obviously
$\map(BX/\cZ,B\Out(B\cZ))_{0} \xrightarrow{ev} B\Out(B\cZ)$ is a
homotopy equivalence. We thus have a fibration sequence
$$
\map(BX/\cZ,B\Aut_1(B\cZ))_{[f]} \to \map(BX/\cZ,B\Aut(B\cZ))_{f} \to
B\Out(B\cZ)
$$
where $[f]$ denotes the components which map to the component of $f$.
Since $B\Aut_1(B\cZ) \simeq B^2\cZ$ is a
loop space, 
$$
\pi_1(\map(BX/\cZ,B\Aut_1(B\cZ))_g) \cong
\pi_1(\map(BX/\cZ,B\Aut_1(B\cZ))_0) = [BX/\cZ,B\cZ] = 0,
$$
for any $g: BX/\cZ \to B\Aut_1(B\cZ)$, using that
$BX/\cZ$ is simply connected and $\pi_2(BX/\cZ)$ is
finite. Hence $\pi_1( \map(BX/\cZ,B\Aut(B\cZ))_f) \to \Out(B\cZ)$ is
injective. But, since $[\phi]$ by assumption maps to the identity in
$\Out(B\cZ)$, we conclude that $[\phi]$ is the identity, and in
particular $\phi$ is homotopic to the identity as wanted.
\end{proof}

\begin{prop}[Reduction to center-free case] \label{centerred} 
Let $X$ be a connected $p$-compact group with center $\cZ$.
\begin{enumerate}
\item \label{center1}
If $X/\cZ$ is determined by $\D_{X/\cZ}$ and $\AM: \Out(BX/\cZ) \to
\Out(\D_{X/\cZ})$ is surjective, then $X$ is determined by $\D_X$.
\item \label{center2}
If $\AM: \Out(BX/\cZ) \to \Out(\D_{X/\cZ})$ is injective then so is
$\AM: \Out(BX) \to \Out(\D_{X})$.
\item \label{center3}
If $\AM:
\Out(BX/\cZ) \to \Out(\D_{X/\cZ})$ is surjective, then so is $\AM:
\Out(BX) \to \Out(\D_{X})$.
\end{enumerate}
\end{prop}

\begin{proof}
The proof of \eqref{center1} follows the outline of the corresponding
statement for odd primes \cite[Lem.~6.8(1)]{AGMV03}, but with the
important additional input that we need to keep track of the root
subgroups: Suppose that $X$ and $X'$ are connected $p$-compact groups
with the same $\Z_p$-root datum $\D$. By \cite[Prop.~1.10]{DW05} and \cite[Thm.~1.2]{kksa:thesis} $X$
and $X'$ have isomorphic maximal torus normalizers $\N_\D$, cf.\
Section~\ref{rootdatumappendix}. By
\cite[Thm.~3.2(2)]{AG05automorphisms} we can choose monomorphisms $j:
B\N_\D \to BX$ and $j': B\N_\D \to BX'$ such that the root subgroups in
$B\N_\D$ with respect to $j$ and $j'$ agree. Furthermore, the centers
of $X$ and $X'$ agree, and can be viewed as a subgroup $\cZ$ of
$\N_\D$. Now, $\N_\D/\cZ$ will be a maximal torus normalizer for both
$X/\cZ$ and $X'/\cZ$ (see e.g.,~\cite[Thm.~1.2]{moller99}) and the
root subgroups in $B\dN_\D/\dZ$ with respect to  $j/\cZ$ and $j'/\cZ$
also agree by construction. By our assumptions there
hence exists a homotopy equivalence $f: BX/\cZ \to BX'/\cZ$ such that
$$
\xymatrix{
 & B\N_\D/\cZ \ar[dl]_-{j/\cZ} \ar[dr]^-{j'/\cZ} \\
BX/\cZ \ar[rr]^-f & & BX'/\cZ
}
$$
commutes up to homotopy. We need to see that $f$ is a map over
$B^2\cZ$, since this implies that $f$ induces a homotopy equivalence
$BX \to BX'$ as wanted. This is a short argument given as the last
part of \cite[Pf.\ of Lem.~6.8(1)]{AGMV03}.

To see \eqref{center2} consider $\varphi\in \Aut(BX)$ corresponding to
an element in the kernel of $\AM: \Out(BX) \to \Out(\D_X)$. By
definition, cf.\ Recollection~\ref{adamsmahmud}, this means that if $i: BT
\to BX$ is a maximal torus then the diagram
$$
\xymatrix{
 & BT \ar[dl]_-i \ar[dr]^-i \\
BX \ar[rr]^-\varphi & & BX
}
$$
commutes up to homotopy. Since the center $B\cZ \to BX$ factors
through $i:BT \to BX$ \cite[Thm.~1.2]{dw:center},
Construction~\ref{adjointconstruction} combined with our assumption shows
that we get a homotopy commutative diagram
$$
\xymatrix{
 & B\cZ \ar[dr] \ar[dl] & \\
BX \ar[dr]_-q \ar[rr]^-\varphi & & BX \ar[dl]^-q \\
& BX/\cZ
}
$$
so Lemma~\ref{autoinjprop} gives the desired conclusion.

We now embark on showing \eqref{center3}, i.e., that $\AM: \Out(BX)
\to \Out(\D_X)$ is surjective, which requires some preparation. Recall
that for any connected $p$-compact group $Y$, $\widetilde{Y}$ is the
$p$-compact group whose classifying space $B\widetilde{Y}$ is the
fiber of the fibration $BY \xrightarrow{q} P_2BY$, where $P_2BY$ is
the second Postnikov section. Let $i:
BT\to BY$ be a maximal torus, which we can assume to be a fibration,
and let $B\widetilde{T}$ denote the fiber of the fibration
$q\circ i:BT \to P_2BY$. Since $\pi_2(i):
\pi_2(BT) \to \pi_2(BY)$ is surjective (see
Proposition~\ref{constrpi1map}), the long exact sequence of homotopy
groups shows that $B\widetilde{T}$ is a $p$-compact torus, and
furthermore $B\widetilde{T} \to B\widetilde{Y}$ is a maximal torus by
the diagram
$$
\xymatrix{
{\widetilde{Y}}/{\widetilde{T}} \ar[r] \ar[d] & {B\widetilde{T}}
\ar[r] \ar[d] & {B\widetilde{Y}} \ar[d] \\
Y/T \ar[r] \ar[d] & BT \ar[r]^-i \ar[d]^-{q\circ i} & BY \ar[d]^-q \\
\ast \ar[r] & P_2BY \ar@{=}[r] & P_2BY.
}
$$
Any self-homotopy equivalence $f: BY\to BY$ lifts to a self-homotopy
equivalence $\widetilde{f}: B\widetilde{Y} \to B\widetilde{Y}$, by
taking fibers, and it is clear that the assignment $f\mapsto
\widetilde{f}$ induces a homomorphism $\Out(BY) \to \Out(B\widetilde{Y})$.

For a $Y$ which satisfies $\pi_1(\D_Y) \xrightarrow{\cong} \pi_1(Y)$,
Proposition~\ref{datum-of-topcover} shows that $\D_{\widetilde{Y}}
\cong \widetilde{\D_Y}$. Hence, the Adams-Mahmud map $\AM$, recalled in
Recollection~\ref{adamsmahmud}, together with
Proposition~\ref{datum-structure} provides the maps in the following
diagram
\begin{equation}\label{AMdiag10}
\xymatrix{
\Out(BY) \ar[d] \ar[r]^-\AM & \Out(\D_Y) \ar[d] \\
\Out(B\widetilde{Y}) \ar[r]^-\AM & \Out(\widetilde{\D_Y}).
}
\end{equation}
The diagram commutes, since for a given $f: BY \to BY$, both
compositions give a map $B\widetilde{T} \to B\widetilde{T}$ over
$\widetilde{f}:B\widetilde{Y} \to B\widetilde{Y}$, and hence they give
the same element in $\Out(\widetilde{\D_Y})$.

By Proposition~\ref{datumofquotient}\eqref{grouppart}, $\widetilde{\D_X}
\xrightarrow{\cong} \widetilde{\D_{X/\cZ}}$, and, chasing through the
definitions, $B\widetilde{X} \xrightarrow{\simeq}
B\widetilde{X/\cZ}$. By the fundamental group formula,
Theorem~\ref{fundgrp},  $\pi_1(\D_{X/\cZ}) \xrightarrow{\cong}
\pi_1(X/\cZ)$. (Note that for the proof of
Theorem~\ref{conn-classification} we can assume that $X/\cZ$ is
determined by $\D_{X/\cZ}$ making this reference to
Theorem~\ref{fundgrp} unnecessary.)
Hence applying diagram
\eqref{AMdiag10} with
$Y=X/\cZ$ and using the aforementioned
identifications for $Y = X/\cZ$ produces the diagram
$$
\xymatrix@C40pt{
\Out(BX/\cZ) \ar[d] \ar@{->>}[r]^(0.47)\AM & \Out(\D_{X/\cZ}) \ar[d]^-\cong \\
\Out({B\widetilde{X}}) \ar[r]^(0.47)\AM & \Out({\widetilde{\D_X}}).
}
$$
Here the right-hand vertical map is an isomorphism by
Proposition~\ref{datumofquotient}\eqref{grouppart} and
Corollary~\ref{datum-adjisocover}, and the top horizontal map is
surjective by assumption. Hence $\AM: \Out({B\widetilde{X}}) \to
\Out(\D_{\widetilde{X}})$ is also surjective.

By \cite[Thm.~5.4]{MN94} there is a short exact sequence $BA \xrightarrow{i}
BX' \to BX$, $BX' = B\widetilde{X} \times B\cZ(X)_1$, where $A$ is a finite
$p$-group and $i: BA \to BX'$ is
central. Proposition~\ref{centercentralizer}\eqref{centerpart2} shows
that $X'$ has $\Z_p$-root datum $\D_{X'} = (W, L_0, \{\Z_p b_\sigma\})
\times (1, L^W, \emptyset)$ and we have the identification $\D_X \cong
\D_{X'}/A$ as in Proposition~\ref{datum-structure}.

We are now ready to show that $\Out(BX) \to \Out(\D_X)$ is surjective,
by lifting an arbitrary element $\alpha \in \Out(\D_X)$ to
$\Out(BX)$. By Proposition~\ref{datum-structure}, $\alpha$ identifies
with an element $\alpha' \in \Out(\D_{X'})$ with $\alpha'(A) =
A$. Since $\AM: \Out({B\widetilde{X}}) \to \Out(\D_{\widetilde{X}})$
is surjective it follows from Proposition~\ref{rootdatumauto} that
$\AM: \Out(BX') \to \Out(\D_{X'})$ is surjective, so we can find a
self-homotopy equivalence $\varphi$ of $BX'$ with $\AM(\varphi) =
\alpha'$. Since $\alpha'(A) = A$ there exists a lift $\varphi': BA \to
BA$ of $\varphi$ fitting into a homotopy commutative diagram
$$
\xymatrix{
BA \ar[d]_-i \ar[r]^-{\varphi'} & BA \ar[d]^-i \\
BX' \ar[r]^-\varphi & BX'.
}
$$
Finally, since $X \cong X'/A$, Construction~\ref{adjointconstruction}
now gives a self-homotopy equivalence $\overline{\varphi}$ of $BX$ with the
property that $\AM(\overline{\varphi}) = \AM(\varphi)/A = \alpha'/A =
\alpha$ as desired.
\end{proof}


\section{Preliminaries on self-equivalences of non-connected $p$-compact groups}

In this short section we prove a fact about
detection of self-equivalences of non-connected $p$-compact groups on
maximal torus normalizers, which we need in the proof of the main
theorem, where non-connected groups occur as centralizers of
elementary abelian $p$-groups in connected groups. 

\begin{prop}\label{disconn}
Let $X$ be a (not necessarily connected) $p$-compact group with maximal
torus normalizer $\N$ and identity component $X_1$. If $\AM: \Out(BX_1)
\xrightarrow{\cong} \Out(\D_{X_1})$ is an isomorphism, then $\AM:
\Out(BX) \to \Out(B\N)$ is injective.
\end{prop}

\begin{proof}
Let $j: B\N \to BX$ be a normalizer inclusion map, which we can assume
is a fibration. Let $f: BX \to
B\pi_0(X) = B\pi$ be the canonical fibration and set $q = f \circ j$. 

We first argue that we can make the identification
$$
\pi_0(\Aut(q)) \xrightarrow{\cong} \{\phi \in \Out(B\N)\,|\,
\phi(\ker(\pi_1(q))) = \ker(\pi_1(q))\}.
$$
Surjectivity is obvious, so we have to see injectivity, where we first
observe that we can pass to a discrete approximation $\breve q:B\dN
\to B\pi$, where $B\dN$ and $B\pi$ are the standard bar construction
models. The simplicial maps $B\dN \to B\dN$ are exactly the group
homomorphisms, so any map $\phi: \breve q \to \breve q$ with $B\dN \to
B\dN$ homotopic to the identity is induced by conjugation by an
element in $\dN$. Hence $\phi$ is homotopic to the identity as a map
of fibrations, proving the claim.

Since $B\Aut(f) \xrightarrow{\simeq} B\Aut(BX)$, we have the following 
diagram, with horizontal maps fibrations, where
$B\N_1$ denotes the fiber of $B\N \to B\pi$:
$$
\xymatrix{
\map(B\pi,B\Aut(BX_1))_{C(f)} \ar[d] \ar[r] & B\Aut(BX) \ar[r] \ar[d]
& B\Aut(B\pi) \ar@{=}[d]\\
\map(B\pi,B\Aut(B\N_1))_{C(q)} \ar[r] & B\Aut(q) \ar[r] & B\Aut(B\pi).
}
$$
Here the horizontal fibrations are established in \cite{DKS89} (see also
\cite[Prop.~11.9]{dw:center}) and the vertical maps are induced by
Adams-Mahmud maps, cf.\ \cite[Lem.~4.1]{AGMV03}, so the diagram is homotopy
commutativity by the naturality of these maps.

To establish the proposition it is enough to verify that $\pi_1(B\Aut(BX)) \to
\pi_1(B\Aut(q))$ is injective since we already saw that 
$\pi_1(B\Aut(q))$ injects into $\Out(B\N)$.
By \cite[Thm.~C]{AG05automorphisms} the map $\AM: B\Aut(BX_1) \to
B\Aut(B\N_1)$ factors through the covering space $Y$ of $B\Aut(B\N_1)$
with respect to the subgroup $\Out(B\N_1,\{(B\N_1)_\sigma\})$, and
$B\Aut(BX_1) \to Y$ has left homotopy inverse. Since $Y \to
B\Aut(B\N_1)$ is a covering, map $\map(B\pi,Y) \to
\map(B\pi,B\Aut(B\N_1))$ is likewise a covering map over each
component where it is surjective, and hence induces a monomorphism on
$\pi_1$ on all components of $\map(B\pi,Y)$. Since $B\Aut(BX_1) \to Y$ has
a homotopy retract,
$\map(B\pi,B\Aut(BX_1)) \to \map(B\pi,B\Aut(B\N_1))$ also induces an injection
on $\pi_1$ for all choices of base-point. 
Hence the five-lemma and the long-exact sequence in homotopy groups
applied to the pair of fibration above guarantees that $\Out(BX) =
\pi_1(B\Aut(BX)) \to \pi_1(B\Aut(q))$ is injective as wanted.
\end{proof}

\begin{prop}\label{disconn2}
Suppose that $X$ is a (not necessarily connected) $p$-compact group
such that $\AM: \Out(BX_1) \xrightarrow{\cong} \Out(\D_{X_1})$. 
Let $i: B\N_p \to BX$ be the inclusion of a
$p$-normalizer of a maximal torus, and let $\phi: BX \to BX$ be a
self-homotopy equivalence.
If $\phi \circ i$ is homotopic to $i$ then $\phi$ is homotopic to the
identity map. 
\end{prop}

\begin{proof}
Let $j: B\N \to BX$ be the inclusion of a maximal torus normalizer,
turned into a fibration. By construction of the Adams-Mahmud map, cf.\
\cite[Lem.~4.1]{AGMV03}, $\phi$ lifts to a map $\phi': B\N \to B\N$, making
the diagram
$$
\xymatrix{
B\N \ar[r]^-{\phi'} \ar[d]^-j & B\N \ar[d]^-j\\
BX \ar[r]^-\phi & BX
}
$$
(strictly) commute, and the space of such lifts is contractible.            

We want to see that $\phi'$ is homotopic to the identity, since
Proposition~\ref{disconn}  then implies that $\phi$ is homotopic to
the identity, as wanted, using the assumption that $\AM: \Out(BX_1)
\xrightarrow{\cong} \Out(\D_{X_1})$.

Lift $i: B\N_p \to BX$ to a map $k: B\N_p \to B\N$.
By
\cite[Pf.~of Lem.~4.1]{AGMV03} the space
of such lifts is contractible, so since $\phi \circ i$ is
homotopic to $i$, we conclude that $\phi' \circ k$ is homotopic to $k$.
Replacing $k$ and $\phi'$ by
discrete approximations, we get the following diagram which commutes
up to homotopy
$$
\xymatrix{
 & B\dN_p \ar[dr]^-{\breve k} \ar[dl]_-{\breve k} \\
B\dN \ar[rr]^-{\breve \phi'} & & B\dN.
}
$$
This is a diagram of $K(\pi,1)$'s, so after changing the spaces and
maps up to homotopy we can assume that all maps are induced by group
homomorphisms, and that the diagram commutes {\em strictly}. But now $\breve
\phi'$ is a group homomorphism which is the identity on $\im(\breve
k)$, and in particular it is the identity on $\dT$. Hence $\breve
\phi'/\dT: W \to W$ is the identity on the Weyl group $W_1$ of $X_1$,
since $W_1$ acts
faithfully on $\dT$. But $W$ is generated by $W_1$ and the image of
$\im(\breve k)$, so we conclude that $\breve \phi'/\dT: W \to W$
is the
identity as well. Hence $\breve \phi'$ is in the image of
$\Der(W;\dT) \to \Aut(\dN)$, cf.\ \cite[Pf.~of~Prop.~5.2]{AGMV03}, and
the homotopy class of $\breve \phi'$ is the image of an element in
$H^1(W;\dT)$ under the homomorphism $H^1(W;\dT) \to \Out(\dN)$. Since
$\breve \phi'$ is the identity on $\im(\breve k)$, this element
restricts trivially to $H^1(W_p;\dT)$, where $W_p$ is a Sylow
$p$-subgroup in $W$. By a transfer argument the element in
$H^1(W;\dT)$ is therefore also trivial, and $\breve \phi'$
is homotopic to the identity map, as wanted.
\end{proof}


\section{First part of the proof of Theorem~\ref{conn-classification}:
  Maps on centralizers} \label{proofI}

In this section we carry out the first part of the proof of
Theorem~\ref{conn-classification}, by constructing maps on certain
centralizers. We have chosen to be quite explicit about when we
replace spaces by homotopy equivalent spaces, since some of these
issues become important later on in the proof, where we want to
conclude that various constructions really take place in certain over-
or undercategories.

Recall that our setup is as follows: Let $X$ and $X'$ be connected
center-free simple $p$-compact groups with isomorphic $\Z_p$-root data
$\D_X$ and $\D_{X'}$. We want to prove that $BX$ is homotopy
equivalent to $BX'$, by induction on cohomological dimension, where by
\cite[Lem.~3.8]{dw:split} the cohomological dimension $\cd Y$ of a
connected $p$-compact group $Y$ depends only on $\D_Y$. We make the
following inductive hypothesis:
\begin{itemize}
\item[$(\star)$]
For all connected $p$-compact groups $Y$ with $\cd Y < \cd X$, $Y$ is
determined by its $\Z_p$-root datum $\D_Y$ and $\AM : \Out(BY) \rightarrow
\Out(\D_Y)$ is an isomorphism.
\end{itemize}

Let $\D$ be a fixed $\Z_p$-root datum, isomorphic to both $\D_X$ and
$\D_{X'}$ and let $\N = \N_\D$
denote the associated normalizer; see
Section~\ref{rootdatumsub}. By \cite[Thm.~3.2(2)]{AG05automorphisms} we
can choose maps $j: B\N \to BX$ and $j': B\N \to BX'$ making $\N$ a
maximal torus normalizer in both $X$ and $X'$ in such a way that the
root subgroups $B\N_\sigma$ of $B\N$ with respect to $j: B\N \to BX$
and $j':B\N \to BX'$ agree.

By Theorem~\ref{fundgrp},
$X$ and $X'$ have canonically isomorphic fundamental groups, which both
identify with $\pi_1(\D)$ via the inclusions $j$ and $j'$.
Applying the second Postnikov section $P_2$ to $BX$ and $BX'$ hence
gives us a diagram
$$
\xymatrix@R10pt{
 & B\N \ar[dr]^-{j'} \ar[dl]_-{j} \\
BX \ar[dr] & & BX'\ar[dl] \\
& B^2\pi_1(\D)
}
$$
where we by changing $B\N$ and the maps up to homotopy can assume that
all maps are fibrations, and that the diagram commutes strictly. After
doing this, we now leave these maps fixed throughout the proof.

Suppose that $\nu:BV \to BX$ is a rank one elementary abelian
$p$-subgroup of $X$, and let $\mu: BV \to BT \to B\N$ denote the
factorization of $\nu$ through the maximal torus $T$, which exists by
\cite[Prop.~5.6]{DW94}, and is unique up to conjugacy in $\N$ by
\cite[Prop.~3.4]{dw:split}. Set $\nu' = j' \circ \mu$ for short.

Taking centralizers, these maps produce the following diagram
\begin{equation}\label{centmaps}
\xymatrix{
 & B\cC_\N(\mu) \ar[dr]^-{j'} \ar[dl]_-{j} \\
B\cC_X(\nu) & & B\cC_{X'}(\nu').
}
\end{equation}
where we by a slight abuse of notation keep the labeling $j$ and
$j'$. Now the fundamental groups of $B\cC_X(\nu)$ and
$B\cC_{X'}(\nu')$ identify via $j$
and $j'$ with a certain quotient group $\pi$ of $\pi_1(B\cC_\N(\mu))$,
explicitly described in \cite[Rem.~2.11]{dw:center} and
Proposition~\ref{centercentralizer}\eqref{centralizerpart}. Passing to
the universal cover of $B\cC_X(\nu)$ and $B\cC_{X'}(\nu')$ and the cover
of $B\cC_\N(\mu)$ determined by the kernel of $\pi_1(B\cC_\N(\mu)) \to
\pi$ produces a diagram
\begin{equation} \label{conndiag0}
\xymatrix{
 & B{\cC_\N(\mu)_1} \ar[dr]^-{j'_{\nu'}} \ar[dl]_-{j_\nu} & \\
B{\cC_X(\nu)_1} & & B{\cC_{X'}(\nu')_1}
}
\end{equation}
where the maps, which are the covers of $j$ and $j'$, are
$\pi$-equivariant with respect to the natural free action of $\pi$ on
all three spaces in the diagram.

Note that in general if $Y$ is a space with a map $f: Y \to BG$, with
$BG$ the classifying space of a simplicial group $G$, a specific model
for the homotopy fiber $\widetilde Y$ of $f$ is given by the subspace
of $Y \times EG$, consisting of pairs whose images in $BG$
agree. In particular it has a canonical free $G$-action,
via the action on the second coordinate and the projection map
$\widetilde{Y} \to Y$ induces a homotopy equivalence $\widetilde{Y}/G
\to Y$. We use this model $\widetilde{( \cdot )}$ for the homotopy
fiber in what follows. Note that if $Y$ has a free $H$-action, then
$\widetilde Y$ has a free $G\times H$-action.

The spaces in \eqref{conndiag0} all have maps to $B^2\pi_1(\D)$ making
the obvious diagrams commute, so we can take homotopy
fibers of these maps by
pulling back along the map $EB\pi_1(\D) \to B^2\pi_1(\D)$, as
described above. This produces the diagram
\begin{equation} \label{conndiag}
\xymatrix{
 & B\widetilde{\cC_\N(\mu)_1} \ar[dr]^-{\widetilde{j'_{\nu'}}}
\ar[dl]_-{\widetilde{j_\nu}} & \\
B\widetilde{\cC_X(\nu)_1} & & B\widetilde{\cC_{X'}(\nu')_1}.
}
\end{equation}
Note that by construction $K = \pi \times B\pi_1(\D)$ acts freely on
the spaces in \eqref{conndiag}, and the maps are $K$-equivariant.
By Propositions~\ref{centercentralizer}\eqref{centralizerpart} and
\ref{datum-of-topcover}, $B\widetilde{\cC_X(\nu)_1}$ and
$B\widetilde{\cC_{X'}(\nu')_1}$ have isomorphic $\Z_p$-root data and
strictly smaller cohomological dimension than $X$, so the inductive
assumption ($\star$) guarantees that they are homotopy equivalent.
Furthermore by construction $\widetilde{j_\nu}$ and $\widetilde{j'_{\nu'}}$
are both maximal torus normalizers and they define the same root
subgroups in $B\widetilde{\cC_\N(\mu)_1}$, since this is true for $j$
and $j'$. Therefore, by the above and the inductive hypothesis ($\star$)
there exists a map $\phi: B\widetilde{\cC_X(\nu)_1} \to
B\widetilde{\cC_{X'}(\nu')_1}$, unique up to homotopy, making the
above diagram \eqref{conndiag} homotopy commute. We now want to argue
that this map can be chosen to be $K$-equivariant so that passing to a
quotient of diagram \eqref{conndiag} with $\phi$ inserted produces a
left-to right map making the diagram \eqref{centmaps} homotopy commute.

Consider the Adams-Mahmud-like zig-zag 
\begin{equation}\label{amlike}
\Xi:
\map(B\widetilde{\cC_X(\nu)_1},B\widetilde{\cC_{X'}(\nu')_1})_{\phi}
\xleftarrow{\simeq} \map(\widetilde{j_\nu},\widetilde{j'_{\nu'}})_{\phi}
\xrightarrow{\simeq}
\map(B\widetilde{\cC_\N(\mu)_1},B\widetilde{\cC_\N(\mu)_1)_1}
\end{equation}
where the first map is a homotopy equivalence by
\cite[Pf.~of~Lem.~4.1]{AGMV03}. That the composite is a homotopy
equivalence will follow once we know that the center of
$\cC_X(\nu)_1$ agrees with the center of $\cC_\N(\mu)_1$, and this
follows from Lemma~\ref{centerformula-maxrank} applied to the subgroup
$\widetilde{\cC_X(\nu)_1}$ of $\widetilde{X}$, where the assumptions
are satisfied since $\pi_1(\D_{\widetilde{X}})=0$ by
Proposition~\ref{datum-of-topcover} and Theorem~\ref{fundgrp}.

By construction the maps in \eqref{amlike} are equivariant with respect to the
$K$-actions. Likewise, since the action on the sources in the mapping
spaces is already free, taking homotopy fixed-points agrees up to
homotopy with taking actual fixed-points, so the maps
in \eqref{amlike} induce homotopy equivalences between the
fixed-points. This produces homotopy equivalences 
$$
\map_K(B\widetilde{\cC_X(\nu)_1},B\widetilde{\cC_{X'}(\nu')_1})_{[\phi]}
\xleftarrow{\simeq} \map_K(\widetilde{j_\nu} ,\widetilde{j'_{\nu'}})_{[\phi]}
\xrightarrow{\simeq}
\map_K(B\widetilde{\cC_\N(\mu)_1},B\widetilde{\cC_\N(\mu)_1})_{[\text{1}]},
$$
where the subscript $[\phi]$ denotes that we are taking all components
of maps non-equivariantly homotopy equivalent to $\phi$.
We can therefore pick an {\em equivariant} map $\psi \in
\map_{K}(B\widetilde{\cC_X(\nu)_1},B\widetilde{\cC_{X'}(\nu')_1})_{[\phi]}$
corresponding to $1 \in
\map_K(B\widetilde{\cC_\N(\mu)_1},B\widetilde{\cC_\N(\mu)_1})_{[\text{1}]}$.
Define $\widetilde{h_\nu}$ as the composite
$$
\widetilde{h_\nu}: B\widetilde{\cC_X(\nu)} \xleftarrow{\simeq}
(B\widetilde{\cC_X(\nu)_1})/\pi \xrightarrow[\simeq]{\psi/\pi}
(B\widetilde{\cC_{X'}(\nu')_1})/\pi \xrightarrow{\simeq}
B\widetilde{\cC_{X'}(\nu')}
$$
and similarly define $h_\nu: B\cC_X(\nu) \xleftarrow{\simeq} (B\widetilde{\cC_X(\nu)_1})/K
\xrightarrow[\simeq]{\psi/K}  (B\widetilde{\cC_{X'}(\nu')_1})/K
\xrightarrow{\simeq} B\cC_{X'}(\nu').$

By construction the maps $\widetilde{h_\nu}$ and $h_\nu$ fit into
following homotopy commutative diagram
$$
\xymatrix{
 & B\widetilde{\cC_\N(\mu)} \ar[dr]^-{\widetilde{j'_{\nu'}}}
\ar[dl]_-{\widetilde{j_\nu}} \ar[d]\\
B\widetilde{\cC_X(\nu)} \ar[d] \ar@/_3pc/[rr]_\simeq^{\widetilde
  h_\nu} & B\cC_\N(\mu) \ar[dr]^(0.4){j'_{\nu'}} \ar[dl]_(0.4){j_\nu} &
B\widetilde{\cC_{X'}(\nu')} \ar[d]\\
B{\cC_X(\nu)} \ar@/_3pc/[rr]_\simeq^{h_\nu} & & B{\cC_{X'}(\nu')}\\
}
$$
and are in fact uniquely determined up to homotopy  by this property
by Proposition~\ref{disconn} and the inductive assumption $(\star)$.

Let $\widetilde{\phi_\nu} : B\widetilde{\cC_X(\nu)} \to
B\widetilde{X'}$ be the composite of $\widetilde{h_\nu}$ with the
evaluation $B\widetilde{\cC_{X'}(\nu')} \to B\widetilde{X'}$ and
similarly define $\phi_\nu: B\cC_X(\nu) \to BX'$ as the composite of
$h_\nu$ with the evaluation $B\cC_{X'}(\nu') \to BX'$.

We define $\phi_\nu$ and $\widetilde{\phi_\nu}$ when $\nu: BE \to BX$
is an elementary abelian $p$-subgroup of rank
greater that one, by restricting to a rank one subgroup $V\subseteq E$
and using
adjointness as follows:  As
before $\nu|_V$ factors through $T$, uniquely as a map to $\N$, and we
let $\mu: BV \to B\N$ denote the resulting map to $\N$.
Restriction produces a map
$$
\phi_{\nu,V}: B\cC_X(\nu) \to B\cC_X(\nu|_V) \xrightarrow[\simeq]{h_{\nu|_V}}
B\cC_{X'}({j'\mu}) \to BX'.
$$
Similarly, letting $B\widetilde{\cC_X(\nu)}$ denote the homotopy fiber of the
map $B\cC_X(\nu) \to BX \to B^2\pi_1(\D)$ as in the rank one case, we
define $\widetilde{\phi_{\nu,V}}$ as
$$
\widetilde{\phi_{\nu,V}}: B\widetilde{\cC_X(\nu)} \to
B\widetilde{\cC_X(\nu|_{V})} \xrightarrow[\simeq]{\widetilde{h_{\nu|_V}}}
B\widetilde{\cC_{X'}(j'\mu)} \to B\widetilde{X'}.
$$
By construction $\phi_{\nu,V}$ and $\widetilde{\phi_{\nu,V}}$ fit
together in that the diagram
$$
\xymatrix{
B\widetilde{\cC_X(\nu)} \ar[r]^-{\widetilde{\phi_{\nu,V}}} \ar[d] &
B\widetilde{X'} \ar[d]\\
B{\cC_X(\nu)} \ar[r]^-{\phi_{\nu,V}} & B{X'}
}
$$
commutes up to homotopy, and is in fact a homotopy pull-back square by
construction.
This concludes the construction of the maps
on centralizers which we will use in the next sections to construct
our equivalence $BX \to BX'$. We will in particular prove that
$\phi_{\nu,V}$ and $\widetilde{\phi_{\nu,V}}$ are independent of
the choice of the rank one subgroup $V\subseteq E$, after which we
will drop the subscript $V$.


\section{Second part of the proof of
  Theorem~\ref{conn-classification}: The element in $\lim^0$}
\label{independence}

In this section we prove that the maps $\widetilde{\phi_{\nu,V}}$
constructed in the previous section are independent of the
choice of rank one subgroup $V\subseteq E$ and give coherent maps into
$B\widetilde{X'}$. More specifically we prove the following.

\begin{thm}\label{collectlemma1}
Let $X$ and $X'$ be two connected simple center-free $p$-compact
groups with isomorphic $\Z_p$-root data, and assume the inductive
hypothesis $(\star)$.
Then the maps
$$
\widetilde{\phi_{\nu,V}}: B\widetilde{\cC_X(\nu)} \to B\widetilde{X'}
$$
constructed in Section~\ref{proofI} are independent of the
choice of $V$ and together form an element in
$\lim^0_{\nu \in \A(X)} [B\widetilde{\cC_X(\nu)},B\widetilde{X'}]$.
\end{thm}

Here $\A(X)$ is the Quillen category of $X$ with objects the elementary
abelian $p$-subgroups $\nu: BE \to BX$ of $X$ and morphisms given by
conjugation (i.e., the morphisms from $(\nu: BE\to BX)$ to $(\nu': BE'
\to BX)$ are the linear maps $\phi:E \to E'$ such
that $\nu$ is freely homotopic to $\nu' \circ B\phi$).

We need the following proposition, whose proof we postpone to after
the rest of the proof of Theorem~\ref{collectlemma1}.

\begin{prop}\label{ranktwo}
Let $X$ be a connected simple center-free $p$-compact group. If $\nu: BE
\to BX$ is a non-toral elementary abelian $p$-subgroup of rank two,
then $\cC_X(\nu)_1$ is non-trivial or $\D_X \cong \D_{\PU(p)\pcom}$.
\end{prop}

\begin{proof}[Proof of Theorem~\ref{collectlemma1}]
We divide the proof into two steps.  Step 1 verifies the independence
of the choice of $V$, and the shorter Step 2 then uses this to
construct the element in $\lim^0$.

\smallskip

\noindent{\em Step 1: The maps $\widetilde{\phi_{\nu,V}}$ and
  $\phi_{\nu,V}$ are independent of the choice of rank one subgroup
  $V$:} We divide this step into three substeps a--c. Step 1a assumes
$\nu$ toral, Step 1b assumes $\nu$ rank two non-toral, and finally
Step 1c considers the general case.

\noindent{\em {Step 1a:} Assume $\nu: BE \to BX$ is a toral
elementary abelian $p$-subgroup.}
By assumption $\nu$ factors through $BT$ to give a map $\mu: BE \to
B\N$, unique up to conjugation in $\N$, and as in the rank one  case
we let $\nu' = j' \mu$.
We want to say that the map $\widetilde{\phi_{\nu,V}}$ does not depend
on $V$, basically since it is a map suitably under
$B\widetilde{\cC_\N(\mu)}$, and hence uniquely determined, independently
of $V$. This will follow by adjointness, analogously to
\cite[Pf.~of~Thm.~2.2]{AGMV03}, although a bit of care has to be
taken, since we have to verify that this happens over $B^2\pi_1(\D)$
in order to be able to pass to the cover $\widetilde{(\cdot)}$, as we
now explain.

Recall that by construction the map $h_{\nu|_V}$ is the bottom
left-to-right composite in the following diagram
\begin{equation}
\xymatrix@C8pt{
 & & B\cC_\N(\mu|_V) \ar[ddll]_-{j_{\nu|_V}} \ar[ddrr]^-{j'_{\nu'|_V}}
 & & \\
 & & (B\widetilde{\cC_\N(\mu|_V)_1})/K \ar[u]_-\simeq \ar[dl] \ar[dr] \\
B\cC_X(\nu|_V) & **[r] (B\widetilde{\cC_X(\nu|_V)_1})/K \ar[l]_-\simeq
\ar[rr]_-{\simeq}^-{\psi/K} & & **[l]
(B\widetilde{\cC_{X'}({\nu'|_V})_1})/K \ar[r]^-\simeq &
B\cC_{X'}(\nu'|_V)
}
\end{equation}
where we notice that all subdiagrams commute up to homotopy {\em over $BK$}.

Since $BE$ maps into these spaces via $\nu$ and $\mu$,
adjointness produces the following diagram
\begin{equation}\label{cooldiag1}
\xymatrix{
 & & B\cC_\N(\mu) \ar[ddll]_-{j_\nu} \ar[ddrr]^-{j'_{\nu'}} & & \\
 & & B\cC_{\cC_\N(\mu|_{V})}(\mu) \ar[u]_-\simeq \ar[dl] \ar[dr] \\
B\cC_X(\nu) & B\cC_{\cC_X(\nu|_{V})}(\nu) \ar[l]_(0.55)\simeq &
Z \ar[l]_-\simeq \ar[r]_-{\simeq}^-{\cC(\psi/K)} &
B\cC_{\cC_{X'}(\nu'|_{V})}({\nu'}) \ar[r]^(0.55)\simeq & B\cC_{X'}(\nu')
}
\end{equation}
where $Z = \map(BE,(B\widetilde{\cC_X(\nu|_V)_1})/K)_\nu$
and $\cC(\psi/K)$ is the map induced by $\psi/K$ on mapping spaces.

Since the diagram \eqref{cooldiag1} homotopy commutes {\em as a
diagram over $B^2\pi_1(\D)$} we get a homotopy commutative diagram by
passing to homotopy fibers:
\begin{equation}\label{cooldiag2}
\xymatrix{
 & & B\widetilde{\cC_\N(\mu)} \ar[ddll]_-{\widetilde{j_\nu}}
\ar[ddrr]^-{\widetilde{j'_{\nu'}}} & & \\
 & & B\widetilde{\cC_{\cC_\N(\mu|_{V})}(\mu)} \ar[u]_-\simeq \ar[dl] \ar[dr] \\
B\widetilde{\cC_X(\nu)} & B\widetilde{\cC_{\cC_X(\nu|_{V})}(\nu)}
\ar[l]_(0.55)\simeq &
{\widetilde{Z}} \ar[l]_-\simeq \ar[r]_-{\simeq}^-{\widetilde{\cC(\psi/K)}} &
B\widetilde{\cC_{\cC_{X'}(\nu'|_{V})}({\nu'})} \ar[r]^(0.55)\simeq &
B\widetilde{\cC_{X'}(\nu')}.
}
\end{equation}
Denote the bottom left-to-right homotopy equivalence in
\eqref{cooldiag2} by $\cC_{\nu}(\widetilde{h_{\nu|_V}})$, justified by
the fact that by construction the following diagram homotopy commutes
\begin{equation}\label{cooldiag3}
\xymatrix{
B\widetilde{\cC_X(\nu)} \ar[r]^-{\cC_{\nu}(\widetilde{h_{\nu|_V}})}_-\simeq
\ar[d] & B\widetilde{\cC_{X'}(\nu')} \ar[d] \ar[dr] \\
B\widetilde{\cC_X(\nu|_{V})}
\ar[r]^-{\widetilde{h_{\nu|_{V}}}}_-\simeq &
B\widetilde{\cC_{X'}(\nu'|_{V})} \ar[r] & B\widetilde{X'}.
}
\end{equation}
Diagram
\eqref{cooldiag2}, together with the inductive assumption $(\star)$
and Proposition~\ref{disconn} shows that the homotopy class of
$\cC_{\nu}(\widetilde{h_{\nu|_V}})$ does not depend on $V$. Hence by diagram
\eqref{cooldiag3} the same is true for $\widetilde{\phi_{\nu,V}}$
which is what we wanted.
(The key point in the above argument is that we can
choose $\mu$ once and for all, such that $\mu|_{V}$ is a factorization
of $\nu|_V$ through $BT$ for every $V \subseteq E$.)

Note that the {\em construction} of $\cC_{\nu}(\widetilde{h_{\nu|_V}})$
in \eqref{cooldiag2} does {\em not} depend on $\nu$ being toral, as
long as $\nu'$ in that case is defined as $BE \xrightarrow{\nu}
B\cC_X(\nu) \xrightarrow{\phi_{\nu,V}} BX'$ (instead of $j'\mu$),
which makes sense in this more general setting---the top part
of diagram \eqref{cooldiag2} is only needed to conclude the
independence of $V$. In Step 1b below we will also use the notation
$\cC_{\nu}(\widetilde{h_{\nu|_V}})$ for non-toral $\nu$. 

\smallskip
\noindent{\em {Step 1b:} Assume $\nu: BE \to BX$ is a rank two non-toral
elementary abelian $p$-subgroup.} By Proposition~\ref{ranktwo} either 
$\cC_X(\nu)_1$ is non-trivial or  $\D_X \cong \D_{\PU(p)\pcom}$.

Assume first that $\D_X \cong \D_{\PU(p)\pcom}$. Since uniqueness for
this group is well known, both for $p$ odd and $p=2$ by \cite{DMW86},
\cite{moller97puppreprint}, \cite{BV99} and \cite{AGMV03}, the
statement of course follows for this reason. But one can also argue
directly, using a slight modification of the proof of
\cite[Lem.~3.2]{AGMV03} which we quickly sketch: For $\alpha \in
W_X(\nu)$ we have the following diagram
$$
\xymatrix@C=40pt{
B\widetilde{\cC_X(\nu)} \ar[r] \ar[d] & B\widetilde{\cC_X(\nu|_V)}
\ar[d] \ar[r]^-{\widetilde{h_{\nu|_V}}}_-{\simeq} &
B\widetilde{\cC_{X'}(j'\mu)} \ar[d] \ar[r] & B\widetilde{X'}
\ar@{=}[d] \\
B\widetilde{\cC_X(\nu)} \ar[r] & B\widetilde{\cC_X(\nu|_{\alpha(V)})}
\ar[r]^-{\widetilde{h_{\nu|_{\alpha(V)}}}}_-{\simeq} &
B\widetilde{\cC_{X'}(j'\mu\circ(\alpha^{-1}|_{\alpha(V)}))} \ar[r] &
B\widetilde{X'}. \\
}
$$
Here all the non-identity vertical maps are given on the level of
mapping spaces (i.e., without the tilde) by $f \mapsto f \circ
\alpha^{-1}$, which induces a map on the indicated spaces by taking
homotopy fibers of the map to $B^2\pi_1(\D)$.
The left-hand and right-hand squares obviously commute and the middle
square commutes up to homotopy by our inductive assumption ($\star$)
and Proposition~\ref{disconn}. We thus conclude that
$\widetilde{\varphi_{\nu,\alpha(V)}} \circ
\widetilde{B\cC_X(\alpha^{-1})} \simeq \widetilde{\varphi_{\nu,V}}$
for all $\alpha\in W_X(\nu)$.
Now, since $\D  \cong \D_{\PU(p)\pcom}$ we have $\pi_1(\D) \cong
\Z/p$, and because $E$ is non-toral, we see that $B\cC_X(\nu) \simeq
BE$ and $B\widetilde{\cC_X(\nu)} \simeq BP$, where $P$ is the
extra-special group $p^{1+2}_+$ of order $p^3$ and exponent $p$ if $p$
is odd and $Q_8$ if $p=2$. Now \cite[Prop.~3.1]{AGMV03} shows that
$W_X(\nu)$ contains $\SL(E)$, and the same is true for $X'$. The proof
of \cite[Prop.~3.1]{AGMV03} furthermore shows that
$\widetilde{\varphi_{\nu,V}} \circ \widetilde{B\cC_X(\alpha)} \simeq
\widetilde{\varphi_{\nu,V}}$ for $\alpha \in \SL(E)$. (Apply the
argument there to the $p$-group $P$ instead of $E$.) Since $\SL(E)$
acts transitively on the rank one subgroups of $E$, combining the
above gives that $\widetilde{\varphi_{\nu,V'}} \simeq
\widetilde{\varphi_{\nu,V}}$ for any rank one subgroup $V'\subseteq E$ as
desired.

We can therefore assume that
$\cC_X(\nu)_1$ is non-trivial, and the proof in this case is an
adaptation of \cite[Pf.~of~Lem.~3.3]{AGMV03} to our new setting:
Choose a rank one elementary abelian $p$-subgroup
$\eta: BU = B\Z/p \to B\cC_X(\nu)_1$ in the center of a $p$-normalizer of
a maximal torus in $\cC_X(\nu)$. Let $\eta \times \nu: BU \times BE \to
BX$ be the map defined by adjointness, and for any rank one subgroup
$V$ of $E$, consider the map $\eta \times \nu|_{V} : BU \times BV \to BX$
obtained by restriction. By construction $\eta \times \nu|_{V}$ is the
adjoint of the composite $BU \xrightarrow{\eta} B\cC_X(\nu)_1
\xrightarrow{\res} B\cC_X(\nu|_{V})_1$, so $\eta \times \nu|_{V} : BU
\times BV \to BX$ factors through a maximal torus in $X$ by
\cite[Prop.~5.6]{DW94}. It is furthermore straightforward to check that $\eta
\times \nu$ is a monomorphism (compare \cite[Pf.~of~Lem.~3.3]{AGMV03}).

Now consider the following diagram
$$
\xymatrix{
 & B\cC_X(\nu|_V) \ar[rd]^-{\varphi_{\nu|_V}} \\
BU \times BE \ar[r] \ar[ru] \ar[rd] & B\cC_X(\eta\times \nu|_V) \ar[u]
\ar[d] & BX' \\
 & B\cC_X(\eta) \ar[ru]_-{\varphi_\eta}
}
$$
Here the left-hand side of the diagram is constructed by taking adjoints of
$\eta \times \nu$ and hence it commutes.
The right-hand side
homotopy commutes by Step 1a (using the inductive assumption
$(\star)$),  since $\eta \times \nu|_{V}$ is toral of rank two.
We can hence without ambiguity define $(\eta \times \nu)'$ as either the top
left-to-right composite (for some rank one subgroup $V \subseteq E$)
or the bottom left-to-right composite. Denote by $\nu'$ the
restriction of $(\eta \times \nu)'$ to $BE$.

By construction of the map $\cC_{\eta \times
  \nu}(\widetilde{h_{\nu|_V}})$ as the bottom composite in
\eqref{cooldiag2}, the diagram
\begin{equation}\label{tauto}
\xymatrix@C60pt{
B\widetilde{\cC_X(\eta \times \nu)}
\ar[r]^-{\cC_{\eta \times \nu}(\widetilde{h_{\nu|_V}})}_-\simeq \ar[d] &
B\widetilde{\cC_{X'}((\eta \times \nu)')} \ar[d] \\
B\widetilde{\cC_X(\nu)} \ar[r]^-{\cC_{\nu}(\widetilde{h_{\nu|_V}})}_-\simeq
& B\widetilde{\cC_{X'}(\nu')}
}
\end{equation}
commutes. Furthermore, since $\eta \times \nu|_{V}$ is toral, diagram
\eqref{cooldiag2}, applied with $\mu$ equal to a factorization of $\eta \times
\nu|_{V}$ through $BT$, shows that the top horizontal map in
\eqref{tauto} agrees with $\cC_{\eta \times
  \nu}(\widetilde{h_{\eta}})$, and in particular it is independent of
$V$, again using Proposition~\ref{disconn} and our inductive
assumption $(\star)$.

We claim that this forces the same to be true for the bottom
horizontal map in \eqref{tauto}:
By our choice of $\eta$, the centralizer $\cC_X(\eta \times \nu)$ contains a
$p$-normalizer of a maximal torus in $\cC_X(\nu)$, and hence
$\widetilde{\cC_X(\eta \times \nu)}$ contains a
$p$-normalizer of a maximal torus in $\widetilde{\cC_X(\nu)}$.
Proposition~\ref{disconn2} and our inductive
assumption $(\star)$ therefore shows that the bottom
map in \eqref{tauto} is independent of $V$, so
$\widetilde{\varphi_{\nu,V}}: B\widetilde{\cC_X(\nu)}
\xrightarrow{\cC_{\nu}(\widetilde{{h_{\nu|_V}}})}
B\widetilde{\cC_{X'}(\nu')} \to B\widetilde{X'}$ is also
independent of $V$ as wanted.

\smallskip
\noindent{\em {Step 1c:} Assume $\nu: BE \to BX$ is
an elementary abelian $p$-subgroup of rank $\geq 3$.}
The fact that $\widetilde{\varphi_{\nu,V}}$
is independent of $V$ when $E$ has
rank two implies the statement in general: Let $\nu : BE \to BX$
be an elementary abelian $p$-subgroup of rank at least three, and
suppose that $V_1\neq V_2$ are two rank
one subgroups of $E$. Setting $U = V_1 \oplus V_2$ we get the following diagram
$$
\xymatrix{
 & B\widetilde{\cC_X(\nu|_{V_1}) }
 \ar[dr]^-{\widetilde{\varphi_{\nu|_{V_1}}}} & \\
B\widetilde{\cC_X(\nu)} \ar[ur] \ar[r] \ar[dr] &
B\widetilde{\cC_X(\nu|_{U})} \ar[u] \ar[d] & B\widetilde{X'} \\
 & B\widetilde{\cC_X(\nu|_{V_2})} \ar[ur]_-{\widetilde{\varphi_{\nu|_{V_2}}}}
}
$$
The left-hand side of this diagram is constructed by adjointness and
hence commutes, and the right-hand side of the diagram commutes up to
homotopy by Steps 1a and 1b. Thus the top left-to-right composite
$\widetilde{\varphi_{\nu,V_1}}$ is homotopic to the bottom
left-to-right composite $\widetilde{\varphi_{\nu,V_2}}$, i.e., the map
$\widetilde{\varphi_{\nu,V}}$ is independent of the choice of rank one subgroup
$V$ as claimed.

\smallskip
\noindent{\em { Step 2:} An element in $\lim^0$.}
With the above preparations in place it is easy to see, as in
\cite[Pf.~of~Thm.~2.2]{AGMV03}, that the maps $\widetilde{\phi_\nu} :
B\widetilde{\cC_X(\nu)} \to B\widetilde{X'}$ fit together to form an
element in
$$
\lim{}_{\nu \in \A(X)}^0
[B\widetilde{\cC_X(\nu)},B\widetilde{X'}].
$$
In order not to cheat the reader of the finale, we repeat the short
argument from \cite[Pf.~of~Thm.~2.2]{AGMV03}: For
any morphism  $\rho: (\nu: BE \to BX) \to (\xi: BF \to BX)$
in $\A(X)$ we need to verify that
$$
\xymatrix{
B\widetilde{\cC_X(\xi)} \ar[rr]^-{B\widetilde{\cC_X(\rho)}}
\ar[rd]_-{\widetilde{\varphi_{\xi}}} &
& B\widetilde{\cC_X(\nu)} \ar[ld]^-{\widetilde{\varphi_{\nu}}} \\
& B\widetilde{X'}
}
$$
commutes up to homotopy.
If $F$ has rank one, then $\rho$ is an isomorphism, 
and we let $\mu: BF \to BT \to B\N$
be a factorization of $\xi$ through $BT$. In this case the claim follows since
\begin{equation}\label{steptworankone}
\xymatrix@C40pt{
B\widetilde{\cC_X(\xi)} \ar[r]_-{\simeq}^-{\widetilde{h_{\xi}}}
\ar[d]^-{\simeq}_-{B\widetilde{\cC_X(\rho)}} &
B\widetilde{\cC_{X'}(j'\mu)}
\ar[d]^-{\simeq}_-{B\widetilde{\cC_{X'}(\rho)}} \ar[dr] & \\
B\widetilde{\cC_X(\xi \rho)}
\ar[r]^-{\widetilde{h_{\xi \rho}}}_-\simeq  &
B\widetilde{\cC_{X'}(j'\mu\rho)} \ar[r] & B\widetilde{X'}
}
\end{equation}
commutes up to homotopy, since we can view the diagram
as taking place under $B\widetilde{\cC_\N(\mu)} \xrightarrow{\simeq}
B\widetilde{\cC_\N(\mu\rho)}$, up to homotopy, using the inductive
assumption ($\star$) and Proposition~\ref{disconn} as in Step 1a.
The case where $F$ has arbitrary rank follows from the independence of
the choice of rank one subgroup $V$, established in Step 1, together
with the rank one case: For a rank one subgroup $V$ of $E$ set $V' =
\rho(V)$ and consider the diagram
$$
\xymatrix@C40pt{
B\widetilde{\cC_X(\xi)} \ar[r] \ar[d]_-{B\widetilde{\cC_X(\rho)}} &
B\widetilde{\cC_X({\xi}|_{V'})} \ar[d]_-{B\widetilde{\cC_X(\rho|_V)}} \ar[dr]
& \\
B\widetilde{\cC_X(\nu)} \ar[r] & B\widetilde{\cC_X(\nu|_{V})} \ar[r]
& B\widetilde{X'}
}
$$
The left-hand side commutes by construction and the right-hand side
commutes since the diagram \eqref{steptworankone} commutes, proving the
claim. 

This constructs an element $[\vartheta] \in \lim_{\nu \in
  \A(X)}^0[B\widetilde{\cC_X(\nu)},B\widetilde{X'}]$ as
wanted. 
\end{proof}

We now give the proof of Proposition~\ref{ranktwo}, used in the proof
of Theorem~\ref{collectlemma1}, which we postponed. The proof uses
case-by-case arguments on the level of $\Z_p$-root data.

\begin{proof}[Proof of Proposition~\ref{ranktwo}]
Assume first that $(W_X, L_X)$ is an exotic $\Z_p$-reflection group. We
claim that the centralizer of any rank one elementary abelian
$p$-subgroup of $X$ is connected, which in particular implies that there
can be no rank two non-toral elementary abelian $p$-subgroups:
For $p>2$ this follows from \cite[Thms.~11.1,
12.2(2) and 7.1]{AGMV03} combined with Dwyer-Wilkerson's formula for
the Weyl group of a centralizer \cite[Thm.~7.6]{dw:center}. For $p=2$
we have $\D_X \cong \D_{\DI(4)}$ and hence $W_X \cong \Z/2 \times
\GL_3(\F_2)$ where the central $\Z/2$ factor acts by $-1$ on $L_X$ and
$\GL_3(\F_2)$ acts via the natural representation on $L_X\otimes \F_2$,
cf.\ \cite[Pf~of~Thm.~11.1]{AGMV03} or \cite[Rem.~7.2]{DW05}. In
particular it follows
directly (cf.\ \cite[Pf.~of~Prop.~9.12, $\DI(4)$ case]{DW05}) that $X$
contains a unique elementary abelian $2$-subgroup of rank one up to
conjugation and that the centralizer of this subgroup is
connected.

By the classification of $\Z_p$-root data,
Theorem~\ref{rootdata-classification}\eqref{rdc-part1}, we
may thus assume that $\D_X$ is of the form $\D_{G\pcom}$ for some
simple compact connected Lie group $G$, and since $X$ is center-free
we may assume that so is $G$.
If $\pi_1(G)$ has no $p$-torsion,
\cite[Thm.~2.27]{steinberg75} implies that the centralizer in $G$ of
any element of order $p$ is connected. By the formula for
the Weyl group of a centralizer \cite[Thm.~7.6]{dw:center}
(cf.\ Proposition~\ref{centercentralizer}\eqref{centralizerpart}), it then
follows that $\cC_X(\eta)$ is
connected for any rank one elementary abelian $p$-subgroup $\eta:
B\Z/p \to BX$ of $X$, and hence $X$ does not have any rank two
non-toral elementary abelian $p$-subgroups.

We can thus furthermore assume that $\pi_1(G)$ has $p$-torsion, which implies
that $\D_X \cong \D_{G\pcom}$ for one of the following $G$: $G=\PU(n)$
for $p|n$; $G=\SO(2n+1)$, $n\geq 2$ for $p=2$; $G=\PSp(n)$, $n\geq 3$
for $p=2$; $G=\PSO(2n)$, $n\geq 4$ for $p=2$; $G=PE_6$ for $p=3$ and
$G=PE_7$ for $p=2$.  We want to see that if $\eta: B\Z/p \to BX$ has
rank one and $\D_X \not\cong \D_{\PU(p)\pcom}$, then $\cC_X(\eta)$ is
not a $p$-compact toral group, since then \cite[Cor.~1.4]{DM97}
implies that for any elementary abelian $p$-subgroup $\nu$ of rank
two, $\cC_X(\nu)_1$ is non-trivial.

By \cite[Thm.~7.6]{dw:center} it is enough to see this for the
corresponding Lie group $G$. So, let $V\subseteq G$ be a rank one
elementary abelian $p$-subgroup of $G$. If $C_G(V)_1$ is a torus then
$W_{C_G(V)_1}=1$ and hence $W_{C_G(V)} \cong \pi_0(C_G(V))$. By
\cite[\S 5, Ex.~3(b)]{bo9} or \cite[Thm.~9.1(a)]{steinberg68}
$\pi_0(C_G(V))$ is isomorphic to a subgroup of $\pi_1(G)$, so
$|W_{C_G(V)}|$ divides $|\pi_1(G)|$. In particular, if $x$
is a generator of $V$ then the number of elements in a maximal torus
which are conjugate to $x$ in $G$ is at least $|W_G|/|\pi_1(G)|$
(since two elements in a maximal torus are conjugate in $G$ if and
only if they are conjugate by a Weyl group element).
In particular
$$
p^n-1 \geq \frac{|W_G|}{|\pi_1(G)|},
$$
where $n$ is the rank of $G$. A direct case-by-case check of the above
cases shows that this
inequality can only hold when $G=\PU(n)$, $p|n$. In this case, a
generator of $V$ must have the form $x=[\diag(\lambda_1, \ldots,
\lambda_n)]$. For $n>p$, some of the $\lambda$'s must agree, so
$C_G(V)_1$ is not a torus in this case. This proves the claim.
\end{proof}


\section{Third and final part of the proof of
  Theorem~\ref{conn-classification}: Rigidification} \label{rigidification-section}

In this section we finish the proof of
Theorem~\ref{conn-classification} by showing that our element in
$\lim^0$ from Theorem~\ref{collectlemma1} rigidifies to produce a
homotopy equivalence $BX \to BX'$. We first need a lemma:

\begin{lemma} \label{centricity}
Suppose that we have a homotopy pull-back square of $p$-compact groups
$$
\xymatrix{
BX' \ar[d]^-{g'} \ar[r]^-{f'} & BY' \ar[d]^-g \\
BX \ar[r]^-f & BY
}
$$
where $f:BX \to BY$ is a centric monomorphism (i.e., $\map(BX,BX)_1
\xrightarrow{\simeq} \map(BX,BY)_f$ and $Y/X$ is $\F_p$-finite)
and $g:BY' \to BY$ is an epimorphism (i.e., $Y/Y'$ is the classifying
space of a $p$-compact group).
Then $f': BX' \to BY'$ is a centric monomorphism, and $g'$ is an epimorphism.
\end{lemma}

\begin{proof}
It is clear from the pull-back square that $f'$ is a monomorphism and
that $g'$ is an epimorphism. To see that $f'$ is centric, observe that
we have a map of fibrations
$$
\xymatrix{
(Y/Y')^{hX'} \ar[r] \ar@{=}[d] & \map(BX',BX')_{[g']}
\ar[r] \ar[d] & \map(BX',BX)_{g'} \ar[d] \\
(Y/Y')^{hX'} \ar[r] & \map(BX',BY')_{[f \circ g']} \ar[r] &
\map(BX',BY)_{f \circ g'}
}
$$
where e.g., the subscript $[g']$ denotes the components of the
mapping space mapping to the component of $g'$. Since the wanted map
$\map(BX',BX')_1 \to \map(BX',BY')_{f'}$ is the restriction of the map of total
spaces to a component, it is hence enough to see that the map between
the base spaces is an equivalence, which follow since we have equivalences
$$
\xymatrix{
\map(BX,BX)_1 \ar[r]^-\simeq \ar[d]^-\simeq & \map(BX',BX)_{g'} \ar[d]\\
\map(BX,BY)_{f} \ar[r]^-\simeq & \map(BX',BY)_{f \circ g'}.
}
$$
Here the vertical equivalence follows from the centricity of $f$ and
the horizontal equivalences follows from \cite[Prop.~3.5]{dwyer96cent}
combined with \cite[Prop.~10.1]{dw:center} (in
\cite[Prop.~3.5]{dwyer96cent} taking $E = BX'$, $B = BX$, and $X =
BX$ and $BY$ respectively).
\end{proof}

\begin{proof}[Proof of Theorem~\ref{conn-classification}]
By the results of Section~\ref{reductionsection}
(Propositions~\ref{productreduction} and \ref{centerred}) we are reduced
to the case where we consider two simple, center-free
$p$-compact groups $BX$ and $BX'$ with the same root datum
$\D$. Furthermore, in Theorem~\ref{collectlemma1} of the previous
section we constructed an element $[\vartheta] \in \lim^0_{\nu \in
  \A(X)} [B\widetilde{\cC_X(\nu)},B\widetilde{X'}]$.
We want to use this element to construct the map from $BX$ to $BX'$,
and show that it is an equivalence. By the classification of $\Z_p$-root data,
Theorem~\ref{rootdata-classification}\eqref{rdc-part1}, $\D$ is either
exotic or $\D \cong \D_{G\pcom}$ for a compact connected Lie group $G$ and we
handle these two cases separately.

Suppose that $\D$ is exotic, and notice that in this case $\pi_1(\D) =
0$ by Theorem~\ref{rootdata-classification}\eqref{rdc-part2}.
If $p$ is odd we are exactly in the situation covered by rather easy
arguments in \cite{AGMV03} (see \cite[Pf.~of~Thm.~2.2]{AGMV03} and
\cite[Prop.~9.5]{AGMV03}).
For $p=2$ uniqueness of $B\DI(4)$, as well
as the statements about self-maps, are already well known
by the work of Notbohm \cite{notbohm00di4preprint}, but we
nevertheless quickly remark how a proof also falls out of the current
setup, noticing that the arguments from
\cite{AGMV03} from this point on, in the special case of $\DI(4)$,
carries over verbatim: We have $\D \cong \D_{\DI(4)}$ and since
$\pi_1(\D) =0$, the functor
$\nu \mapsto \pi_i(\map(B\widetilde{\cC_X(\nu)},B\widetilde{X'})_\nu)$
identifies with $\nu \mapsto \pi_i(B\cZ(\cC_X(\nu)))$, as explained in
detail in \cite[Pf.~of~Thm.~2.2]{AGMV03}. Since we are considering $X = \DI(4)$
we know by Dwyer-Wilkerson \cite[Prop.~8.1]{DW93} that
$\lim^j_{\A(X)}\pi_i(B\cZ(\cC_X(\nu))) = 0$ for all $i,j \geq 0$. (The
proof of this is a Mackey functor argument, and relies on the
regular structure of the Quillen category of $\DI(4)$, due to the fact
that its classifying space, like the exotic $p$-compact groups for $p$
odd, has polynomial $\F_p$-cohomology ring.)
The centralizer decomposition
theorem \cite[Thm.~8.1]{dw:center} now produces a map $BX \to BX'$, which by
standard arguments given in \cite[Pf.~of~Thm.~2.2]{AGMV03} is seen to be an
equivalence. The statement about self-maps also follows as in
\cite[Pf.~of~Thm.~2.2]{AGMV03}, using $\Out(B\N,\{B\N_\sigma\})$
instead of $\Out(B\N)$.

Now, suppose that $\D \cong \D_{G\pcom}$, for a simple center-free
compact Lie group $G$, with universal cover $\widetilde{G}$. Let
$\bO_p^r(\widetilde{G})$ be the full subcategory of the orbit category
with objects the $\widetilde{G}$-sets $\widetilde{G}/\widetilde{P}$ with
$\widetilde{P}$ a $p$-radical subgroup (i.e, $\widetilde{P}$ is an
extension of a torus by a finite $p$-group, such that
$N_{\widetilde{G}}(\widetilde{P})/\widetilde{P}$ is finite and
contains no non-trivial normal $p$-subgroups).

For the $p$-radical homology decomposition \cite[Thm.~4]{JMO92} of
the compact Lie group $\tilde G$, one considers the functor $F:
\bO_p^r(\widetilde G) \to \Spaces$ given by $\widetilde G/\widetilde P
\mapsto E{\tilde G} \times_{\tilde G} {\tilde G}/{\widetilde
  P}$, where $\Spaces$ denotes the category of topological spaces.
Viewed as a functor to the {\em homotopy category} of spaces, $\text{Ho}(\Spaces)$,
this functor is isomorphic to the functor $F': \bO_p^r(\widetilde G)
\to \text{Ho}(\Spaces)$ given on objects by $\tilde G/\tilde P \mapsto
B\widetilde P$ and on morphisms by sending the ${\tilde G}$-set map
$\widetilde{G}/\widetilde{P} \xrightarrow{f} \widetilde{G}/\widetilde{Q}$
to the map $c_{g^{-1}}:B{\tilde P} \to B{\tilde Q}$, where
$g\widetilde{Q} = f(e\widetilde{P})$, via the canonical equivalences $B\tilde P =
(E{\tilde P})/{\tilde P} \xrightarrow{\simeq} E{\tilde G}
\times_{\tilde G} {\tilde G}/{\tilde P}$. (Note that $F'$ is
not well-defined as a functor to $\Spaces$,
since the element $g$ is just an arbitrary coset representative for
the coset $g{\tilde Q}$.) We can hence in what follows replace
$E{\tilde G} \times_{\tilde G} {\tilde G}/\tilde P$ by $B{\tilde P}$
in this way, whenever we are working in the homotopy category.

Since $\widetilde P$ and $\widetilde Q$ are $p$-radical, the same is
true for their images $P$ and $Q$ in $G$, and $C_G(P) = Z(P)$ and
likewise for $Q$ (see \cite[Prop.~1.6(i) and Lem.~1.5(ii)]{JMO92}).
Hence there is a well defined
induced morphism $c_g: {}_pZ(Q) \to {}_pZ(P)$ as well as a well-defined
(free) homotopy class of maps $c_{g^{-1}}: B\widetilde{C_G({}_pZ(P))}
\to B\tilde{C_G({}_pZ(Q))}$. Consider the diagram
$$
\xymatrix@C60pt{
B{\tilde P}\pcom \ar[r] \ar[d]^-{c_{g^{-1}}} &
(B\widetilde{C_G({}_pZ(P))})\pcom \ar[r]^-\simeq \ar[d]^-{c_{g^{-1}}} &
B\widetilde{\cC_{G\pcom}(i_{({}_pZ(P))})}
\ar[dr]^-{{\widetilde{\phi_{i_{({}_pZ(P))}}}}}
\ar[d]_-{B\cC_{G\pcom}(c_g)} & \\
B{\tilde Q}\pcom \ar[r] & (B\widetilde{C_G({}_pZ(Q))})\pcom
\ar[r]^-\simeq & B\widetilde{\cC_{G\pcom}(i_{({}_pZ(Q))})}
\ar[r]_-{{\widetilde{\phi_{i_{({}_pZ(Q))}}}}} & B\widetilde{X'}
}
$$
where $i_{V}: BV \to BG\pcom$ denotes the map induced by the inclusion
of a subgroup $V\subseteq G$, and where the horizontal maps in the
middle square are given by lifting the standard homotopy equivalence given
by adjointness to the covers.
The first two squares are homotopy commutative by
construction, and the right-hand triangle commutes since
$[\vartheta] \in \lim_{\nu \in \A(X)}^0
[B\widetilde{\cC_X(\nu)},B\widetilde{X'}]$.
Hence this diagram produces an element
$[\zeta] \in \lim^0_{\widetilde{G}/\widetilde{P} \in
\bO_p^r(\widetilde{G})}[B{\tilde P}\pcom, B\widetilde{X'}]$. 
Denote the composition $B\widetilde{P}\pcom \to
B\widetilde{C_G({}_pZ(P))}\pcom \to B\widetilde{\cC_{G\pcom}(i_{{}_pZ(P)})}
\xrightarrow{\phi_{i_{({}_pZ(P))}}} B\widetilde{X'}$ by $\tilde
\psi_{\tilde G /\tilde P}$ (i.e., the coordinate of $[\zeta]$
corresponding to $\widetilde{G}/\widetilde{P}$), and let $\psi_{G/P}:
BP\pcom \to  BX'$ be the map constructed analogously, using
$\phi_{i_{({}_pZ(P))}}$ instead.
We want to lift $[\zeta]$ to a map
$$
\hocolim_{\widetilde{G}/\widetilde{P} \in \bO_p^r(\widetilde{G})}
(E\widetilde{G} \times_{\widetilde{G}} \widetilde{G}/\widetilde{P})\pcom \to
B\widetilde{X'}.
$$
By \cite[Prop.~XII.4.1~and~XI.7.1]{bk} (see also
\cite[Prop.~3]{wojtkowiak87} or \cite[Prop.~1.4]{JO96}) the
obstructions to doing this lie in
$$
\lim{}^{i+1}_{\widetilde{G}/\widetilde{P} \in \bO_p^r(\widetilde{G})}
\pi_i(\map(B\widetilde{P}\pcom,
B\widetilde{X'})_{\tilde \psi_{\tilde G / \tilde P}}), \quad i \geq 1.
$$
By construction, $\widetilde{\psi}_{\widetilde{G}/\widetilde{P}}$ and
$\psi_{G/P}$ fit into a homotopy pull-back square
$$
\xymatrix{
B\widetilde{P}\pcom \ar[r]^-{\widetilde{\psi}_{\tilde G / \tilde P}}
\ar[d] & B\widetilde{X'} \ar[d]\\
 BP\pcom \ar[r]^-{\psi_{G/P}} & B{X'}.
}
$$
By \cite[Lem.~4.9(1)]{CLN} the map $\psi_{G/P}$ is centric, so
Lemma~\ref{centricity} implies that $\widetilde{\psi}_{\tilde G /
  \tilde P}$ is centric as well. Hence by centricity and naturality,
the functor $\widetilde{G}/\widetilde{P} \mapsto
\pi_i(\map(B\widetilde{P}\pcom,
B\widetilde{X'})_{\widetilde{\psi}_{\tilde G / \tilde P}})$
identifies with the functor $\widetilde{G}/\widetilde{P} \mapsto
\pi_{i-1}(Z(\widetilde{P})\pcom)$. Since $\widetilde{G}$ is simple and
simply connected it now follows from the fundamental calculations of
Jackowski-McClure-Oliver \cite[Thm.~4.1]{JMO92} that
$$
\lim{}^{i+1}_{\bO^r_p(\widetilde{G})}\pi_{i-1}(Z(\widetilde{P})\pcom)=0
\text{, for }
i\geq 1.
$$
Hence by the homology decomposition theorem \cite[Thm.~4]{JMO92} we get a map
$$
B\widetilde{G}\pcom \xleftarrow{\simeq}
(\hocolim_{\widetilde{G}/\widetilde{P} \in \bO^r_p(\widetilde{G})}
(E\widetilde{G} \times_{\widetilde{G}}
\widetilde{G}/\widetilde{P})\pcom)\pcom \to B\widetilde{X'}
$$
which by construction is a map under $B\widetilde{\N_p}$, the
$p$-normalizer of a maximal torus in $B\widetilde{G}\pcom$. Dividing out by
$\cZ(\D)$ as explained in Construction~\ref{adjointconstruction}
produces the homotopy commutative diagram
$$
\xymatrix{
 & B\N_p \ar[dr] \ar[dl]\\
BG\pcom \ar[rr] & & BX'.
}
$$
It is now a short argument, given in detail
in \cite[Pf.~of~Thm.~2.2]{AGMV03}, to see that $BX = BG\pcom \to BX'$ is
a homotopy equivalence as wanted.

We want to show that $\Out(BX) \to \Out(\D_X)$ is an isomorphism: To see
surjectivity, note that if $\alpha \in \Out(\D_X)$, then by
\cite[Thm.~C]{AG05automorphisms} and \cite[Prop.~5.1]{AGMV03},
$\alpha$ corresponds to a unique map $\alpha' \in \Out(B\N,\{B\N_\sigma\})$.
Hence if $j: B\N \to BX$ is a maximal
torus normalizer, then repeating the above argument with respect to
the two maps $j:B\N \to BX$ and $j \circ \alpha': B\N \to BX$ gives a
map $BX \to BX$ realizing $\alpha \in \Out(\D_X)$. Finally we show
injectivity, essentially repeating the argument of
Jackowski-McClure-Oliver, cf.\ \cite[Pf.~of~Thm.~4.2]{JMO92}: We are
assuming that $BX \simeq BG\pcom$, for some compact connected
center-free Lie group $G$. As in the proof of
Proposition~\ref{centerred}\eqref{center3},
we have the following commutative diagram
$$
\xymatrix{
\Out(BG\pcom) \ar[d] \ar[r] &
\Out(\D_{G\pcom}) \ar[d] \\
\Out(B\widetilde{G}\pcom) \ar[r] & \Out(\D_{\widetilde{G}\pcom}).
}
$$
(Compare diagram \eqref{AMdiag10}, and note that we use that $\widetilde{
  \D_{G\pcom}} = \D_{\widetilde{G\pcom}}$ for the right-hand vertical
map, which uses Theorem~\ref{fundgrp} for compact Lie groups.)
The left-hand vertical map in the above diagram is injective since
factoring out by the center (via the quotient construction recalled in
Construction~\ref{adjointconstruction}) provides a left inverse (in
fact an actual inverse, though we do not need this here).
So we just have to see that $\Out(B\widetilde{G}\pcom) \to
\Out(\D_{\widetilde{G}\pcom})$ is injective. By
\cite[Thm.~C]{AG05automorphisms}, this map factors through
$\Out(B\widetilde{\N},\{B\widetilde{\N}_\sigma\})$ and
$\Out(B\widetilde{\N}, \{B\widetilde{\N}_\sigma\}) \to
\Out(\D_{\widetilde{G}\pcom})$ is an isomorphism, where
$\widetilde{\N}$ denotes the
maximal torus normalizer in $\widetilde{G}\pcom$. By the homology
decomposition theorem \cite[Thm.~4]{JMO92} and obstruction theory
\cite[Thm.~3.9]{JMO92} (cf.\ \cite[Prop.~3]{wojtkowiak87}), injectivity of
$\Out(B\widetilde{G}\pcom) \to \Out(B\widetilde{\N})$ follows from
Jackowski-McClure-Oliver's calculation of higher limits \cite[Thm.~4.1]{JMO92}
$$
\lim{}^i_{\widetilde{G}/\widetilde{P} \in
  \bO^r_p(\widetilde{G})}\pi_{i-1}(Z(\widetilde{P})\pcom) = 0,\quad i\geq
1.
$$
We conclude that $\Out(BG\pcom) \to \Out(\D_{G\pcom})$ is
also injective as claimed.

Finally, the last statement in Theorem~\ref{conn-classification} about
the homotopy type of $B\Aut(BX)$ follows by combining
\cite{AG05automorphisms} with what we have proved so far: The
Adams-Mahmud map factors as $\AM: B\Aut(BX) \to
B\Aut(B\N,\{B\N_\sigma\}) \to B\Aut(B\N)$,
where $B\Aut(B\N,\{B\N_\sigma\})$ is the
covering of $B\Aut(B\N)$ with respect to the subgroup
$\Out(B\N,\{B\N_\sigma\})$ of the fundamental group. Furthermore
\cite[\S 5]{AG05automorphisms} explains how killing elements in
$\pi_2(B\Aut(B\N,\{B\N_\sigma\}))$ constructs a space denoted
$B\aut(\D_X)$, whose universal cover is $B^2\cZ(\D_X)$, where
$B\cZ(\D_X) = (B\dZ(\D_X))\pcom$. It it furthermore shown there that the
fibration $B\aut(\D_X) \to B\Out(\D_X)$ is split, i.e., $B\aut(\D_X)$ has
the homotopy type of $(B^2\cZ(\D_X))_{h\Out(\D_X)}$. Composition gives a
map $B\Aut(BX) \to B\aut(\D_X)$, which is an isomorphism on $\pi_i$
for $i>1$ by construction, and an isomorphism on $\pi_1$ by what we
have shown in the first part of the theorem.

This concludes the proof of the main Theorem~\ref{conn-classification}.
\end{proof}

\begin{rem}[The fundamental group of a $p$-compact group] \label{pi1remark}
In \cite{AGMV03}, with M{\o}ller and Viruel, we established the
fundamental group formula for $p$-compact groups,
Theorem~\ref{fundgrp}, for $p$ odd, as a consequence of the
classification. In this remark we sketch how one by expanding
the proof of Theorem~\ref{conn-classification} somewhat can avoid the
reliance on Theorem~\ref{fundgrp}, hence proving Theorem~\ref{fundgrp}
also for $p=2$ this way---this was the strategy which we had originally
envisioned before Dwyer-Wilkerson \cite{DW06} provided an alternative
direct proof of Theorem~\ref{fundgrp} as we were writing this paper.

The fundamental group formula was not used in the reduction to simple
center-free groups, except for a reference in
the proof of Proposition~\ref{centerred}\eqref{center3}, where it, for
the purpose
of inductively proving Theorem~\ref{conn-classification}, was only
needed in the well-known case of compact Lie groups.
Hence  we can
assume that $X$ and $X'$ are simple center-free $p$-compact groups
with the same $\Z_p$-root datum $\D$, and that we know the fundamental
group formula for $X$ and Theorem~\ref{conn-classification} for
connected $p$-compact groups of lower cohomological dimension than
$X$. In the last section we constructed an element
$[\vartheta] \in \lim_{\nu \in \A(X)}^0
[B\widetilde{\cC_X(\nu)},B\widetilde{X'}]$.
By not passing to a universal cover, and hence not using the
fundamental group formula, a simplified version
of the same argument gives an element 
$[\bar \vartheta] \in \lim_{\nu \in \A(X)}^0
[B{\cC_X(\nu)},B{X'}]$.
The only non-obvious change is that since the center formula in
Lemma~\ref{centerformula-maxrank} does not hold in the presence of
direct factors isomorphic to $\D_{\SO(2n+1)\twocom}$ in the
$\Z_2$-root datum, one has to take the target of $\Xi$ in
\eqref{amlike} to be
$B\cZ(\D_{\cC_X(\nu)_1})$ (obtained as a quotient of
$B\cZ(\cC_\N(\mu)_1)$; cf.\ \cite[Lem.~5.1]{AG05automorphisms})
to obtain a homotopy equivalence. With this change the rest of
Section~\ref{proofI} proceeds as before, but ignoring everything
on the level of covers, and one constructs a  map like before without
any choices. Section~\ref{independence} has to be modified in the
following way: Instead of having maps being under the maximal torus
normalizer (which we now do not a priori know), we utilize instead
that potentially different maps agree in
$$
M = (\map(B\pi,B\aut(\D)))_{h\Aut(B\pi)},
$$
where $B\aut(\D) = (B^2\cZ(\D))_{h\Out(\D)}$,
which means that they are homotopy equivalent, by the description of
self-equivalences of non-connected groups, cf.,
Theorem~\ref{nonconnclassification}, which we know by induction.
From the element in $\lim^0$ one can easily get an isomorphism between
fundamental groups: Consider the diagram
$$
\xymatrix@C50pt{
 & H_*^{\Z_p}(B\N_p) \ar@{->>}[ld]_-j \ar[d] \ar@{->>}[rd]^-{j'} & \\
H_*^{\Z_p}(BX) & \colim_{\nu\in \A(X)} H_*^{\Z_p}(B\cC_X(\nu))
\ar[l]_(0.6){\cong} \ar[r]^(0.6){[\bar\vartheta]} & H_*^{\Z_p}(BX')
}
$$
where the vertical map is given by choosing a central rank one
subgroup $\rho$ of the $p$-normalizer $\N_p$ and considering the
corresponding inclusion
$\N_p \to \cC_X(j\circ \rho)$. (Here $H^{\Z_p}_n(Y) = \lim
H_n(Y;\Z/p^n)$.) Note that the maps $j$ and $j'$ are surjective by a
transfer argument, and that the indicated isomorphism follows by the
centralizer homology decomposition theorem \cite[Thm.~8.1]{dw:center}.
This shows that the kernel of $j$ is contained in the kernel of
$j'$. However, since we could reverse the role of $X$ and $X'$ in all
the previous arguments (we have not used any special model for $X$),
we conclude by symmetry that the kernel of $j$ equals the kernel of
$j'$, and in particular the kernels of the maps $j: \pi_1(\D)
\to \pi_1(X)$ and $j': \pi_1(\D) \to \pi_1(X')$ agree, and they are
surjective by Proposition~\ref{constrpi1map}.
If $\D$ is exotic then $\pi_1(\D) = 0$ by
Theorem~\ref{rootdata-classification}\eqref{rdc-part2}, so there is nothing to
prove. If $\D$ is not exotic, then by
Theorem~\ref{rootdata-classification}\eqref{rdc-part1},  $\D = \D_G
\otimes_{\Z} \Z_p$ for a compact connected Lie group $G$ and so we can assume $BX \simeq
BG\pcom$. Hence $j: \pi_1(\D) \xrightarrow{\cong} \pi_1(X)$ by Lie
theory (cf., \cite[\S 4, no.\ 6, Prop.\ 11]{bo9}), so the same holds for $X'$.
\end{rem}


\section{Proof of the corollaries of Theorem~\ref{conn-classification}}\label{consequences-section}

In this section we prove the Theorems 1.3--1.6 from the introduction.
The first theorem is the maximal torus conjecture:

\begin{proof}[Proof of Theorem~\ref{maxtorus}]
The proof in the connected case is an extension of
\cite[Pf.~of~Thm.~1.10]{AGMV03}, where a partial result excluding the
case $p=2$ was given: Assume that
$(X,BX,e)$ is a {\em connected} finite loop space with a maximal torus
$i: BT \to BX$. Let $W$ denote the set of conjugacy classes of
self-equivalences $\varphi$ of $BT$ such that $i\varphi$ is conjugate
to $i$. It is straightforward to see (consult e.g.,
\cite[Pf.~of~Thm.~1.10]{AGMV03}) that $BT\pcom \to
BX\pcom$ is a maximal torus for the $p$-compact group $X\pcom$ and
that $\F_p$-completion allows us to identify $(W_{X\pcom}, L_{X\pcom})$
with $(W, L\otimes \Z_p)$, where $L = \pi_1(T)$. In particular
$(W, L)$ is a finite $\Z$-reflection group and all reflections have order
$2$. Furthermore, for a fixed $\Z$-reflection group $(W,L)$, there is a
bijection between $\Z$-root data with underlying $\Z$-reflection
group $(W, L)$ and $\Z_2$-root data with underlying
$\Z_2$-reflection group $(W, L\otimes \Z_2)$ given by the assignments
$\D \mapsto \D \otimes_{\Z} \Z_2$ and $(W, L\otimes \Z_2, \{\Z_2b_\sigma\})
\mapsto (W, L, \{L \cap \Z_2b_\sigma\})$, as is seen by examining the
definitions. Let $\D$ be the $\Z$-root datum with
underlying $\Z$-reflection group $(W,L)$ corresponding to
$\D_{X\twocom}$.

By the classification of compact connected
Lie groups, cf.\ \cite[\S 4, no.\ 9, Prop.\ 16]{bo9}, there is a
(unique) compact
connected Lie group $G$ with maximal torus $i': T \to G$ inducing an
isomorphism of $\Z$-root data $\D_G \cong \D$.
The $\Z_p$-root data of $X\pcom$ and $G\pcom$ are isomorphic at all
primes, since the root data at odd primes are determined by $(W,L
\otimes \Z_p)$. Theorem~\ref{conn-classification} hence implies that
for each $p$, we have a homotopy equivalence $\phi_p: BX\pcom \to
BG\pcom$ such that
$$
\xymatrix{
 & BT\pcom \ar[dl]_-{i\pcom} \ar[dr]^-{i'\pcom} \\
BX\pcom \ar[rr]^-{\phi_p} & & BG\pcom
}
$$
commutes.
As in
\cite[Pf.~of~Thm.~1.10]{AGMV03} we see that $H^*(BX;\Q) \to
H^*(BT;\Q)^W$ is an isomorphism and since the same is true for $BG$ we
also have a map $BX_\Q \to BG_\Q$ under $BT$. We have the following diagram
\begin{equation}\label{sullivandiag}
\xymatrix{
{\prod_{p}} {BX\pcom} \ar[r] \ar[d]^-\simeq & (\prod_{p} BX\pcom)_\Q
\ar[d]^-\simeq & BX_{\Q} \ar[l] \ar[d]^-\simeq \\
{\prod_{p}} {BG\pcom} \ar[r] & (\prod_{p} BG\pcom)_\Q & BG_{\Q}. \ar[l]
}
\end{equation}
The left-hand square in this diagram is homotopy commutative by
construction. For the right-hand side note that, since all maps in
\eqref{sullivandiag} are under $BT$, the following diagram commutes
$$
\xymatrix{
H_*(\prod_pBX\pcom;\Q) \ar[dd]^-\cong & & H_*(BX;\Q) \ar[ll] \ar[dd]^-\cong\\
  & H_*(BT;\Q)_W \ar[ul] \ar[ur]_-\cong \ar[dr]^-\cong \ar[dl]\\
H_*(\prod_pBG\pcom;\Q) & & H_*(BG;\Q). \ar[ll]
}
$$
This implies commutativity on the level of homotopy groups, and since
the involved spaces are all products of rational Eilenberg-Mac\,Lane
spaces (since they have homotopy groups only in even degrees) this
implies that the diagram \eqref{sullivandiag} homotopy commutes. By
changing the maps up to homotopy we can hence arrange that the diagram
strictly commutes, and taking homotopy pull-backs produces a homotopy
equivalence $BX \to BG$, by the Sullivan arithmetic
square \cite[VI.8.1]{bk}, as wanted.
This
proves that every connected finite loop space $BX$ with a
  maximal torus is homotopy
equivalent to $BG$ for some compact connected Lie group $G$, and in
the course of the analysis, we furthermore saw that $G$ is unique
(since we can read off the $\Z$-root datum of $G$ from $BX$).

We now give the description of $B\Aut(BG)$, also providing a quick
description of how the $p$-completed results are used. By
\cite[Thm.~1.4]{dw:center} $B^2Z(G) \simeq B\Aut_1(BG)$. (This is a
consequence of that $B^2\cZ(G\pcom) \simeq B\Aut_1(BG\pcom)$.) By
\cite[Cor.~3.7]{JMO95} $\Out(G) \cong \Out(BG)$. (Since $(G/T)^{hT}$
is homotopically discrete with components the Weyl group, as in
\cite[Lem.~4.1]{AGMV03} and \cite[Thm.~2.1]{notbohm91}, any map $f: BG
\to BG$ gives rise to a map $\phi: BT \to BT$ over $f$, unique up to
Weyl group conjugation and since $\phi\pcom \in \Aut(\D_{G\pcom})$ for
all $p$ one sees that $\phi \in \Aut(\D_G)$; now $\phi$ determines the
collection $\{ \phi\pcom \}$ which by $p$-complete results
determines $\{ f\pcom\}$, which again determines $f$ by the arithmetic
square.) Finally by a theorem of de Siebenthal \cite[Ch.~I, \S 2, no.\
2]{desiebenthal56} (see also \cite[\S 4, no.\ 10]{bo9}), the short
exact sequence $1 \to G/Z(G) \to \Aut(G) \to \Out(G) \to 1$ is split,
so the fibration $B\Aut(BG) \to B\Out(BG)$ has a section. (See also
\cite[\S 6]{AG05automorphisms}.) Taken together these facts establish
that $B\Aut(BG) \simeq (B^2Z(G))_{h\Out(G)}$ as claimed.

The non-connected case follows easily from the connected, using our
knowledge of self-maps of classifying spaces of compact connected Lie
groups: Suppose that $X$ is potentially non-connected, let $X_1$ be
the identity component, and
$\pi$ the component group. Note that homotopy classes of spaces $Y$
such that $P_1(Y)$ is homotopy equivalent to $B\pi$ and the fiber
$Y\langle 1 \rangle$ is homotopy equivalent to $BX_1$ can
alternatively be described as homotopy classes of fibrations $Z \to Y
\to K$, with $Z$ is homotopic to $BX_1$ and $K$ is homotopic to
$B\pi$. (A homotopy equivalence of fibrations means a compatible
triple of homotopy equivalences between fibers, total spaces, and
base spaces.)

By the first part of the proof we know that $BX_1 \simeq BG_1$ for a
unique compact connected Lie group $G_1$. Homotopy classes of
fibrations with base homotopy equivalent to $B\pi$ and fiber homotopy
equivalent to $BX_1$ are classified by $\Out(\pi)$-orbits on the set
of free homotopy classes $[B\pi, B\Aut(BG_1)]$. By the results in the
connected case the space $B\Aut(BG_1)$ sits in a split fibration
$$
B^2 Z(G_1) \to B\Aut(BG_1) \to B\Out(G_1).
$$
Hence $\Out(\pi)$-orbits on the set $[B\pi,B\Aut(BG_1)]$ correspond to
$(\Out(\pi) \times \Out(G_1))$-orbits on the set
$$
\coprod_{\alpha \in \Hom(\pi,\Out(G_1))} H^2_\alpha(\pi;Z(G_1)).
$$
This agrees with the classification of isomorphism classes of group
extensions of the form $1 \to H \to \text{?} \to K \to 1$, where $H$
is isomorphic to $G_1$ and $K$ is isomorphic to $\pi$. (An isomorphism
of group extensions is here a compatible triple of isomorphisms.)
Since the identity component $H$ is necessarily a characteristic
subgroup, isomorphism classes of group extensions as above are in
one-to-one correspondence with isomorphism classes of compact Lie
groups with identity component isomorphic to $G_1$ and component group
isomorphic to $\pi$. These equivalences put together give that there
is a one-to-one correspondence
between homotopy classes of spaces $Y$ such that $P_1(Y) \simeq B\pi$
and $Y\langle 1 \rangle \simeq BG_1$ and isomorphism classes of compact
Lie groups with identity component isomorphic to $G_1$ and component
group isomorphic to $\pi$. Hence our $BX$ is homotopy equivalent to
$BG$ for a unique compact Lie group $G$, completing the proof of the
theorem.
\end{proof}

\begin{proof}[Proof of Theorem~\ref{Steenrod}]
Let $Y$ be a space such that $H^*(Y;\F_2)$ is a graded polynomial
algebra of finite type. Let $V=H_1(Y;\F_2)$ (dual to $H^1(Y;\F_2)$)
and let $Y'$ denote the
fiber of the classifying map $Y\to BV$. Clearly $Y'$ is connected and
since $\pi_1(BV)$ is a finite $2$-group it follows from \cite{dwyer74}
that the Eilenberg-Moore spectral sequence for the fibration $Y' \to Y\to BV$
converges strongly to $H^*(Y';\F_2)$. The map $H^*(BV;\F_2) \to
H^*(Y;\F_2)$ is an isomorphism in degree $1$ and hence injective since
$H^*(Y;\F_2)$ is a polynomial algebra, so $H^*(Y;\F_2)$ is free over
$H^*(BV;\F_2)$. Hence the spectral sequence collapses and we get an
isomorphism of rings (but not necessarily of algebras over the
Steenrod algebra) $H^*(Y;\F_2) \cong H^*(Y';\F_2) \otimes
H^*(BV;\F_2)$. In particular $H_1(Y';\F_2)=0$ so by
\cite[Prop.~VII.3.2]{bk} $Y'$ is $\F_2$-good and
$\pi_1(Y'\twocom)=0$. So to
prove the theorem, we can without restriction assume that $Y'$ is
$\F_2$-complete and simply connected.

Write $\pi_2(Y') \cong F \oplus T$, where $F$ is a finitely generated
free $\Z_2$-module and $T$ is a finite abelian $2$-group, and let
$Y''$ be the fiber of the map $Y' \to B^2F$. The induced homomorphism
$H_2(Y';\F_2) \to H_2(B^2F;\F_2)$ is an epimorphism so $H^*(B^2F;\F_2)
\to H^*(Y';\F_2)$ is injective. As above we obtain an isomorphism
$H^*(Y';\F_2) \cong H^*(Y'';\F_2) \otimes H^*(B^2F;\F_2)$ as rings.

By construction $Y''$ is simply connected, $\pi_2(Y'')$ is finite, and
by the fibre lemma \cite[Lem.~II.5.1]{bk} $Y''$ is
$\F_2$-complete. Since $H^*(Y'';\F_2)$ is polynomial, the
Eilenberg-Moore spectral sequence shows that $H^*(\Omega Y'';\F_2)$ is
$\F_2$-finite, so $Y''$ is the classifying space of a connected
$2$-compact group. The first part of Theorem~\ref{Steenrod} now
follows from the classification Theorem~\ref{splitobj}. The second
part follows from this using the calculation of the mod $2$ cohomology
of the simple simply connected Lie groups, cf.\ \cite[Thm.~5.2]{kono75}.
\end{proof}

\begin{rem} \label{konoremark}
In addition to the list in Theorem~\ref{Steenrod}, the only
polynomial rings arising as $H^*(BG;\F_2)$ for a simple compact
connected Lie group $G$ are $\F_2[x_2, x_3, \ldots,
x_n]$ for $G=\SO(n)$, $n\geq 5$, and $\F_2[x_2, x_3, x_8, x_{12}, \ldots,
x_{8n+4}]$ for $G=\PSp(2n+1)$, $n\geq 0$, cf.\ \cite[Thm.~5.2]{kono75}.
It is conceivable that any graded polynomial algebra of finite type
which is the mod $2$ cohomology ring of a space, is a tensor product
of these factors and the ones listed in Theorem~\ref{Steenrod}.
\end{rem}

\begin{proof}[Proof of Theorem~\ref{cohomologicaluniqueness}]
The first statement claims that a polynomial $\F_2$-algebra $A^*$ with
given action of the Steenrod algebra $\mathcal{A}_2$, can be realized
by at most one space $Y$, up to $\F_2$-equivalence, if $A^*$ has all
generators in degree $\geq 3$. For this notice that, as in the proof
of Theorem~\ref{Steenrod},
the assumptions assure that $Y\twocom \simeq BX$ for a simply
connected $2$-compact group $X$. Using
the classification Theorem~\ref{splitobj}, the statement can
now easily be checked as done in Proposition~\ref{algebraicweyl} below.

We now prove the second statement, that for any polynomial $\F_2$-algebra
$P^*$, there are only finitely many spaces $Y$, up to $\F_2$-equivalence,
with $H^*(Y;\F_2) \cong P^*$ as rings. By the proof of Theorem~\ref{Steenrod}, 
we can assume that $Y$ is $2$-complete, and any such $Y$ fits in a fibration
sequence $Y' \to Y \to BV$ where $V= H_1(Y;\F_2)$. It also follows that
$H^*(Y';\F_2)$ is a polynomial ring, uniquely determined as an algebra
over the Steenrod algebra from $H^*(Y;\F_2)$. In particular $Y' \simeq
BX$ for a connected $2$-compact group $X$. By Proposition~\ref{replemma},
$\Out(\D_X)$ only contains finitely many $2$-subgroups up to conjugation.
Hence the description of
$B\Aut(BX)$ in Theorem~\ref{conn-classification} implies that
$[BV,B\Aut(BX)]$ is finite, so it is enough
to see that there are only a finite number of possibilities for $Y'$
given $H^*(Y';\F_2)$ as an algebra over the Steenrod algebra. This
again follows easily from the classification of $2$-compact groups:
The rank of $X$ is bounded above by the Krull dimension of the cohomology
ring, so by the classification of $p$-compact groups,
Theorem~\ref{conn-classification}, it is hence enough to see that there
are only a finite number of $\Z_2$-root data with rank less than a fixed
rank. This is the result of Proposition~\ref{finiteoffixedrank}.
\end{proof}

We now give a proof of the auxiliary uniqueness result referred to in the proof
of Theorem~\ref{cohomologicaluniqueness}.

\begin{prop}\label{algebraicweyl}
Suppose $X$ is a $2$-compact group of the form $BX \simeq BG\twocom
\times BDI(4)^s$, for a simply connected compact Lie group $G$ and $s
\geq 0$, such that $H^*(BX;\F_2)$ is a polynomial algebra. Then $X$
has a unique maximal elementary abelian $2$-subgroup
$\nu: BE \to BX$, and the Weyl group $(W(\nu),E)$ together
with the homomorphism $H^6(BX;\F_2)
\to H^6(BE;\F_2)$ is an invariant of $H^*(BX;\F_2)$ as an algebra over
the Steenrod algebra $\mathcal{A}_2$, which uniquely determines $BX$ up to
homotopy equivalence. In particular $BX$ is uniquely determined up to
homotopy equivalence by $H^*(BX;\F_2)$ as an $\mathcal{A}_2$-algebra.
\end{prop}

\begin{proof}
By Lannes theory \cite[Thm.~0.4]{lannes92} homotopy classes of maps
from classifying spaces of elementary abelian
$2$-groups to $BX$ are determined by $H^*(BX;\F_2)$ as an
algebra over the Steenrod algebra. Furthermore, the fact that
$H^*(BX;\F_2)$ is
assumed to be a polynomial algebra, guarantees  that there is only one
maximal elementary abelian $2$-subgroup $\nu: BE \to BX$, up to
conjugation, by \cite[Cor.~10.7]{quillen71cohom} together with the fact
that this is true for $B\DI(4)$. (This can also be deduced using the
unstable algebra techniques of \cite{AW80} and \cite{HLS93}, or simply
by inspecting the calculations of Griess \cite{griess} below.)
Let $(W(\nu),E)$ denote its Weyl group, which also only depends on
$H^*(BX;\F_2)$, and we view this as a pair with $W(\nu)$ a subgroup of
$\GL(E)$.
By \cite[Thm.~5.2]{kono75}, $G$ is a direct product of the groups
$\SU(n)$, $\Sp(n)$, $\Spin(7)$, $\Spin(8)$, $\Spin(9)$, $G_2$ and
$F_4$. In these cases, if $E$ is the maximal elementary abelian
$2$-subgroup of $G$, $W_{G\twocom}(\nu) = N_G(E)/C_G(E)$ whose
structure is well-known in these cases, e.g., by computations of
Griess \cite[\S 5, Thm.~6.1 and Thm.~7.3]{griess}. (See also
\cite[Prop.~3.2]{VV02} for details on the cases $\Spin(8)\subseteq
\Spin(9)\subseteq F_4$.) Since $W_{\DI(4)}(\nu) = \GL_4(\F_2)$ by
construction \cite{DW93}, it follows that the
group $(W(\nu),E)$ is a direct product of the following (with matrix
groups acting on columns):
$$
W_{\SU(n)}(\nu) = (\Sigma_n, V'_{n-1}),\qquad W_{\Sp(n)}(\nu) =
(\Sigma_n, V_n),
$$
$$
W_{G_2}(\nu) = \GL_3(\F_2),\qquad W_{\DI(4)}(\nu) = \GL_4(\F_2),
$$
$$
W_{\Spin(7)}(\nu) = \left[
\begin{array}{c|c}
1 & \begin{array}{ccc}
* & * & *
\end{array} \\
\hline
\begin{array}{c}
0 \\
0 \\
0
\end{array} & \GL_3(\F_2)
\end{array}
\right],\qquad
W_{\Spin(8)}(\nu) = \left[
\begin{array}{c|c}
\begin{array}{cc}
1 & 0\\
0 & 1
\end{array} &
\begin{array}{ccc}
* & * & *\\
* & * & *
\end{array} \\
\hline
\begin{array}{cc}
0 & 0 \\
0 & 0 \\
0 & 0
\end{array} & \GL_3(\F_2)
\end{array}
\right],
$$
$$
W_{\Spin(9)}(\nu) = \left[
\begin{array}{c|c}
\begin{array}{cc}
 1 & *\\
0 & 1
\end{array} &
\begin{array}{ccc}
* & * & *\\
* & * & *
\end{array} \\
\hline
\begin{array}{cc}
0 & 0 \\
0 & 0 \\
0 & 0
\end{array} & \GL_3(\F_2)
\end{array}
\right],\qquad
W_{F_4}(\nu) = \left[
\begin{array}{c|c}
\GL_2(\F_2) & \begin{array}{ccc}
* & * & *\\
* & * & *
\end{array} \\
\hline
\begin{array}{cc}
0 & 0 \\
0 & 0 \\
0 & 0
\end{array} & \GL_3(\F_2)
\end{array}
\right].
$$
Here $V_n$ is the $n$-dimensional permutation module for
$\F_2[\Sigma_n]$ and $V'_{n-1}$ is the $(n-1)$-dimensional submodule
consisting of elements with coordinate sum $0$. 
The Weyl group $W_{\Sp(n)}(\nu) =
(\Sigma_n, V_n)$ decomposes as $(\Sigma_n, V'_{n-1}) \times (1,
L)$, $L = {V_n}^{\Sigma_n} \cong \F_2$, when $n$ is odd, $n\geq 3$.
However, after this decomposition, all the listed pairs satisfy that
$V$ is indecomposable as an $\F_2[W]$-module. (Note that this is a priori
stronger than just saying that $(W,V)$ does not split as a product;
however, it follows from unstable algebra techniques
\cite[Secs.~5 and~7]{neusel06} that $(W,V)$ is an $\F_2$-reflection
group, and hence the two notions are actually equivalent.)
By the Krull-Schmidt theorem \cite[6.12(ii)]{CR81} this
decomposition $V =  V_1 \oplus \ldots \oplus V_m$ as $\F_2[W]$-modules is
unique up to permutation, and since $W_i$ is characterized as the
pointwise stabilizer of $V_1 \oplus \ldots \oplus \widehat{V_i} \oplus
\ldots \oplus V_m$, we get that the decomposition of $(W,V)$ as a
product is unique as well, up to permutation.
Thus the structure of $(W(\nu),E)$ as a product of finite
indecomposable groups, is an invariant which almost
characterizes $BX$ up to homotopy equivalence, except that for $n$ odd, $n\geq
3$, the group $(\Sigma_n, V'_{n-1})$ arises from
both $\SU(n)\twocom$ and $\Sp(n)\twocom$. However since
$H^6(B\SU(n);\F_2) \to H^6(BV'_{n-1};\F_2)$ is injective and
$H^6(B\Sp(n);\F_2) \to H^6(BV'_{n-1};\F_2)$ is trivial, we conclude
that $(W(\nu), E)$ together with the homomorphism $H^6(BX;\F_2) \to
H^6(BE;\F_2)$ characterizes $BX$ up to homotopy equivalence. This proves the
proposition since both data are determined by the
$\mathcal{A}_2$-action on $H^*(BX;\F_2)$.
\end{proof}

\begin{proof}[Proof of Theorem~\ref{nonconnclassification}]
It is obvious that there is a one-to-one correspondence between
isomorphism classes of  $p$-compact groups with identity component
isomorphic to $X_1$ and component group
isomorphic to $\pi$, and equivalence classes of fibration sequences $F
\to E \xrightarrow{p} B$ with $F$ homotopy equivalent to $BX_1$ and $B$
homotopy equivalent to $B\pi$.  It is likewise obvious that in this
case $B\Aut(E)$ is homotopy equivalent to $B\Aut(p)$, where $\Aut(p)$
is the space of self-homotopy equivalences of the fibration $p$.

By the classification of fibrations (see \cite{DKS89}),
equivalence classes of such fibrations are in one-to-one correspondence
with $\Out(B)$-orbits on $[B,B\Aut(F)]$, and the space $B\Aut(p)$
equals $(\map(B,B\Aut(F))_{C(p)})_{h\Aut(B)}$, where $C(p)$ denotes the
$\Out(B)$-orbit on $[B,B\Aut(F)]$ of the element classifying $p: E \to
B$. 

The above considerations completely reduces the proof of the theorem,
to our classification theorem for connected $p$-compact groups,
Theorem~\ref{conn-classification}, except for the finiteness statement.
For this note that Proposition~\ref{replemma} implies that
$[B\pi,B\Out(\D)]$ is finite so the finiteness of $[B\pi,B\aut(\D)]$
follows.
\end{proof}

\begin{rem} \label{applrem}
As stated in the introduction, our classification also shows that
Bott's theorem on the cohomology of $X/T$, the Peter-Weyl theorem, as
well as Borel's characterization of when centralizer of elements of
order $p$ are connected, stated as Theorems 1.5, 1.6 and 1.9 in
\cite{AGMV03}, hold verbatim for $2$-compact groups. To prove these
results it suffices by Theorem~\ref{splitobj} to check them for
$\DI(4)$, since they are well known for compact Lie groups. For
$\DI(4)$ one argues as follows: Bott's theorem \cite[Thm.~1.5]{AGMV03}
follow from \cite[Thm.~1.8(2)]{DW93}, the
Peter-Weyl theorem \cite[Thm.~1.6]{AGMV03} is a result of Ziemia{\'n}ski
\cite{ziemianski06}, and it is trivial to check that
\cite[Thm.~1.9]{AGMV03} hold.
\end{rem}


\section{Appendix: Properties of $\Z_p$-root data} \label{rootdatumappendix}

The purpose of this section is to establish some general results about
$\Z_p$-root data of $p$-compact groups, needed in the proof of the
main theorem. The analogous results for $\Z$-root data and compact Lie
groups are often well known; see \cite{bo9,demazure62}. We build on
the paper \cite{DW05} by Dwyer-Wilkerson and our earlier paper
\cite{AG05automorphisms}.

We briefly recall the definition of root data from the introduction:
For an integral domain $R$, an {\em $R$-reflection group} is a pair
$(W,L)$ where $L$ is a finitely generated free $R$-module and $W$ is a subgroup of
$\Aut_R(L)$ generated by reflections (i.e.\ elements $\sigma\in
\Aut_R(L)$ such that $1-\sigma\in \End_R(L)$ has rank one). If $R$ is a
principal ideal domain, we define an {\em $R$-root
  datum} to be a triple $\D = (W,L,\{R b_\sigma\})$ where $(W,L)$ is a
finite $R$-reflection group and for each reflection $\sigma\in W$, $R
b_\sigma$ is a rank one submodule of $L$ with $\im(1-\sigma) \subseteq
R b_\sigma$ and $w(R b_\sigma) = R b_{w\sigma w^{-1}}$ for all $w\in
W$. If $R \to R'$ is a monomorphism of integral domains, and $\D$
an $R$-root datum, we can define an $R'$-root datum by $\D \otimes_{R}
R'= (W,L \otimes_{R} R',\{Rb_\sigma \otimes_R {R'}\})$.

The element $b_\sigma$, defined up to a unit in $R$, is called the
{\em coroot} associated to $\sigma$. By definition $b_\sigma$
determines a unique linear map $\beta_\sigma : L \to R$ called the
associated {\em root} such that $\sigma(x) = x + \beta_\sigma(x)
b_\sigma$ for $x\in L$. Define the {\em coroot
  lattice} $L_0\subseteq L$ as the sublattice spanned by the coroots
$b_\sigma$ and the {\em fundamental group} of $\D$ by $\pi_1(\D) =
L/L_0$. In general $R b_\sigma \subseteq \ker(N)$,
where $N = 1 + \sigma + \ldots + \sigma^{|\sigma|-1}$ is the norm
element (cf.\ the proof of Lemma~\ref{decomp} below),
so giving an $R$-root datum with underlying reflection group
$(W,L)$ corresponds to choosing a cyclic $R$-submodule of
$H^1(\left<\sigma\right>;L)$ for each conjugacy class of reflections
$\sigma$. In particular for $R = \Z_p$, $p$ odd, the notions of a
$\Z_p$-reflection group and a $\Z_p$-root datum agree. If $R$ has
characteristic zero, an $R$-root datum, or an $R$-reflection group, is
called {\em irreducible} if the representation $W \to \GL(L \otimes_R
K)$ is irreducible, where $K$ denotes the quotient field of $R$, and
it is said to be {\em exotic} if furthermore the values of the
character of this representation are not all contained in $\Q$.

\smallskip

We are now ready to state the classification of $\Z_p$-root data,
which follows easily from the classification of finite
$\Z_p$-reflection groups \cite[Thm.~11.1]{AGMV03}. This classification
is again based on the classification of finite $\Q_p$-reflection
groups \cite{CE74} \cite{DMW92} which states that for a fixed prime
$p$, isomorphism classes of finite irreducible $\Q_p$-reflection
groups are in natural one-to-one correspondence with isomorphism
classes of finite irreducible $\C$-reflection groups $(W,V)$ \cite{ST54}
for which the values of the character of $W \to \GL(V)$ are embeddable
in $\Q_p$; see e.g., \cite[Table~1]{kksa:thesis} for an explicit list
of groups and primes.

\begin{thm}[The classification of $\Z_p$-root data; splitting version]
  \label{rootdata-classification}
\InsertTheoremBreak
\begin{enumerate}
\item \label{rdc-part1}
Any $\Z_p$-root datum is isomorphic to a $\Z_p$-root datum of the form $(\D_1
\otimes_{\Z} \Z_p) \times \D_2$, where $\D_1$ is a $\Z$-root datum and $\D_2$
is a direct product of exotic $\Z_p$-root data.
\item \label{rdc-part2}
There is a one-to-one correspondence between isomorphism classes of
exotic $\Z_p$-root data and isomorphism classes of exotic
$\Q_p$-reflection groups given by the assignment $\D = (W, L, \{\Z_p
b_\sigma\}) \rightsquigarrow (W, L\otimes_{\Z_p} \Q_p)$. Moreover
$\pi_1(\D)=0$ for any exotic $\Z_p$-root datum $\D$.
\end{enumerate}
\end{thm}

\begin{proof}
For any $\Z_p$-root datum $\D = (W, L, \{\Z_p b_\sigma\})$,
\cite[Thm.~11.1]{AGMV03} gives a splitting $(W, L) \cong (W_1,
L_1\otimes_{\Z} \Z_p) \times (W_2, L_2)$ of $(W,L)$, where $(W_1, L_1)$ is a finite
$\Z$-reflection group and $(W_2, L_2)$ is a direct product of exotic
$\Z_p$-reflection groups. It follows by definition that there are
unique $\Z_p$-root data $\D'$ and $\D_2$ with underlying reflection
groups $(W_1, L_1\otimes_{\Z} \Z_p)$ and $(W_2, L_2)$ such that $\D \cong
\D' \times \D_2$, and by the same argument $\D_2$ splits as a direct
product of exotic $\Z_p$-root data. Furthermore, writing $\D' = (W_1,
L_1\otimes_{\Z} \Z_p, \{\Z_p b_\sigma\})$ it is clear that $\D' \cong
\D_1 \otimes_{\Z} \Z_p$ where $\D_1 = (W_1, L_1, \{L_1 \cap \Z_p
b_\sigma\})$. This proves \eqref{rdc-part1}.

By \cite[Thm.~11.1]{AGMV03}, the assignment $(W, L) \rightsquigarrow
(W, L\otimes_{\Z_p} \Q_p)$ establishes a one-to-one correspondence
between exotic $\Z_p$-reflection groups up to isomorphism and exotic
$\Q_p$-reflection groups up to isomorphism. To prove the first part of
\eqref{rdc-part2} it thus suffices to show that any exotic
$\Z_p$-reflection group $(W,L)$ can be given a unique $\Z_p$-root
datum structure. For $p>2$ this holds since $H^1(\left<\sigma
\right>;L) = 0$, cf.\ the discussion in the beginning of this
section. For $p=2$, $(W,L) \cong (W_{\DI(4)}, L_{\DI(4)})$  where the
claim follows by direct inspection (cf.\ \cite[Rem.~7.2]{DW05}).

For any $\Z_p$-root datum $\D = (W, L, \{\Z_p b_\sigma\})$, the
formula $\sigma(x) = x + \beta_\sigma(x) b_\sigma$ shows that the
coroot lattice $L_0$ contains the lattice spanned by the elements
$(1-w)(x)$, $w\in W$, $x\in L$. Hence the final claim follows from the
fact that $H_0(W;L)=0$ for any exotic $\Z_p$-reflection group $(W,L)$
\cite[Thm.~11.1]{AGMV03}.
\end{proof}

\subsection{Root datum, normalizer extension, and root subgroups of a $p$-compact group}
\label{rootdatumsub}

For any connected $p$-compact group $X$ with maximal torus $T$, the
Weyl group $W_X$ acts naturally on $L_X = \pi_1(T)$ as a finite
$\Z_p$-reflection group \cite[Thm.~9.7(ii)]{DW94}. For $p$ odd,
$H^1(\left<\sigma\right>;L) = 0$ for any reflection $\sigma$, so the
finite $\Z_p$-reflection group $(W_X, L_X)$ gives rise to a unique
$\Z_p$-root datum $\D_X$. The construction of root data for connected
$2$-compact groups, in the present form, is due to Dwyer-Wilkerson
\cite[\S 9]{DW05}: Let $\dT$ be
the discrete approximation to $T$, $\dN_X$ the discrete approximation to
the maximal torus normalizer $\N_X$ and $\sigma\in W_X$ a
reflection. Define $\dT^+(\sigma) = \dT^{\left<\sigma\right>}$ and let
$\dT_0^+(\sigma)$ denote its maximal divisible subgroup. Then
$X(\sigma) = \cC_X(\dT^+_0(\sigma))$ is a connected $2$-compact group
with Weyl group $\left<\sigma\right>$ and $\dN(\sigma) =
C_{\dN_X}(\dT^+_0(\sigma))$ is a discrete approximation to its maximal
torus normalizer. Furthermore, let $\dT^-_0(\sigma)$ denote the
maximal divisible subgroup of $\dT^-(\sigma) = \ker(\dT
\xrightarrow{1+\sigma} \dT)$ and define the {\em root subgroup}
$\dN_{X,\sigma}$ by
$$
\dN_{X,\sigma} = \{ x\in \dN(\sigma) \,\mid\, \exists\, y\in \dT^-_0(\sigma)
: \text{$x$ is conjugate to $y$ in $X(\sigma)$} \}.
$$
Then there is a short exact sequence
\begin{equation} \label{Nsigmaext}
1 \to \dT_0^-(\sigma) \to \dN_{X,\sigma} \to \left<\sigma\right> \to 1,
\end{equation}
and we define
$$
\Z_2 b_\sigma =
\begin{cases}
\im(L_X \xrightarrow{1-\sigma} L_X) & \text{if \eqref{Nsigmaext}
  splits,} \\
\ker(L_X \xrightarrow{1+\sigma} L_X) & \text{otherwise}.
\end{cases}
$$
The root datum of $X$ is then the $\Z_2$-root datum $\D_X = (W_X,
L_X, \{\Z_2 b_\sigma\})$; see \cite[\S 6 and 9]{DW05} and
\cite{AG05automorphisms} for a further discussion.

Conversely, the maximal torus normalizer and the root subgroups of a
connected  $p$-compact group can be reconstructed from its root datum: For a
$\Z_p$-root datum $\D = (W, L, \{\Z_p b_\sigma\})$ there is
\cite[Def.~6.15]{DW05} \cite[\S 3]{AG05automorphisms} an algebraically
defined extension
$$
1 \to \dT \to \dN_{\D} \to W \to 1
$$
called the {\em normalizer extension} with a subextension $1 \to
\dT^-_0(\sigma) \to \dN_{\D,\sigma} \to \left<\sigma\right> \to 1$ for
each reflection $\sigma\in W$. For a connected $p$-compact group $X$
with $\Z_p$-root datum $\D_X$ there is an isomorphism of extensions
\cite[Prop.~1.10]{DW05} \cite[Thm.~3.2(2)]{AG05automorphisms}
$$
\xymatrix{
1 \ar[r] & \dT \ar@{=}[d] \ar[r] & \dN_{\D} \ar[d]^-\cong \ar[r] & W
\ar@{=}[d] \ar[r] & 1\\
1 \ar[r] & \dT \ar[r] & \dN_X \ar[r] & W \ar[r] & 1
}
$$
sending the {\em root subgroups} $\dN_{\D,\sigma}$ to $\dN_{X,\sigma}$
for all reflections $\sigma$, and any such isomorphism is unique up to
conjugation by an element in $\dT$.

We define $B\N_\D$ as the fiber-wise $\F_p$-completion \cite[Ch.~I, \S
8]{bk} of  $B\dN_\D$ and likewise introduce the (non-discrete) root subgroups
$B\N_{X,\sigma}$ and $B\N_{\D,\sigma}$ by fiber-wise $\F_p$-completion of the
corresponding discrete versions $B\dN_{X,\sigma}$ and
$B\dN_{\D,\sigma}$.

\begin{recol}[The Adams-Mahmud map]\label{adamsmahmud}
By \cite[Lem.~4.1]{AGMV03} we have an ``Adams-Mahmud''
homomorphism $\AM: \Out(BX) \to \Out(\dN)$, given by associating to $f:BX \to BX$
the homomorphism $\AM(f): \dN \to \dN$, unique up to
conjugation, such that the diagram
$$
\xymatrix{
B\dN \ar[r]^-{B\AM(f)} \ar[d] & B\dN \ar[d] \\
BX \ar[r]^-f & BX
}
$$
commutes up to homotopy. 

By \cite[Thm.~C]{AG05automorphisms}, $\Phi$
factors through $\Out(\dN,\{\dN_\sigma\}) = \{ [\varphi] \in \Out(\dN)
\,|\, \varphi(\dN_\sigma) = \dN_{\varphi(\sigma)} \}$, which is isomorphic to
$\Out(\D_X)$ via restriction to $\dT$.
We again denote this map by $\Phi$.
Likewise, as e.g., explained in \cite[Prop.~5.1]{AGMV03}, fiber-wise
$\F_p$-completion \cite[Ch.~I, \S 8]{bk} induces a natural isomorphism
$\Out(\dN) \xrightarrow{\cong} \Out(B\N)$, and we can hence
equivalently view $\Phi(f)$ as an element in $\Out(B\N)$, and we will
not notationally distinguish between the two cases. We denote the
subgroup of $\Out(B\N)$ corresponding to $\Out(\dN,\{\dN_\sigma\})$ by
$\Out(B\N,\{B\N_\sigma\})$.
\end{recol}

\subsection{Centers and fundamental groups}

If $\D=(W, L, \{\Z_p b_\sigma\})$ is a $\Z_p$-root datum, we define a
{\em subdatum} of $\D$ to be a $\Z_p$-root datum of the form
$(W', L, \{\Z_p b_\sigma\}_{\sigma\in \Sigma'})$ where $(W', L)$ is a
reflection subgroup of $(W, L)$ and $\Sigma'$ is the set of
reflections in $W'$. For the next result, recall
\cite[Def.~4.1]{dw:center} that a homomorphism $f:BX \to BY$ is called
a {\em monomorphism of maximal rank} if the homotopy fiber $Y/X$ is
$\F_p$-finite and $X$ and $Y$ have the same rank.

\begin{prop} \label{subdata}
Let $X$ and $Y$ be connected $p$-compact groups. If $f: BX \to BY$ is
a monomorphism of maximal rank then $\D_X$ naturally identifies with a
subdatum of $\D_Y$.
\end{prop}

\begin{proof}
By definition there is a maximal torus $i: BT \to BX$ such that
$f\circ i:BT \to BY$ is a maximal torus for $Y$. Thus we can identify
$L_X = L_Y = \pi_1(T)$ and by \cite[Lem.~4.4]{dw:center} we have an
induced monomorphism $W_X \to W_Y$. This proves the result for $p$
odd, since in that case the $\Z_p$-root data $\D_X$ and $\D_Y$ are uniquely
determined by their underlying reflection groups $(W_X, L_X)$ and
$(W_Y, L_Y)$. In the case $p=2$ the result follows from the
construction of the $\Z_2$-root data of $X$ and $Y$, cf.\
\cite[Pf.~of~Lem.~9.16]{DW05}.
\end{proof}

For a $\Z_p$-root datum $\D = (W, L, \{\Z_p b_\sigma\})$,
we let $\dT = L \otimes_{\Z_p} \Z/p^\infty$ be the associated discrete
torus and for a reflection $\sigma\in W$, we define $h_\sigma =
b_\sigma/2\in \dT$. Clearly $h_\sigma$ is independent
of the choice of $b_\sigma$ and conversely $h_\sigma$ determines $\Z_p
b_\sigma$, cf.\ \cite[\S 2~and~\S 6]{DW05}. So instead of $(W, L,
\{\Z_p b_\sigma\})$ we might as well work with $(W, \dT,
\{h_\sigma\})$; we will use these two viewpoints interchangeably
without further mention. Also note that $h_\sigma = 1$ for $p$
odd. When $\sigma\in W$ is a reflection, we define the {\em singular
  set} $S(\sigma)$ by
$$
S(\sigma) = \left< \dT^+_0(\sigma), h_\sigma \right> =
\ker(\beta_\sigma \otimes_{\Z_p} \Z/p^\infty : \dT \to \Z/p^\infty),
$$
cf.\ \cite[Def.~7.3]{dw:center} and
\cite[(3.2)]{AG05automorphisms}. Define the {\em discrete
  center} $\dZ(\D)$ of $\D$ as $\bigcap_{\sigma} S(\sigma)$, where the
intersection is taken over all reflections $\sigma\in W$. In other
words, letting $M_0$ denote  the {\em root lattice}, i.e., the
$\Z_p$-sublattice of $L^*$ spanned by the roots $\beta_\sigma$,
we have the identification
\begin{equation}\label{singularset2}
\dZ(\D) = \ker \left(\dT = \Hom_{\Z_p}(L^*, \Z/p^\infty) \onto
  \Hom_{\Z_p}(M_0, \Z/p^\infty)\right).
\end{equation}
The following proposition translates into the present language the
results of Dwyer-Wilkerson \cite[\S 7]{dw:center} on how to compute
the center and centralizers of toral subgroups of $X$.

\begin{prop} \label{centercentralizer}
Let $X$ be a connected $p$-compact group with $\Z_p$-root datum $\D_X = (W,
L, \{\Z_p b_\sigma\})$.
\begin{enumerate}
\item \label{centerpart1}
The center $B\cZ(X)$ of $X$ is canonically homotopy
equivalent to the {\em center} $B\cZ(\D) = (B\dZ(\D))\pcom$ of $\D$.
\item \label{centerpart2}
The identity component $\cZ(X)_1$ of the center has $\Z_p$-root datum $(1,
L^W, \emptyset)$.
\item \label{centralizerpart}
For $A\subseteq \dT$, let $W(A)$ be the pointwise stabilizer of $A$ in
$W$, $\Sigma_A$ the set of reflections $\sigma$ with $A\subseteq
S(\sigma)$, and $W(A)_1$ the subgroup of $W$ generated by $\Sigma_A$.
Then the centralizer $\cC_X(A)$ has Weyl group $W(A)$ and its identity
component $\cC_X(A)_1$ has $\Z_p$-root datum equal to the
subdatum $\D_A = (W(A)_1,L,\{\Z_p b_\sigma\}_{\sigma\in \Sigma_A})$ of
$\D$.
\end{enumerate}
\end{prop}

\begin{proof}
The first part is \cite[Thm.~7.5]{dw:center}. The second part follows
easily from this since the maximal divisible subgroup of $\dZ(\D)$
equals
$$
\bigcap_{\sigma} \dT^+_0(\sigma) = \bigcap_{\sigma} \left(L^+(\sigma)
\otimes_{\Z_p} \Z/p^\infty\right) = L^W \otimes_{\Z_p} \Z/p^\infty.
$$
Part \eqref{centralizerpart} follows by combining
\cite[Thm.~7.6]{dw:center} and Proposition~\ref{subdata} since
$B\cC_X(A)_1 \to BX$ is a monomorphism of maximal rank
\cite[Prop.~4.3]{dw:center}.
\end{proof}

Let $\D=(W,L,\{\Z_p b_\sigma\})$ be a $\Z_p$-root datum. Recall that the
coroot lattice $L_0\subseteq L$ is the $\Z_p$-lattice spanned by the
coroots $b_\sigma$ and that the fundamental group $\pi_1(\D)$ is the
quotient $L/L_0$. Define $H^{\Z_p}_n(X) = \lim_k
H_n(X;\Z/p^k)$. The following proposition, which refines
\cite[Prop.~10.2]{AGMV03}, constructs a canonical epimorphism
$\pi_1(\D_X) \rightarrow \pi_1(X)$ which will be shown to be an
isomorphism in Theorem~\ref{fundgrp} below.

\begin{prop} \label{constrpi1map}
Let $X$ be a connected $p$-compact group with maximal torus $T$
and $\Z_p$-root datum $\D_X = (W, L, \{\Z_p b_\sigma\})$. Then the
homomorphism $L = H_2^{\Z_p}(BT) \to H_2^{\Z_p}(BX) \cong \pi_1(X)$
factors through $\pi_1(\D_X)$ and the induced homomorphism
$\pi_1(\D_X) \to \pi_1(X)$ is surjective with finite kernel.
\end{prop}

\begin{proof}
For $p$ odd we have $\im(1-\sigma) = \Z_p b_\sigma$ for all $\sigma$,
so $L/L_0 \cong H_0(W;L)$ and the result follows from
\cite[Prop.~10.2]{AGMV03}.

Assume now that $p=2$ and let $\sigma\in W$ be a reflection. To see
the first part we have to show that the homomorphism $L = H_2^{\Z_2}(BT)
\to H_2^{\Z_2}(BX) \cong \pi_1(X)$ vanishes on the coroots
$b_\sigma$. This follows from the construction of the root datum of
$X$:
$X(\sigma) = \cC_X(\dT_0^+(\sigma))$ is
a connected  $2$-compact group and
Proposition~\ref{centercentralizer}\eqref{centralizerpart} shows that
$\D_{X(\sigma)}$ equals the subdatum $(\left<\sigma\right>, L,
\{\Z_2 b_\sigma\})$ of $\D_X$. The commutative diagram
$$
\xymatrix{
L = H_2^{\Z_2}(BT) \ar[r] \ar@{=}[d] & H_2^{\Z_2}(BX(\sigma)) \cong
\pi_1(X(\sigma)) \ar[d] \\
L = H_2^{\Z_2}(BT) \ar[r] & H_2^{\Z_2}(BX) \cong \pi_1(X)
}
$$
shows that it suffices to prove the claim for $X(\sigma)$. However
since $X(\sigma)$ is a connected $2$-compact group of semi-simple rank $1$ it
follows (cf.\ \cite[p.~1369--1370]{DW05}) that $X(\sigma) \cong G\twocom$ for
$G = \SU(2)\times (S^1)^{r-1}$, $\SO(3) \times (S^1)^{r-1}$ or $\U(2) \times
(S^1)^{r-2}$ where $r$ is the rank of $X$. The
well-known formula for the fundamental group of a compact connected
Lie group in terms of its root datum, cf.\ \cite[\S 4,
no.\ 6, Prop.\ 11]{bo9} or \cite[Thm.~5.47]{adams69}, now established
the first part of the proposition.

Since $\im(1-\sigma) \subseteq \Z_2 b_\sigma$ by definition,
 $L \rightarrow L/L_0 = \pi_1(\D_X)$ factors through
$H_0(W;L)$, so the final claim now follows from \cite[Prop.~10.2]{AGMV03}.
\end{proof}

We will also be using the following formula for the fundamental group
of a $p$-compact group, proved by Dwyer-Wilkerson \cite{DW06} by a
transfer argument as this paper was being written (see also
\cite[Rem.~10.3]{AGMV03}). The formula was previously known for $p$
odd, by our classification \cite{AGMV03}, and we sketch in
Remark~\ref{pi1remark} how one can bypass the use of this formula also
in the classification for $p=2$ by a more cumbersome argument which we
had originally envisioned using in this paper; in particular providing
an independent proof.

\begin{thm}[{Dwyer-Wilkerson \cite{DW06}} and Remark~\ref{pi1remark}]
\label{fundgrp}
Let $X$ be a connected $p$-compact group. Then $\pi_1(\D_X) \xrightarrow{\cong}
\pi_1(X)$ induced by the maximal torus $T \to X$.
\end{thm}

\begin{proof}
By Theorem~\ref{rootdata-classification}\eqref{rdc-part1} we may write
$\D_X = \D_1 \times \D_2$, where $\D_1$ is of the form $\D'\otimes_{\Z}
\Z_p$ for a $\Z$-root datum $\D'$, and $\D_2$ is
a direct product of exotic $\Z_p$-root data. By
\cite[Thm.~1.4]{dw:split} this induces a splitting $BX \simeq BX_1 \times
BX_2$ with $\D_{X_i} \cong \D_i$. We have to show that the
kernel of $L=\pi_1(T) \to \pi_1(X)$ equals the coroot lattice
$L_0$; by the above it suffices to treat the case where $\D_X$ is
exotic and the case where $\D_X$ is of the form $\D'\otimes_{\Z} \Z_p$ for
a $\Z$-root datum $\D'$.

In the first case, Theorem~\ref{rootdata-classification}\eqref{rdc-part2} shows that
$\pi_1(\D_X)=0$ so the result follows from Proposition~\ref{constrpi1map}.

In the second case we have $\D_X \cong \D_{G\pcom}$ for some compact
connected Lie group $G$. By the result of Dwyer-Wilkerson
\cite[Thm.~1.1]{DW06} the kernel of $L = \pi_1(T) \to \pi_1(X)$ equals
the kernel of $H_2^{\Z_p}(BT_X) \to H_2^{\Z_p}(B\N_X)$. Since the
maximal torus normalizer may be reconstructed from the root datum by
\cite[Prop.~1.10]{DW05} for $p=2$ and \cite[Thm.~1.2]{kksa:thesis} for
$p$ odd, we may identify the homomorphism $H_2^{\Z_p}(BT_X) \to
H_2^{\Z_p}(B\N_X)$ with the homomorphism $H_2(BT_G;\Z) \to
H_2(BN_G;\Z)$ tensored by $\Z_p$. The result now follows from the
corresponding result for compact Lie groups, cf.\ \cite[\S 4, no.\ 6,
Prop.\ 11]{bo9} or \cite[Thm.~5.47]{adams69}.
\end{proof}

\subsection{Covers and quotients}

We now start to address how the root datum behaves upon taking covers
and quotients of a $p$-compact group.

\begin{lemma} \label{centerformula-maxrank}
Let $f: BX\to BY$ be a monomorphism of maximal rank between connected
$p$-compact groups $X$ and $Y$. If $\pi_1(\D_Y)=0$ then
$B\cZ(X) \simeq B\cZ(\N_X)$.
\end{lemma}

\begin{proof}
For $p$ odd, the conclusion holds for any connected $p$-compact group
$X$ \cite[Rem.~7.7]{dw:center}, so we may suppose $p=2$.
Since $\pi_1(\D_Y) = 0$, $\D_Y$ does not
have any direct factors isomorphic to $\D_{\SO(2n+1)\twocom}$, so
\cite[Lem.~5.1]{AG05automorphisms}
implies that the singular set $S_Y(\sigma)$ with respect to $Y$ equals
$\dT^+(\sigma)$ for any reflection $\sigma\in W_Y$. By
Proposition~\ref{subdata}, $\D_X$ identifies with a subdatum of
$\D_Y$ and hence $S_X(\sigma) = \dT^+(\sigma)$ for all reflections
$\sigma\in W_X$. The claim now follows from
Proposition~\ref{centercentralizer}\eqref{centerpart1}.
\end{proof}

\begin{lemma} \label{decomp}
Let $\D = (W, L, \{\Z_p b_\sigma\})$ be a $\Z_p$-root datum with
coroot lattice $L_0$. Then $L_0\cap L^W = 0$ and $L_0 \oplus L^W$ has
finite index in $L$. In particular $W$ acts faithfully on $L_0$.
\end{lemma}

\begin{proof}
The $\Q_p[W]$-module $V = L \otimes_{\Z_p} \Q_p$ decomposes as $V = V^W \oplus U$
where $U^W=0$. Writing $b_\sigma = x + y$ with $x\in V^W$ and $y\in U$
we have $\sigma(b_\sigma) = x + \sigma(y)$ and hence $\sigma(b_\sigma)
- b_\sigma \in U$. Also $\sigma(b_\sigma)\neq b_\sigma$: If $N=1+
\sigma + \ldots + \sigma^{|\sigma|-1}$ is the norm element, then
$\beta_\sigma(x) N(b_\sigma) = N (\sigma-1) (x) = 0$ for all $x\in
L$. Since $\sigma\neq 1$ we have $\beta_\sigma\neq 0$ so $N
b_\sigma=0$. Thus $\sigma(b_\sigma) \neq b_\sigma$ since otherwise $N
b_\sigma = |\sigma| b_\sigma\neq 0$. This proves that
$(\sigma-1)(b_\sigma) = r b_\sigma$ with $r\neq 0$ so $b_\sigma \in
U$. Thus $L_0\subseteq U$ and $L_0\cap L^W = 0$ as desired.  

Since $|W| x = \left(\sum_{w\in W} wx\right) + \sum_{w\in W}(x-wx) \in
L^W + L_0$, for any $x \in L$ we see that $L_0 \oplus L^W$ has finite
index in $L$, and in particular $W$ acts faithfully on $L_0$.
\end{proof}

Let $\D=(W, L, \{\Z_p b_\sigma\})$ be a $\Z_p$-root datum and let $L_0$
be the coroot lattice. If $L'$ is a $\Z_p$-lattice with $L_0\subseteq
L'\subseteq L$, the formula $\sigma(x) = x + \beta_\sigma(x) b_\sigma$
shows that $L'$ is $W$-invariant. By Lemma~\ref{decomp}, $W$ acts
faithfully on $L_0$ and hence also on $L'$, so $(W, L', \{\Z_p
b_\sigma\})$ is a $\Z_p$-root datum. We define a {\em cover} of $\D$
to be any $\Z_p$-root datum of this form. In particular the {\em
universal cover} $\widetilde{\D}$ of $\D$ is defined by
$\widetilde{\D} = (W, L_0, \{\Z_p b_\sigma\})$. Note that by definition,
$\pi_1(\widetilde{\D})=0$. 
For the reduction in Section~\ref{reductionsection} we need the
following result which does not rely on the fundamental group
formula Theorem~\ref{fundgrp}.

\begin{prop} \label{datum-of-topcover}
Let $X$ be a connected $p$-compact group with $\Z_p$-root datum $\D_X$ and let
$H$ be a subgroup of $\pi_1(X)$. Let $Y\to X$ be the cover of $X$
corresponding to $H$. Then $\D_Y$ is the cover of $\D_X$ corresponding
to the kernel of the composition $L_X \to \pi_1(\D_X) \to \pi_1(X) \to
\pi_1(X)/H$.
\end{prop}

\begin{proof}
By construction $BY$ is the fiber of $BX
\to B^2(\pi_1(X)/H)$. Let $BT \to BX$ be a maximal torus of $X$ and
let $B\N_X \to BX$ the maximal torus normalizer. Now consider the
following diagram obtained by pulling the fibration $BY \to BX \to
B^2(\pi_1(X)/H)$ back along $BT \to B\N_X \to BX$
\begin{equation}\label{firstdiag}
\xymatrix{
BT' \ar[r] \ar[d] & B\N' \ar[r] \ar[d] & BY \ar[d] \\
BT \ar[r] \ar[d] & B\N_X \ar[r] \ar[d] & BX \ar[d] \\
B^2(\pi_1(X)/H) \ar@{=}[r] & B^2(\pi_1(X)/H) \ar@{=}[r] & B^2(\pi_1(X)/H).
}
\end{equation}
Thus $BT' \to BY$ is a maximal torus and $B\N' \to BY$ is a maximal
torus normalizer by \cite[Thm.~1.2]{moller99}. The above diagram shows
that the Weyl group of $Y$ identifies with the Weyl group $W$ of
$X$. For a reflection $\sigma\in W$ we have the diagram
$$
\xymatrix{
B\dT' \ar[r] \ar[d] & B\dN'(\sigma) \ar[d] \ar[r] & BY(\sigma) \ar[d] \\
B\dT \ar[r] & B\dN_X(\sigma) \ar[r] & BX(\sigma),\\
}
$$
where $X(\sigma) = \cC_X(\dT_0^+(\sigma))$, $\dN_X(\sigma)
= C_{\dN_X}(\dT_0^+(\sigma))$ and similarly for $Y(\sigma)$
and $\dN'(\sigma)$. It now follows by definition (cf.\ \cite[\S
9]{DW05}) that the image of $h'_\sigma \in \dT'$ equals $h_\sigma \in
\dT$. Diagram \eqref{firstdiag} produces the short exact sequence $0 \to L_Y
\to L_X \to \pi_1(X)/H \to 0$ so $\Z_p b'_\sigma\subseteq L_Y$ maps to
$\Z_p b_\sigma\subseteq L_X$. This shows the claim.
\end{proof}

We next introduce quotients of root data. Let $\D=(W, \dT,
\{h_\sigma\})$ be a $\Z_p$-root datum, and $A\subseteq \dZ(\D)$
a subgroup of the discrete center. We define a {\em quotient} of $\D$
to be a $\Z_p$-root datum of the form $\D/A = (W, \dT/A,
\{\overline{h_\sigma}\})$ where $\overline{h_\sigma}$ denotes the
image of $h_\sigma$ in $\dT/A$; the fact that this is a $\Z_p$-root
datum is part of the following result.

\begin{prop} \label{datumofquotient}
\InsertTheoremBreak
\begin{enumerate}
\item \label{datumpart}
If $\D=(W, \dT, \{h_\sigma\})$ is a $\Z_p$-root datum and $A\subseteq
\dZ(\D)$, then $\D/A = (W, \dT/A, \{\overline{h_\sigma}\})$ is a
$\Z_p$-root datum. Moreover $\widetilde{\D/A} \cong \widetilde{\D}$
and $\dZ(\D/A) \cong \dZ(\D)/A$.
\item \label{grouppart}
Let $X$ be a connected $p$-compact group with $\Z_p$-root datum $\D_X$
and $A \to X$ a central monomorphism. If ${\breve A}$
denotes the discrete approximation to $A$, then the $\Z_p$-root datum
$\D_{X/A}$ of the $p$-compact group $X/A$ identifies with the quotient
datum $\D_X/{\breve A}$. In particular $\widetilde{\D_{X/A}}
\cong \widetilde{\D_X}$.
\end{enumerate}
\end{prop}

\begin{proof}
Write $\D = (W, L, \{\Z_p b_\sigma\})$ where $L =
\Hom_{\Z_p}(\Z/p^\infty,\dT)$ is
the associated $\Z_p$-lattice and $b_\sigma\in L$ the associated
coroots. 
By \eqref{singularset2} the sequence of
discrete tori $\dT \to \dT/A \to \dT/\dZ(\D)$ corresponds to the
sequence $L \to L_{\dT/A} \to {M_0}^*$ of $\Z_p$-lattices, where $M_0$
is the root lattice spanned by the $\beta_\sigma$. Note that $(W,L^*)$
is a reflection group via the action $\sigma(\alpha) = \alpha +
\alpha(b_{\sigma^{-1}})\beta_{\sigma^{-1}}$ so $(W, L^*, \{\Z_p
\beta_{\sigma^{-1}}\})$ is a $\Z_p$-root datum (the {\em dual} of
$\D$). Applying Lemma~\ref{decomp}  to this $\Z_p$-root datum shows that
$M_0 \oplus (L^*)^W$ has finite index in $L^*$, so the dual
homomorphism $L \to {M_0}^* \oplus \left(\left(L^*\right)^W\right)^*$
is injective. The last summand is isomorphic to
$\Hom_{\Z_p}(H_0(W;L),\Z_p)^*$ which again identifies with $H_0(W;L)$
modulo its torsion subgroup. The formula $\sigma(x) - x =
\beta_{\sigma}(x) b_\sigma$ shows that the image of the composition
$L_0 \to L \to H_0(W;L)$ is torsion, so by the above the composition
$L_0 \to L \to L_{\dT/A} \to {M_0}^*$ is injective. In particular $W$
acts faithfully on $L_{\dT/A}$ by Lemma~\ref{decomp}. The homomorphism
$\beta_\sigma \otimes_{\Z_p} \Z/p^\infty : \dT \to \Z/p^\infty$
factors through $\dT/A$, and we get a corresponding homomorphism
$\beta'_\sigma: L_{\dT/A} \to \Z_p$ such that the composition $L \to
L_{\dT/A} \xrightarrow{\beta'_\sigma} \Z_p$ agrees with
$\beta_\sigma$. We now claim that $\sigma(x) = x + \beta'_\sigma(x)
b_\sigma$ for $x\in L_{\dT/A}$. To see this note that the image of $L$
in $L_{\dT/A}$ has finite index since we have the exact sequence $L
\to L_{\dT/A} \to \text{Ext}^1_{\Z_p}(\Z/p^\infty,A)$ with
$\text{Ext}^1_{\Z_p}(\Z/p^\infty,A)$ finite. Thus the claim follows
from the above by using the corresponding formula $\sigma(x) = x +
\beta_\sigma(x) b_\sigma$ for $x\in L$. This proves that $\D/A = (W,
L_{\dT/A}, \{\Z_p b_\sigma\})$ is a $\Z_p$-root datum. In particular
we obtain $\widetilde{\D/A} = (W, L_0, \{\Z_p b_\sigma\}) =
\widetilde{\D}$ by definition. 

To see the claim about $\dZ(\D/A)$,
note that by the above $\beta_\sigma \otimes_{\Z_p} \Z/p^\infty : \dT
\to \Z/p^\infty$ identifies with the composition $\dT \to \dT/A
\xrightarrow{\beta'_\sigma \otimes_{\Z_p} \Z/p^\infty}
\Z/p^\infty$. Hence the singular sets for $\D$ and $\D/A$ satisfies
$S_{\D}(\sigma)/A = S_{\D/A}(\sigma)$ and the claim follows.

To see part~\eqref{grouppart}, 
let $i: BT \to BX$ be a maximal torus and let $f: BA\to BX$ be the central
monomorphism. Then $f$ factors through $BT$ by \cite[Prop.~8.11]{DW94}
to give a central monomorphism $g: BA \to BT$. Moreover $i$ factors
through the maximal torus normalizer $\N$ of $X$ and we obtain the
diagram
$$
\xymatrix{
BT \ar[r] \ar[d] & B\N \ar[r] \ar[d] & BX \ar[d] \\
BT/A \ar[r] & B\N/A \ar[r] & BX/A
}
$$
cf.\ Construction~\ref{adjointconstruction}. It follows that $T/A$ is a
maximal torus in $X/A$ and $\N/A$ is the maximal torus normalizer. The
Weyl groups of $X$ and $X/A$ identifies naturally,
cf.\ \cite[Thm.~4.6]{MN94}. By construction the
elements $h'_\sigma \in \dT/{\breve A}$ corresponding to $\D_{X/A}$
are the images of the elements $h_\sigma \in \dT$ corresponding to
$X$  (cf.\ the proof of
Proposition~\ref{datum-of-topcover}). This shows that $\D_{X/A} \cong
\D_X/{\breve A}$ as desired.
\end{proof}

As a special case of the quotient construction, we define the {\em
  adjoint} $\D_{\text{ad}}$ of a $\Z_p$-root datum $\D = (W, L, \{\Z_p
b_\sigma\})$ by $\D_{\text{ad}} = \D/\dZ(\D)$. Note that by
Proposition~\ref{datumofquotient}\eqref{datumpart} we have
$\dZ(\D_{\text{ad}})=0$ and that it follows
from the proof that $\D_{\text{ad}} = (W, {M_0}^*, \{\Z_p
b_\sigma\})$, where $M_0$ is the {\em root lattice}, i.e., the
sublattice of $L^*$ spanned by the roots $\beta_\sigma$.

\begin{prop} \label{datasplit}
Any $\Z_p$-root datum with $\dZ(\D) = 0$ or $\pi_1(\D) = 0$ splits as
a direct product $\D \cong \D_1 \times \ldots \times \D_n$ of irreducible
$\Z_p$-root data $\D_i$.
\end{prop}

\begin{proof}
The case $\dZ(\D) = 0$ is essentially proved by Dwyer-Wilkerson
\cite[Pf.~of~Thm.~1.5]{dw:split}, for completeness we briefly sketch
the argument: Let $(W,L)$ be the $\Z_p$-reflection group associated to
$\D$ and write $L\otimes_{\Z_p} \Q_p = V_1 \oplus \ldots \oplus V_n$ be the
decomposition of $L\otimes_{\Z_p} \Q_p$ into irreducible
$\Q_p[W]$-modules. Define $L_i=L\cap V_i$. Since $\dZ(\D) = 0$ we have
$\bigcap_{\sigma} \dT^+_0(\sigma) = 0$ as well, and it follows
\cite[Pf.~of~Thm.~1.5]{dw:split} that the homomorphism $L_1 \times
\ldots \times L_n \to L$ is an isomorphism. Letting $W_i$ denote the
pointwise stabilizer of $L_1 \oplus \ldots \oplus \widehat{L_i}
\oplus \ldots \oplus L_n$ we hence (cf.\ \cite[Prop.~7.1]{dw:split})
get a product decomposition $(W,L) \cong (W_1,L_1) \times \ldots
\times (W_n,L_n)$. It is now clear that there is a unique $\Z_p$-root
datum structure $\D_i$ on $(W_i,L_i)$ such that we get a product
decomposition $\D \cong \D_1 \times \ldots \times \D_n$ into
irreducible $\Z_p$-root data.

The case where $\pi_1(\D) = 0$ is easily reduced to the first case
using the previous results: By
Proposition~\ref{datumofquotient}\eqref{datumpart},
$\dZ(\D_{\text{ad}}) = 0$ so we can write $\D_{\text{ad}} \cong \D_1 \times
\ldots \times \D_n$ where the $\D_i$ are
irreducible. Proposition~\ref{datumofquotient}\eqref{datumpart} now
shows that $\D = \widetilde{\D} \cong \widetilde{\D_{\text{ad}}} \cong
\widetilde{\D_1} \times \ldots \times \widetilde{\D_n}$ as claimed.
\end{proof}

\begin{thm}[The classification of $\Z_p$-root data; structure version]
  \label{data-quotient-structure}
\InsertTheoremBreak
\begin{enumerate}
\item \label{dqs-part1}
Let $\D = (W, L, \{\Z_p b_\sigma\})$ be a $\Z_p$-root datum with
coroot lattice $L_0$, and let $\D' = (W, L_0 \oplus L^W, \{\Z_p
b_\sigma\}) = \widetilde{\D} \times \D_{\text{triv}}$, where
$\D_{\text{triv}} = (1, L^W, \emptyset)$ is a trivial $\Z_p$-root
datum. Then $\D \cong \D'/A$ for a finite central subgroup $A\subseteq
\dZ(\D')$ and there is a splitting $\widetilde{\D} \cong \D_1 \times
\ldots \times \D_n$ of $\widetilde{\D}$ into irreducible $\Z_p$-root
data $\D_i$ with $\pi_1(\D_i) = 0$.
\item \label{dqs-part2}
For $p>2$, the assignment $\D = (W, L, \{\Z_p b_\sigma\})
\rightsquigarrow (W, L\otimes_{\Z_p} \Q_p)$ is a one-to-one
correspondence between isomorphism classes of irreducible $\Z_p$-root
data $\D$ with $\pi_1(\D)=0$ and isomorphism classes of irreducible
$\Q_p$-reflection groups. For $p=2$, the assignment is surjective and
the preimage of every element consists of a single element, except for
$(W_{\Sp(n)},L_{\Sp(n)}\otimes \Q_2) \cong
(W_{\Spin(2n+1)},L_{\Spin(2n+1)}\otimes \Q_2)$ whose preimage consists
of $\D_{\Sp(n)}\otimes_\Z \Z_2$ and $\D_{\Spin(2n+1)}\otimes_\Z \Z_2$
which are non-isomorphic for $n \geq 3$.
\end{enumerate}
\end{thm}

\begin{proof}
The fact that $\D'$ is a $\Z_p$-root datum follows from
Lemma~\ref{decomp}. The short exact sequence $0 \to L_0 \oplus L^W \to
L \to F \to 0$ where $F$ is finite with trivial $W$-action produces a
short exact sequence $1 \to A \to \dT' \to \dT \to 1$ between the
associated discrete tori. Since the roots $\beta'_\sigma$ for $\D'$
are given by $L_0 \oplus L^W \to L \xrightarrow{\beta_\sigma} \Z_p$,
it follows that $A\subseteq S_{\D'}(\sigma) = \ker((L_0\oplus L^W)
\otimes_{\Z_p} \Z/p^\infty \to \Z/p^\infty)$ for any reflection
$\sigma$. Hence $A\subseteq \dZ(\D')$ is central and $\D \cong
\D'/A$. The last part of \eqref{dqs-part1} follows from Proposition~\ref{datasplit}.

We now prove \eqref{dqs-part2}. For any prime $p$,
Theorem~\ref{rootdata-classification}\eqref{rdc-part2} shows that the
assignment in \eqref{dqs-part2} gives a one-to-one correspondence
between isomorphism classes of exotic $\Z_p$-root data and
isomorphism classes of exotic $\Q_p$-reflection groups, and that
$\pi_1(\D) = 0$ for all exotic $\Z_p$-root data $\D$. Hence it
suffices by Theorem~\ref{rootdata-classification}\eqref{rdc-part1} to
show the claim for $\Z_p$-root data of the form $\D_1 \otimes_{\Z}
\Z_p$ for a $\Z$-root datum $\D_1$.

In this case it is clear that if $\pi_1(\D_1 \otimes_{\Z} \Z_p) = 0$,
then we can find a $\Z$-root datum $\D'_1$ with $\pi_1(\D'_1)=0$ and
$\D'_1 \otimes_{\Z} \Z_p \cong \D_1 \otimes_{\Z} \Z_p$ (simply choose
$\D'_1$ to be the universal cover of $\D_1$; this is defined for
$\Z$-root data in the same way as for $\Z_p$-root data). Hence it
suffices to study the assignment $\D = (W, L, \{\Z_p b_\sigma\})
\rightsquigarrow (W, L\otimes_{\Z} \Q)$ from irreducible $\Z$-root
data with $\pi_1(\D)=0$ to irreducible $\Q$-reflection groups. It is
well-known (cf.\ \cite[\S 4]{bo9})
that this assignment is surjective and that it only fails
to be injective in that the $\Z$-root data $\D_{\Sp(n)}$ and
$\D_{\Spin(2n+1)}$ which are non-isomorphic for $n\geq 3$ maps to the
same $\Q$-reflection group.
This proves part \eqref{dqs-part2} since
for $n\geq 3$, the $\Z_p$-root data $\D_{\Sp(n)}\otimes_\Z \Z_p$ and
$\D_{\Spin(2n+1)}\otimes_\Z \Z_p$ are non-isomorphic for $p=2$ and
isomorphic for $p=2$.
\end{proof}

\subsection{Automorphisms}

Recall that an {\em isomorphism} between two $\Z_p$-root data $\D =
(W, L, \{\Z_p b_\sigma \})$ and $\D' = (W', L', \{\Z_p b'_{\sigma'}
\})$ is an isomorphism $\varphi : L \to L'$ with the property that $\varphi
W \varphi^{-1} = W'$ as subgroups of $\Aut(L')$ and $\varphi(\Z_p
b_\sigma) = \Z_p b'_{\varphi\sigma\varphi^{-1}}$ for every reflection
$\sigma\in W$. We denote the automorphism group of $\D$ by $\Aut(\D)$;
clearly $W$ is a normal subgroup of $\Aut(\D)$, and we define the
{\em outer automorphism group} $\Out(\D)$ by $\Out(\D) =
\Aut(\D)/W$. 

Recall that a $\Z_p$-root datum $\D = (W, L, \{\Z_p b_\sigma\})$ is
called irreducible if $L\otimes_{\Z_p} \Q_p$ is an irreducible
$\Q_p[W]$-module. The following proposition is a restatement of
\cite[Prop.~5.4]{AGMV03}.

\begin{prop} \label{rootdatumauto}
Suppose $\D_i = (W_i,L_i,\{\Z_p b_\sigma\}_{\sigma\in \Sigma_i})$,
$i=0, \ldots, k$, is a collection of pairwise non-isomorphic irreducible
$\Z_p$-root data. Assume that $W_0 = 1$ but that $W_i$ is non-trivial
for $i \geq 1$. Let $\D = \prod_{i=0}^k \D_i^{m_i}$, $m_i\geq 1$
denote a product of these $\Z_p$-root data. Then
$$
\GL_{m_0}(\Z_p) \times
\left( \prod_{i=1}^k \Out(\D_i)
\wr \Sigma_{m_i}
\right)
\xrightarrow{\cong}
\Out(\D).
$$
\end{prop}

\begin{proof}
This follows directly by combining \cite[Prop.~5.4]{AGMV03} with
\cite[Rem.~4.5]{AG05automorphisms}.
\end{proof}

The following two results are needed for the reduction in
Section~\ref{reductionsection}.

\begin{prop} \label{datum-structure}
Let $\D = (W, L, \{\Z_p b_\sigma\})$ be a $\Z_p$-root datum with
coroot lattice $L_0$. Let $\D' = (W, L_0 \oplus L^W, \{\Z_p
b_\sigma\}) = \widetilde{\D} \times \D_{\text{triv}}$, where
$\D_{\text{triv}} = (1, L^W, \emptyset)$ is a trivial $\Z_p$-root datum. Then
$\D \cong \D'/A$ for a finite subgroup $A\subseteq \dZ(\D')$
and the restriction
$$
\Aut(\D) \longrightarrow \Aut(\D') = \Aut(\widetilde{\D}) \times
\Aut(\D_{\text{triv}})
$$
is an isomorphism onto the subgroup $\{ \varphi\in \Aut(\D')\,\mid\,
\varphi(A) = A \}$. In particular $\Out(\D)$ identifies with a
subgroup of finite index in $\Out(\D')$.
\end{prop}

\begin{rem}
For a $\Z_p$-root datum $\D = (W, L, \{\Z_p b_\sigma\})$ with
$\pi_1(\D)=0$ we have $\Z_p b_\sigma = \ker(L \xrightarrow{N} L)$,
where $N = 1+\sigma + \ldots + \sigma^{|\sigma|-1}$ is the norm
element: By Theorem~\ref{rootdata-classification}\eqref{rdc-part1}
either $\D \cong \D_1 \otimes_\Z \Z_p$ for a $\Z$-root datum $\D_1$
with $\pi_1(\D_1)=0$ (cf.\ the proof of
Theorem~\ref{data-quotient-structure}\eqref{dqs-part2}) or $\D$ is
exotic. In the first case the result is well known \cite[\S 4]{bo9},
and in the second case the claim follows since
$H^1(\left<\sigma\right>;L)=0$ for all reflections $\sigma$ (cf.\ the
proof of Theorem~\ref{rootdata-classification}\eqref{rdc-part2}).
In particular we obtain $\Out(\D) = N_{\GL(L)}(W)/W$, and this
group is explicitly computed for all irreducible $(W,L)$ in
\cite[Thm.~13.1]{AGMV03}.
\end{rem}

\begin{proof}[Proof of Proposition~\ref{datum-structure}]
The fact that $\D \cong \D'/A$ for a finite subgroup $A\subseteq
\dZ(\D')$ is part of Theorem~\ref{data-quotient-structure}\eqref{dqs-part1}, and the
identification $\Aut(\D') = \Aut(\widetilde{\D}) \times
\Aut(\D_{\text{triv}})$ comes from
Proposition~\ref{rootdatumauto}. Clearly both $L_0$ and $L^W$ are
invariant under $\Aut(\D)$ so we get a restriction homomorphism
$\Aut(\D) \longrightarrow \Aut(\D')$ which is injective since $L_0
\oplus L^W$ has finite index in $L$ by Lemma~\ref{decomp}.
If $\varphi \in \Aut(\D)$, then $\varphi$
corresponds to an automorphism $\varphi : \dT \to \dT$ and the
restriction of $\varphi$ to $L_0 \oplus L^W$ corresponds to a lift
$\widetilde{\varphi} : \dT' \to \dT'$ which sends $A$ to
itself. Conversely, if $\varphi \in \Aut(\D')$ with $\varphi(A) = A$,
then $\varphi$ clearly induces an automorphism of $\D$.
The last claim follows from the fact that the orbit of $A$ under
$\Out(\D')$ is finite since $\dZ(\D')$ has only finitely many
subgroups isomorphic to $A$.
\end{proof}

\begin{cor} \label{datum-adjisocover}
For any $\Z_p$-root datum $\D$ there is a canonical isomorphism
$$
\Aut(\widetilde{\D}) \stackrel{\cong}{\longleftrightarrow}
\Aut(\D_{\text{ad}}).
$$
\end{cor}

\begin{proof}
Write $\D = (W, L, \{\Z_p b_\sigma\})$, and let $M_0 \subseteq L^*$
denote the root lattice, i.e., the lattice spanned by the roots
$\beta_\sigma$. Then $\widetilde{\D} = (W, L_0, \{\Z_p b_\sigma\})$
and the roots for $\widetilde{\D}$ are given by the composition $L_0
\to L \xrightarrow{\beta_\sigma} \Z_p$. Hence we can identify the root
lattice for $\widetilde{\D}$ with $M_0$ and hence
$\widetilde{\D}/\dZ(\widetilde{\D}) =
(\widetilde{\D})_{\text{ad}} \cong (W, {M_0}^*, \{\Z_p
b_\sigma\}) = \D_{\text{ad}}$. The result now follows from
Proposition~\ref{datum-structure}.
\end{proof}

\subsection{Finiteness properties}

In this final subsection we prove that there are only finitely many
$\Z_p$-root data of a given rank, and that for a fixed $\Z_p$-root
datum $\D$, $\Out(\D)$ only contains finitely many finite subgroups up
to conjugation. These results are used for the proof of the finiteness
statements in Theorems~\ref{cohomologicaluniqueness} and
\ref{nonconnclassification} in Section~\ref{consequences-section}.

\begin{prop}\label{finiteoffixedrank}
For any prime $p$ there is, up to isomorphism, only finitely many
$\Z_p$-root data of a fixed rank.
\end{prop}

\begin{proof}
By the classification of finite $\Q_p$-reflection groups \cite{CE74}
\cite{DMW92}, there are only finitely many finite $\Q_p$-reflection
groups of a fixed rank. For each finite $\Q_p$-reflection group
$(W,V)$, there are, up to equivalence, only finitely many
representations $W \to \GL(V)$ which gives rise to the reflection
group $(W,V)$. For each of these representations, the local version of the
Jordan-Zassenhaus theorem \cite[Thm.~24.7]{CR81} shows that, up to
isomorphism, there are only finitely many $\Z_p[W]$-lattices $L$ with
$L\otimes_{\Z_p} \Q_p \cong V$ as $\Q_p[W]$-modules. In particular we
conclude that, up to isomorphism, there are only finitely many finite
$\Z_p$-reflection groups of fixed rank. Finally, choosing a
$\Z_p$-root datum for a finite $\Z_p$-reflection group $(W,L)$
corresponds to choosing a cyclic subgroup of the finite group
$H^1(\left<\sigma\right>;L)$ for each conjugacy class of reflections
$\sigma$. Hence any finite $\Z_p$-reflection group gives rise to only
finitely many $\Z_p$-root data. This proves the result.
\end{proof}

\begin{prop}\label{replemma}
Let $\D$ be a $\Z_p$-root datum. Then $\Out(\D)$ contains only
finitely many conjugacy classes of finite subgroups.
\end{prop}

\begin{proof}
Let $(W,L)$ be the finite $\Z_p$-reflection group underlying $\D$, and
let $n=\rank L$. Note first that the order of a finite subgroup of
$\GL(L)$ is bounded above: If $G\subseteq \GL_n(\Z_p)$ has finite
order, then it is easily seen that the composition $G \to \GL_n(\Z_p)
\to \GL_n(\F_p)$ is injective for $p>2$ and has kernel of order at
most $2^{n^2}$ for $p=2$ (cf.\ \cite[Lem.~11.3]{AGMV03}). Hence $G$
has order at most $\left|\GL_n(\F_p)\right|$ for $p>2$ and $2^{n^2}\cdot
\left|\GL_n(\F_2)\right|$ for $p=2$. In particular there is an upper
bound on the order of finite subgroups of $N_{\GL(L)}(W)$. Since
$\Out(\D)$ is contained in $N_{\GL(L)}(W)/W$ and $W$ is finite, it
follows that there is also an upper bound on the order of finite
subgroups of $\Out(\D)$.

Fix a finite group $G$. By the above it suffices to show that the set
$\Rep(G,\Out(\D))$ (i.e., the set of homomorphisms $G\to \Out(\D)$
modulo conjugation in $\Out(\D)$) is finite. By
Theorem~\ref{data-quotient-structure}\eqref{dqs-part1}, we can write $\D \cong \D'/A$
where $\D' = \widetilde{\D} \times \D_{\text{triv}}$ and
$\D_{\text{triv}}$ is a trivial $\Z_p$-root datum, and this identifies
$\Out(\D)$ with a subgroup of finite index in $\Out(\D')$, cf.\
Proposition~\ref{datum-structure}. Hence it is enough to prove that
$\Rep(G,\Out(\D'))$ is finite since this will imply that
$\Rep(G,\Out(\D))$ is finite.

By Theorem~\ref{data-quotient-structure}\eqref{dqs-part1} and
Proposition~\ref{rootdatumauto}, $\Out(\D')$ is isomorphic
to a direct product of $\GL_{m_0}(\Z_p)$ and groups of the
form $\Out(\D_i)\wr \Sigma_{m_i}$ where the $\D_i$ are
irreducible $\Z_p$-root data. It thus suffices to know that
$\Rep(G,\GL_m(\Z_p))$ and $\Rep(G,\Out(\D)\wr \Sigma_m)$ are finite
for any $m$ and any irreducible $\Z_p$-root datum $\D$. The first
claim follows directly from the local version of the Jordan-Zassenhaus
theorem \cite[Thm.~24.7]{CR81}.

To see the second claim, let $W$ denote the Weyl group of $\D$ and
note that since $\D$ is irreducible, Schur's lemma implies that the
image of the central homomorphism $\Z_p^\times \to \Out(\D)$ equals
the kernel of the canonical homomorphism $\Out(\D) \to \Out(W)$. Since
$\Z_p^\times \cong \Z_p \times C$ where $C$ is finite ($C=\Z/2$ for $p=2$ and
$C=\Z/(p-1)$ for $p>2$) and $\Out(W)$ is finite, it follows that
$\Out(\D)$ is a finite (central) extension of $\Z_p$. In particular
$E=\Out(\D) \wr \Sigma_m$ fits into an extension of the form
\begin{equation} \label{extE}
1 \to (\Z_p)^m \to E \xrightarrow{\pi} Q \to 1,
\end{equation}
where $Q$ is finite. Any homomorphism $\varphi:G \to E$ gives a homomorphism
$\pi\circ\varphi:G \to Q$ by composition, and for a fixed homomorphism
$\alpha:G \to Q$, the set of homomorphisms $\varphi:G \to E$ with
$\pi\circ\varphi = \alpha$ equals the set of splittings of the
pullback $1 \to (\Z_p)^m \to E' \to G \to 1$ of \eqref{extE} along
$\alpha:G \to Q$. The set of such splittings, modulo conjugation by
elements in $(\Z_p)^m$, is in one-to-one correspondence with
$H^1(G;(\Z_p)^m)$, which is finite for all actions of $G$ on
$(\Z_p)^m$. Hence the set of homomorphisms $\varphi: G\to E$ with
$\pi\circ\varphi = \alpha$ is finite modulo conjugation in $(\Z_p)^m
\subseteq E$. Since $Q$ is finite, there are only finitely many
homomorphisms $\alpha:G \to Q$, so we conclude that $\Rep(G,E)$ is
finite as claimed. 
\end{proof}


\bibliographystyle{plain}
\bibliography{../pcg/poddclassification}
\end{document}